\documentclass{amsart}

\usepackage{latexsym,amsmath,amstext,mathrsfs,amsthm,amscd,amsopn,verbatim,amsrefs}
\usepackage[english,french]{babel}

\usepackage{graphicx}
\usepackage[dvips]{color}
\usepackage{epsfig}
\usepackage{psfrag} 
\usepackage[all]{xy}
\usepackage{amscd}
\usepackage{amssymb}

\newtheorem{Thm}{Theorem}[section]
\newtheorem{Prop}[Thm]{Proposition}
\newtheorem{Lem}[Thm]{Lemma}

\theoremstyle{remark}
\newtheorem{Rem}[Thm]{Remark}

\theoremstyle{definition}

\begin{document}

\title[Biquantization of symmetric pairs and the quantum shift]{Biquantization of symmetric pairs and the quantum shift}
\author{A.~S.~Cattaneo}
\email{alberto.cattaneo@math.uzh.ch}
\address{Institut f\"ur Mathematik, Universität Z\"urich, Winterthurerstrasse 190
CH-8057 Z\"urich (Switzerland)}
\author{C.~A. Rossi}
\email{crossi@mpim-bonn.mpg.de}
\address{MPIM Bonn, Vivatsgasse 7, 53111 Bonn (Germany)}
\author{C.~Torossian}
\email{torossian@math.dot.jussieu.dot.fr}
\address{Universit\'e Denis-Diderot-Paris 7, UFR de math\'matiques, Site Chevaleret, Case 7012 , 75205 Paris cedex 13 (France)}

%\classification{16E40; 16E45; 81R60}
%\keywords{Deformation quantization; symmetric pairs}
\thanks{The first author acknowledges partial support by SNF Grant $200020\_131813/1$}

\maketitle

\selectlanguage{french}
\begin{abstract}
Dans~\cite{CT} la biquantification des paires sym\'etriques a \'et\'e \'etudi\'e \`a l'aide du formalisme graphique de M.~Kontsevich.
Dans ce papier on corrige, compte tenu des r\'esultats montr\'es r\'ecemment dans~\cite{CFFR}, une erreur mineure dans~\cite{CT}, qui en tout cas n'invalide pas les r\'esultats plus importants dans~\cite{CT}: cette erreur consiste au fait que les auteurs avaient oubli\'e une contribution provenante d'une certaine composante du bord dans le propagateur \`a quatre couleurs.
La correction qu'on apporte ici a l'avantage de remettre finalement en jeu la ``translation quantique'' des charact\`eres qui appara\^it dans la m\'ethode des orbites, et qui \'etait myst\`erieusement absente dans~\cite{CT}.
En plus, on pr\'esente une comparaison d\'etaill\'ee des deux fa\c{c}ons diff\'erentes de construire la biquantification, {\em i.e.} en utilisant ou le d\'emi-plan de Poincar\'e ou le premier quadrant, ainsi qu'un approche plus conceptuel \`a la biquantification selon~\cite{CFFR} et toutes les corrections dues des r\'esultats dans~\cite{CT} qu'il faut corriger \`a cause de la pr\'esence de la translation quantique.
Finalement on reconsid\`ere la construction de la triquantification d\'ev\'elopp\'ee dans la partie finale de~\cite{CT} pour l'\'etude des charact\`eres compte tenu du m\^eme probl\`eme du bord dans la biquantification.  
\end{abstract}

\selectlanguage{english}
\begin{abstract}
The biquantization of symmetric pairs was studied in~\cite{CT} in terms of Kontsevich-like
graphs. 
This paper, also in view of recent results in~\cite{CFFR}, amends a minor mistake that did
not spoil the main results of the paper. 
The mistake consisted in ignoring a regular term in the boundary contribution of some propagators. 
On the other hand, its correction brings back the quantum shift, present in the approaches by the
orbit method, that was otherwise puzzlingly missing. 
In addition a detailed comparison of the two, equivalent, ways of defining biquantization working on the upper half plane or on one
quadrant is presented, as well as a more conceptual approach to biquantization and the due corrections of some results of~\cite{CT} in view of the aforementioned correction by the quantum shift.
Finally, we review the triquantization construction developed for the treatment of characters by taking into accounts the same boundary problem as for the biquantization.
\end{abstract}

\selectlanguage{english}

\tableofcontents

\section{Introduction}\label{s-0}
In \cite{CT} the biquantization of symmetric pairs was studied in terms of Kontsevich-like graphs. 
A puzzling result was the absence of the quantum shift, otherwise present in the treatments using the orbit method,
by the natural character of the adjoint representation of $\mathfrak k$ on $\mathfrak g/\mathfrak k=\mathfrak p$, 
where $\mathfrak g=\mathfrak k\oplus\mathfrak p$ is the symmetric pair under consideration.
It turns out that due to a (fortunately minor) mistake in \cite{CT} the quantum shift is actually there. Apart from this, the mistake does not spoil 
the other results of the paper.

The whole construction of \cite[Section 1.6]{CT}  relies on the 4-colored propagators introduced in \cite{CFb} for the Poisson sigma model with two branes. 
It was recently observed by
G.~Felder and the second author in the preparation of \cite{CFFR}, that, unlike in Kontsevich \cite{K}, the boundary contributions of the $4$-colored propagators on the first quadrant for the
collapse of the two endpoints may have a regular term in addition to the usual singular one. 
The regular term turns out simply to be the differential of the angle of the
position where the two points collapsed, measured with respect to the origin, up to a sign, which depends on the boundary conditions themselves (roughly speaking, if we consider the same boundary conditions on the two half-lines bounding the first quadrant, then the sign is positive, while, for different boundary conditions on the two half-lines, we have a negative sign).
Recall that these propagators are constructed from the Euclidean propagator
(the differential of the angle of the line joining the two points) by reflecting the second argument with respect to the two boundaries of the first quadrant (producing four distinct closed $1$-forms on the compactified configuration space of two points in the interior of the first quadrant) and then summing them up with signs; concretely,
\[
\omega^{\varepsilon_1,\varepsilon_2}(z_1,z_2)=\frac{1}{2\pi}\left[\mathrm d\ \mathrm{arg}(z_2-z_1)+\varepsilon_1\mathrm d\ \mathrm{arg}(z_2-\overline z_1)+\varepsilon_2 \mathrm d\ \mathrm{arg}(z_2+\overline z_1)+\varepsilon_1\varepsilon_2 \mathrm d\ \mathrm{arg}(z_2+z_1)\right], 
\]
where $\varepsilon_i=\pm$, $i=1,2$.

The contribution where the second argument is reflected w.r.t.\ the origin (corresponding to the situation where the second argument is reflected w.r.t.\ both boundaries of the first quadrant) is responsible for the regular term in all four situations, see Figure 1.
\bigskip
\begin{center}
\resizebox{0.43 \textwidth}{!}{\input{reg_term.pstex_t}}\\
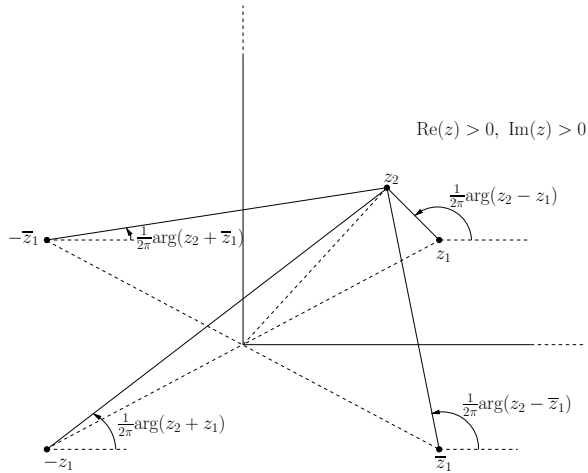
\text{Figure 1 - A geometric explanation of the ``singular term'' and the ``regular term''} \\
\end{center}
\bigskip
The presence of this regular term was mistakenly neglected in \cite{CT}.
Its net effect is that more boundary contributions have to be taken into account and extra terms are needed for cancellation. It turns out \cite{CFFR} that it is
enough to allow for the presence of short loops and to assign each of them the regular term. This also has the pleasant effect of restoring the quantum shift.
Some by-products, in particular in~\cite[Sections 4 and 5]{CT}, have to be modified accordingly.

Regular terms also appear in the 8-colored propagators introduced in \cite[Section 6]{CT} for the three brane case that is needed to show the independence from the choice of polarization. 
Also in this case the introduction of short loops, consistent with what was observed above, saves the game: we will review these aspects, as well as their relationship with the Harish-Chandra homomorphism.

\section{Notation and conventions}\label{s-1}
We work over a ground field $\mathbb K$, which may be $\mathbb R$ or $\mathbb C$.
We consider a finite-dimensional symmetric pair $\mathfrak g$ over $\mathbb K$, {\em i.e.} a Lie algebra $\mathfrak g$, endowed with a Lie algebra automorphism $\sigma$, which is additionally an involution. 
In particular, $\mathfrak g=\mathfrak k\oplus\mathfrak p$, where $\mathfrak k$, resp.\ $\mathfrak p$, is the eigenspace w.r.t.\ the eigenvalue $+1$, resp.\ $-1$, of $\sigma$.
For a Lie subalgebra $\mathfrak h$ of a Lie algebra $\mathfrak g$, we denote by $\mathfrak h^\perp$ its annihilator.

As $\mathfrak k$ is a Lie subalgebra of $\mathfrak g$ and $\mathfrak g/\mathfrak k=\mathfrak p$ is a $\mathfrak k$-module, we introduce the short-hand notation 
\[
\delta(\bullet)=\frac{1}2\mathrm{tr}_{\mathfrak p}(\mathrm{ad}_\mathfrak k(\bullet)),
\]
see {\em e.g.} also~\cite{T1,T2}.

%We refer to~\cite[Section 5]{CFFR} for the construction of the $2$-colored and $4$-colored propagators in the construction of the relative formality theorem and of the $2$-brane formality theorem respectively.

\section{Biquantization in the framework of the $2$-brane formality}\label{s-2}
We consider a finite-dimensional Lie algebra $\mathfrak g$ over $\mathbb K$; further, we consider two Lie subalgebras $\mathfrak h_i$, $i=1,2$.

To these data, we associate a Poisson manifold $X$ and two coisotropic submanifolds $U_i$, $i=1,2$, as follows: we set $X=\mathfrak g^*$, endowed with the linear Kirillov--Kostant--Souriau Poisson bivector $\pi$, and $U_i=\mathfrak h_i^\perp$.
%\begin{Rem}\label{r-KKS}
%In order to be consistent with the convention in~\cite{CT}, the Poisson bivector field on $\mathfrak g^*$ corresponds to half the Lie bracket on $\mathfrak g$.
%\end{Rem}

We want to regard biquantization as analyzed in~\cite{CT} in the more general framework of the $2$-brane Formality Theorem~\cite[Theorem 7.2]{CFFR}: thus, before entering into the details of biquantization, we need to review in some detail the main result of~\cite{CFFR} and draw a bridge between it and the computations in~\cite{CT}.

\subsection{On compactified configuration spaces}\label{ss-2-1}
For the upcoming discussion of the $4$-colored propagators, we need to fix certain issues regarding compactified configuration spaces: in particular, we discuss two types of compactified configuration spaces, which arise in the context of biquantization, and we prove that they are in fact diffeomorphic.
We observe that the following discussion may be viewed as an extension of certain computations in~\cite{CFFR}.

\subsubsection{The compactified configuration space of points in $\mathbb H^+\sqcup \mathbb R$}\label{sss-2-1-1}
For two non-negative integers $m$, $n$, we consider the (open) configuration space $C_{n,m}^+$ of $n$ distinct points in the complex upper half-plane $\mathbb H^+$ and $m$ ordered points on the real axis $\mathbb R$.
Its precise definition is
\[
C_{n,m}^+=\left\{(z_1,\dots,z_n,x_1,\dots,x_m)\in (\mathbb H^+)^n\times \mathbb R^m:\ z_i\neq z_j,\ i\neq j,\ x_1<\cdots<x_m\right\}/G_2,
\]
where $G_2$ is the two-dimensional real Lie group $\mathbb R^+\ltimes \mathbb R$, acting on $\mathbb H^+\sqcup \mathbb R$ by rescalings and real translations. 
Provided $2n+m-2\geq 0$, $C_{n,m}^+$ is a smooth manifold of dimension $2n+m-2$; it is obviously oriented.

We further consider the open configuration space
\[
C_n=\left\{(z_1,\dots,z_n)\in \mathbb C^n:\ z_i\neq z_j,\ i\neq j\right\}/G_3,
\]
where $G_3$ is the $3$-dimensional real Lie group $\mathbb R^+\ltimes \mathbb C$, acting on $\mathbb C$ by rescalings and complex translations.
It is obvious that, provided $2n-3\geq 0$, $C_n$ is a smooth manifold of dimension $2n-3$, which admits an obvious orientation from $\mathbb C^n$ and from the obvious orientability of $G_3$.

Kontsevich~\cite[Subsection 5.1]{K} provides compactifications $\mathcal C_{n,m}^+$ and $\mathcal C_n$ of $C_{n,m}^+$ and $C_n$ respectively in the sense of Fulton--MacPherson~\cite{FMcP} (to be more precise, the smooth version of the algebraic compactification of~\cite{FMcP}, exploited in detail by Axelrod--Singer~\cite{AS}): both compactified configuration spaces admit structures of smooth manifolds with corners ({\em i.e.} locally modeled on $(\mathbb R^+)^k\times\mathbb R^l$), and as such they admit boundary stratifications. 

We observe that the permutation group $\mathfrak S_n$ acts naturally on $C_n$, and it can be proved that its action extends to $\mathcal C_n$: in particular, we may consider more general compactified configuration spaces $\mathcal C_A$, for a finite subset of $\mathbb N$.
Because of similar reasons, we may consider compactified configuration spaces $\mathcal C_{A,B}^+$, for a finite subset $A$ and a finite, ordered subset $B$ of $\mathbb N$.

The stratifications of $\mathcal C_{n,m}^+$ and $\mathcal C_n$ admit a beautiful description in terms of trees~\cite[Subsection 5.1]{K}.
We first consider the boundary stratification of $\mathcal C_n$: for simplicity, we illustrate the boundary strata of codimension $1$, namely such boundary strata are labeled by subsets $A$ of $[n]=\{1,\dots,n\}$ of cardinality $2\leq |A|\leq n$,
\[
\partial_A \mathcal C_n\cong \mathcal C_A\times\mathcal C_{([n]\smallsetminus A)\sqcup \{\bullet\}}, 
\]
where the first, resp.\ second, factor on the right-hand side of the previous identification represents the configuration of distinct points in $\mathbb C$ labeled by $A$ which collapse together in $\mathbb C$ to a single point $\bullet$, resp.\ the final configuration of points after the collapse.
%The symbol ``$\sqcup$'' means disjoint union.

The boundary strata of codimension $1$ of $\mathcal C_{n,m}^+$ are of two types, namely,
\begin{itemize}
\item[$i)$] there exists a subset $A$ of $[n]$, of cardinality $2\leq |A|\leq n$, such that 
\[
\partial_A \mathcal C_{n,m}^+\cong \mathcal C_A\times \mathcal C_{([n]\smallsetminus A)\sqcup \{\bullet\},m}^+,
\]
where the first, resp.\ second, factor on the right-hand side of the previous identification describes the collapse of the points in $\mathbb H^+$ labeled by $A$ to a single point $\bullet$ in $\mathbb H^+$, resp.\ the final configuration of points after the collapse;
\item[$ii)$] there exist a subset $A$ of $[n]$ and an ordered subset of $[m]$ consisting of consecutive non-negative integers, such that $0\leq |A|\leq n$, $0\leq |B|\leq m$, $1\leq |A|+|B|\leq n+m-1$, for which we have
\[
\partial_{A,B}\mathcal C_{n,m}^+\cong \mathcal C_{A,B}^+\times \mathcal C_{[n]\smallsetminus A,([m]\smallsetminus B)\sqcup \{\bullet\}}^+,
\]
where the first, resp.\ second, factor on the right-hand side of the previous identification describes the collapse of the points in $\mathbb H^+$ labeled by $A$ and the consecutive, ordered points on $\mathbb R$ labeled by $B$ to a single point $\bullet$ in $\mathbb R$, resp.\ the final configuration of points after the collapse.
\end{itemize}

\subsubsection{The compactified configuration space of points in $Q^{+,+}\sqcup i\mathbb R^+\sqcup \mathbb R^+$}\label{sss-2-1-2}
For three non-negative integers $l$, $m$ and $n$, we consider the (open) configuration space $C_{n,k,l}^+$ of $n$ distinct points in the first quadrant $Q^{+,+}$, $k$ ordered points on the positive imaginary axis $i\mathbb R^+$ and $l$ ordered points on the positive real axis $\mathbb R$.
We observe that the order of the points on the positive imaginary axis is the opposite of the intuitive one, {\em i.e.} $i x\leq i y$ if and only if $y\leq x$, for $x$, $y$ in $\mathbb R$.

The precise definition of $C_{n,k,l}^+$ is
\[
\begin{aligned}
C_{n,k,l}^+=&\left\{(z_1,\dots,z_n,ix_1,\dots,ix_k,y_1,\dots,y_l)\in (Q^{+,+})^n\times (i\mathbb R^+)^k\times (\mathbb R^+)^l:\ z_i\neq z_j,\ i\neq j,\right.\\
&\left.\ x_k<\cdots<x_1,\ y_1<\cdots< y_l\right\}/G_1,
\end{aligned}
\]
where $G_1=\mathbb R^+$ acts on $Q^{+,+}\sqcup i\mathbb R^+\sqcup \mathbb R^+$ by rescalings. 
Provided $2n+k+l-1\geq 0$, $C_{n,k,l}^+$ is a smooth manifold of dimension $2n+k+l-1$.
It inherits an obvious orientation from the natural orientation of $(Q^{+,+})^n\times (i\mathbb R^+)^k\times(\mathbb R^+)^l$ and the one of $\mathbb R^+$.

We may provide a compactification $\mathcal C_{n,k,l}^+$ in the sense of Fulton--MacPherson~\cite{FMcP} of $C_{n,k,l}^+$ in a way similar to Kontsevich: the compactified configuration space $\mathcal C_{n,k,l}^+$ admits a structure of smooth manifold with corners. 

We now consider the boundary strata of $\mathcal C_{n,k,l}^+$ of codimension $1$, which are of the three following types:
\begin{itemize}
\item[$i)$] there exists a subset $A$ of $[n]$, of cardinality $2\leq |A|\leq n$, such that 
\[
\partial_A \mathcal C_{n,k,l}^+\cong \mathcal C_A\times \mathcal C_{([n]\smallsetminus A)\sqcup \{\bullet\},k,l}^+,
\]
where the first, resp.\ second, factor on the right-hand side of the previous identification describes the collapse of the points in $Q^{+,+}$ labeled by $A$ to a single point $\bullet$ in $Q^{+,+}$, resp.\ the final configuration of points after the collapse;
\item[$ii)$] there exist a subset $A$ of $[n]$ and an ordered subset of $[k]$, resp.\ $[l]$, consisting of consecutive non-negative integers, such that $0\leq |A|\leq n$, $0\leq |B|\leq k$, resp.\ $0\leq |B|\leq l$, for which we have
\[
\partial_{A,B}\mathcal C_{n,k,l}^+\cong \mathcal C_{A,B}^+\times \mathcal C_{[n]\smallsetminus A,([k]\smallsetminus B)\sqcup \{\bullet\},l}^+,\ \text{resp.}\ \partial_{A,B}\mathcal C_{n,k,l}^+\cong \mathcal C_{A,B}^+\times \mathcal C_{[n]\smallsetminus A,k,([l]\smallsetminus B)\sqcup \{\bullet\}}^+ 
\]
where the first, resp.\ second, factor on the right-hand side of the previous identification describes the collapse of the points in $Q^{+,+}$ labeled by $A$ and the consecutive, ordered points on $i\mathbb R^+$ or $\mathbb R^+$ labeled by $B$ to a single point $\bullet$ in $i\mathbb R^+$ or $\mathbb R^+$, resp.\ the final configuration of points after the collapse.
\item[$iii)$] there exist a subset $A$ of $[n]$ and an ordered subset $B=B_1\sqcup B_2$ of $[k]\sqcup [l]$, for which $B_1$ and $B_2$ are ordered subsets of consecutive points in $[k]$ and $[l]$ respectively, such that $k\in B_1$ if $B_1\neq \emptyset$, $1\in B_2$ if $B_2\neq \emptyset$, $0\leq |A|\leq n$, $0\leq |B|\leq k+l$, resp.\ $1\leq |A|+|B|\leq n+k+l-1$, for which we have
\[
\partial_{A,B_1,B_2}\mathcal C_{n,k,l}^+\cong \mathcal C_{A,B_1,B_2}^+\times \mathcal C_{[n]\smallsetminus A,[k]\smallsetminus (B\cap[k]),[l]\smallsetminus (B\cap [l])}^+, 
\]
where the first, resp.\ second, factor on the right-hand side of the previous identification describes the collapse of the points in $Q^{+,+}$ labeled by $A$ and the consecutive, ordered points on $i\mathbb R^+$ and $\mathbb R^+$ labeled by $B=B_1\sqcup B_2$ to the origin of the axes, resp.\ the final configuration after the collapse.
\end{itemize}

\subsubsection{The relationship between $\mathcal C_{n,m}^+$ and $\mathcal C_{n,k,l}^+$}\label{sss-2-1-3}
First of all, we consider the open configuration spaces $C_{n,m}^+$ and $C_{n,k,l}^+$, where $m=k+l+1$.

We observe that the holomorphic function $z\mapsto z^2$ on $\mathbb C$, when restricted to $Q^{+,+}\sqcup i\mathbb R^+\sqcup \mathbb R^+$, gives rise to a biholomorphism to $\mathbb H^+\sqcup (\mathbb R\smallsetminus \{0\})$, whose inverse we denote by $z\mapsto \sqrt z$: in fact, we have to choose a well-suited branch-cut for the complex square root, {\em e.g.} we cut out from the complex plane the negative imaginary axis plus the origin.

For $m$, $k$ and $l$ as before, we choose the $k+1$-st point on $\mathbb R$.
Then, there is an obvious map from $C_{n,m}^+$ to $C_{n,k,l}^+$, which is defined by the following explicit formula:
\begin{equation}\label{eq-square}
\begin{aligned}
&C_{n,m}^+\ni [(z_1,\dots,z_n,x_1,\dots,x_{k+1},\dots,x_m)]\mapsto\\
&\mapsto\left[\left(\sqrt{z_1-x_{k+1}},\dots,\sqrt{z_n-x_{k+1}},\sqrt{x_1-x_{k+1}},\dots,\sqrt{x_k-x_{k+1}},\sqrt{x_{k+1}-x_{k+1}},\dots,\sqrt{x_m-x_{k+1}}\right)\right]\in C_{n,k,l}^+.
\end{aligned}
\end{equation}
First of all, we observe that, because of the order on the points on the real axis, the difference $x_i-x_{k+1}$ is strictly negative, resp.\ positive, if $1\leq i\leq k$, resp.\ $k+2\leq i\leq m$: thus, because of the said choice of a complex square root, $\sqrt{x_i-x_{k+1}}=i\sqrt{x_{k+1}-x_i}$, if $1\leq i\leq k$, or $\sqrt{x_i-x_{k+1}}=\sqrt{x_i-x_{k+1}}$, if $k+2\leq i\leq m$, where now both square roots on the right-hand side of both equalities are real, positive numbers.
Again, the order on the points $x_i$, $1\leq i\leq k$ implies that $\sqrt{x_{k+1}-x_i}>\sqrt{x_{k+1}-x_{i+1}}$, therefore, the natural order on $x_i$, $1\leq i\leq k$, is mapped to the natural order on $i\sqrt{x_{k+1}-x_i}$ discussed in Subsubsection~\ref{sss-2-1-2}.
We may depict the morphism~\eqref{eq-square} graphically {\em via}
\bigskip
\begin{center}
\resizebox{0.6 \textwidth}{!}{\input{sq-comp.pstex_t}}\\
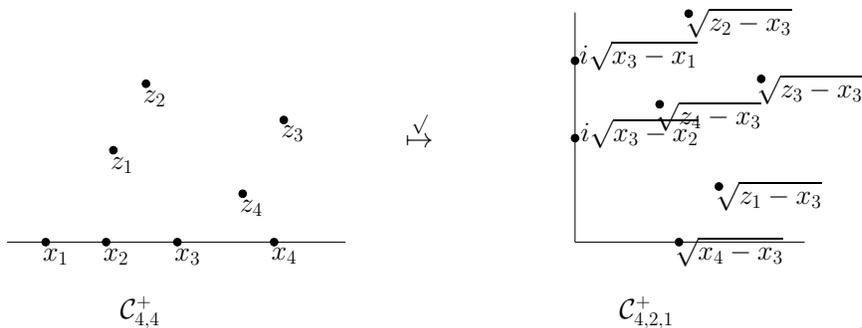
\text{Figure 2 - A pictorial description of the action of the complex square root on $\mathcal C_{4,4}^+$} \\
\end{center}
\bigskip
An easy computation proves that the above morphism is well-defined, {\em i.e.} it does not depend on the choice of representatives; furthermore, the morphism is obviously smooth, and is in fact a diffeomorphism, whose inverse is 
\[
C_{n,k,l}^+\ni [(z_1,\dots,z_n,ix_1,\dots,i x_k,y_1,\dots,y_l)]\mapsto\left[\left(z_1^2,\dots,z_n^2,-x_1^2,\dots,-x_k^2,0,y_1^2,\dots,y_l^2\right)\right]\in C_{n,m}^+.
\] 

The important point is that the complex square function and the chosen inverse (the above complex square root) extend to smooth functions between the compactified configuration spaces $\mathcal C_{n,m}^+$ and $\mathcal C_{n,k,l}^+$.
\begin{Prop}\label{p-square}
For non-negative integers $n$, $m$, $k$, $l$, such that $m=k+l+1$, the smooth manifolds with corners $\mathcal C_{n,m}^+$ and $\mathcal C_{n,k,l}^+$ are diffeomorphic {\em via} the choice of a complex square root with branch cut $i\mathbb R^-\sqcup \{0\}$.
\end{Prop}
\begin{proof}
We prove that the diffeomorphism~\eqref{eq-square} extends to a diffeomorphism on the compactified configuration spaces by computing its expression w.r.t.\ local coordinates for the relevant boundary strata of codimension $1$.
In fact, as sketched in~\cite[Subsection 5.2]{K}, the boundary strata of higher codimension correspond to products with more than two factors of compactified configuration spaces of the same kind, representing configuration of points collapsing together, be it in the complex upper half-plane or on the real axis, resp.\ in the first quadrant or on the positive complex or real axis or on the origin.

It suffices therefore to prove the claim on the interior of the boundary strata of codimension $1$ of $\mathcal C_{n,m}^+$ and $\mathcal C_{n,k,l}^+$: these have been characterized explicitly in Subsubsections~\ref{sss-2-1-1} and~\ref{sss-2-1-2}.
Furthermore, without loss of generality, we may assume $A=[i]$ and $B=[j]$.

We have to prove that the map~\eqref{eq-square} maps diffeomorphically the interior of boundary strata of codimension $1$ of $\mathcal C_{n,m}^+$ to the interior of boundary strata of codimension $1$ of $\mathcal C_{n,k,l}^+$.  

We consider first the boundary stratum of type $i)$ of $\mathcal C_{n,m}^+$ labeled by $A=[i]$, for $2\leq i\leq n$.
Local coordinates of the interior $C_i\times \mathcal C_{n-i+1,m}^+$ are provided by 
\[
C_i\times C_{n-i+1,m}^+\ni \left(\left(e^{i\varphi},z_1,\dots,z_{i-2}\right),\left(e^{it},w_1,\dots,w_{n-i},x_1,\dots,x_k,0,x_{k+2},\dots,x_m\right)\right),
\]
where $\varphi$ is in $[0,2\pi)$, $t$ in $(0,\pi)$, $z_i$ in $\mathbb C$, $w_i$ in $\mathbb H^+$, and all points in $\mathbb C$ and $\mathbb H^+$ are distinct, while the points on the real axis are lexicographically (strictly) ordered.
On the other hand, the interior of the boundary stratum $\mathcal C_i\times\mathcal C_{n-i+1,k,l}^+$ is described by the following local coordinates:
\[
C_i\times C_{n-i+1,k,l}^+\ni \left(\left(e^{i\varphi},z_1,\dots,z_{i-2}\right),\left(e^{it},w_1,\dots,w_{n-i},i x_1,\dots,i x_k,y_1,\dots,y_l\right)\right),
\]
where $\varphi$ is in $[0,2\pi)$, $t$ in $(0,\frac{\pi}2)$, $z_i$ in $\mathbb C$, $w_i$ in $Q^{+,+}$, and all points in $\mathbb C$ and $Q^{+,+}$ are distinct, and $x_1>\cdots x_k>0$, $0<y_1<\cdots<y_l$.

For $\varepsilon>0$ sufficiently small, local coordinates for $\mathcal C_{n,m}^+$, resp.\ $\mathcal C_{n,k,l}^+$, near the interior of the boundary stratum $\mathcal C_i\times\mathcal C_{n-i,m}^+$, resp.\ $\mathcal C_i\times\mathcal C_{n-i+1,k,l}^+$, are given by
\[
\begin{aligned}
&\left[\left(e^{it},e^{it}+\varepsilon e^{i\varphi},e^{it}+\varepsilon z_1,\dots,e^{it}+\varepsilon z_{i-2},w_1,\dots,w_n,x_1,\dots,x_k,0,x_{k+2},\dots,x_m\right)\right],\ \text{resp.}\\
&\left[\left(e^{it},e^{it}+\varepsilon e^{i\varphi},e^{it}+\varepsilon z_1,\dots,e^{it}+\varepsilon z_{i-2},w_1,\dots,w_n,ix_1,\dots,i x_k,y_1,\dots,y_l\right)\right].
\end{aligned}
\]
We apply the morphism~\eqref{eq-square} to the first of the previous expressions, getting 
\begin{equation}\label{eq-sq-1}
\left[\left(\sqrt{e^{it}},\sqrt{e^{it}+\varepsilon e^{i\varphi}},\sqrt{e^{it}+\varepsilon z_1},\dots,\sqrt{e^{it}+\varepsilon z_{i-2}},\sqrt{w_1},\dots,\sqrt{w_n},i\sqrt{-x_1},\dots,i\sqrt{-x_k},\sqrt{x_{k+2}},\dots,\sqrt{x_m}\right)\right].
\end{equation}
We rewrite the terms in the previous expression containing the infinitesimal parameter $\varepsilon$ using the fact that the chosen complex square root is holomorphic on $\mathbb H^+$, thus getting
\[
\sqrt{e^{it}+\varepsilon z}=e^{i\frac{t}2}+\frac{\varepsilon z}{2 e^{i\frac{t}2}}+\mathcal O(\varepsilon^2)=e^{i\frac{t}2}+\frac{\varepsilon}2 e^{-i\frac{t}2}z+\mathcal O(\varepsilon^2),\ t\in (0,\pi),\ z\in \mathbb C.
\]
To compare expressions, we may neglect terms in the expansion of order strictly higher than $2$: rescaling by $\frac{1}2$ the infinitesimal parameter $\varepsilon$, Expression~\eqref{eq-sq-1} can be rewritten as 
\[
\left[\left(e^{i\frac{t}2},e^{i\frac{t}2}+\varepsilon e^{i\left(\varphi-\frac{t}2\right)},e^{i\frac{t}2}+\varepsilon e^{-i\frac{t}2}z_1,\dots,e^{i\frac{t}2}+\varepsilon e^{-i\frac{t}2}z_{i-2},\sqrt{w_1},\dots,\sqrt{w_n},i\sqrt{-x_1},\dots,i\sqrt{-x_k},\sqrt{x_{k+2}},\dots,\sqrt{x_m}\right)\right],
\]
whence it follows immediately that the morphism~\eqref{eq-square} maps the interior of $\mathcal C_i\times\mathcal C_{n-i+1,m}^+$ diffeomorphically to the interior of $\mathcal C_i\times \mathcal C_{n-i+1,k,l}^+$, where the diffeomorphism is explicitly the product of the morphism~\eqref{eq-square} from $C_{n-i+1,m}^+$ to $C_{n-i+1,k,l}^+$ with the obvious diffeomorphism of $C_i$ given by 
\[
C_i\ni \left[\left(z_1,\dots,z_i\right)\right]\mapsto \left[\left(\frac{z_1}{\sqrt{w_1-x_{k+1}}},\dots,\frac{z_i}{\sqrt{w_1-x_{k+1}}}\right)\right]\in C_i,
\]
where $w_1$ and $x_{k+1}$ are taken from $C_{n-i+1,m}^+$.

We now consider the interior of the boundary stratum $\mathcal C_{i,B}^+\times \mathcal C_{n-i,([m]\smallsetminus B)\sqcup \{\bullet\}}^+$, where $B$ is an ordered subset of $[m]$ consisting of consecutive elements, and we assume that $1\leq i+|B|\leq n+m-1$.
We have to further distinguish between two situations: $|B|=0$ (and consequently $1\leq i\leq n$), and $|B|\neq 0$.

We consider the situation $|B|=0$, and we further distinguish between the case, where the new point $\bullet$ on the real axis (corresponding to the cluster of points labeled by $[i]$ in $\mathbb H^+$ approach $\mathbb R$) lies on the left or on the right of the distinguished point $x_{k+1}$.
We do the explicit computations only in the case, where $\bullet$ is on the left of $x_{k+1}$, leaving the other case to the reader.

If $\bullet$ lies on the left of $x_{k+1}$, we may safely assume that $\bullet=x$ lies on the left of $x_1$: then, local coordinates for the interior of $\mathcal C_{i,0}^+\times \mathcal C_{n-i,m+1}^+$ are given by
\[
C_{i,0}^+\times C_{n-i,m+1}^+\ni \left(\left(i,z_1,\dots,z_{i-1}\right),\left(e^{it},w_1,\dots,w_{n-i-1},x,x_1,\dots,x_k,0,x_{k+2},\dots,x_m\right)\right),
\] 
where $t$ in $(0,\pi)$, $z_i$ and $w_j$ are in $\mathbb H^+$, and all points in $\mathbb H^+$ are distinct, while the points on the real axis are lexicographically (strictly) ordered.
Similarly, local coordinates for the interior of $\mathcal C_{i,0}^+\times \mathcal C_{n-i,m+1}^+$ are given by
\[
C_{i,0}^+\times C_{n-i,k+1,l}^+\ni \left(\left(i,z_1,\dots,z_{i-1}\right),\left(e^{it},w_1,\dots,w_{n-i-1},ix,ix_1,\dots,ix_k,y_1,\dots,y_l\right)\right),
\] 
where $t$ in $(0,\frac{\pi}2)$, $z_i$ in $\mathbb H^+$ and $w_i$ in $Q^{+,+}$, all points in $\mathbb H^+$ and $Q^{+,+}$ are distinct, $x>x_1>\cdots>x_k>0$ and $0<y_1<\cdots<y_l$.

For $\varepsilon>0$ sufficiently small, local coordinates for $\mathcal C_{n,m}^+$, resp.\ $\mathcal C_{n,k,l}^+$, near the interior of the boundary stratum $\mathcal C_{i,0}^+\times\mathcal C_{n-i,m+1}^+$, resp.\ $\mathcal C_{i,0}^+\times\mathcal C_{n-i,k+1,l}^+$, are given by
\[
\begin{aligned}
&\left[\left(x+\varepsilon i,x+\varepsilon z_1,\dots,x+\varepsilon z_{i-1},e^{it},w_1,\dots,w_{n-i-1},x_1,\dots,x_k,0,x_{k+2},\dots,x_m\right)\right],\ \text{resp.}\\
&\left[\left(ix+\varepsilon,ix-i\varepsilon z_1,\dots,ix-i\varepsilon z_{i-1},e^{it},w_1,\dots,w_{n-i-1},ix_1,\dots,ix_k,y_1,\dots,y_l\right)\right]
\end{aligned}
\]
The image of the first expression w.r.t.\ the morphism~\eqref{eq-square} is simply
\begin{equation}\label{eq-sq-2}
\left[\left(\sqrt{x+\varepsilon i},\sqrt{x+\varepsilon z_1},\dots,\sqrt{x+\varepsilon z_{i-1}},\sqrt{e^{it}},\sqrt{w_1},\dots,\sqrt{w_{n-i-1}},i\sqrt{-x_1},\dots,i\sqrt{-x_k},\sqrt{x_{k+2}},\dots,\sqrt{x_m}\right)\right].
\end{equation}
We consider the first $i$ entries in the previous expression: once again, using the holomorphy of the chosen complex square root, and recalling that $x<x_1<0$ and that $\varepsilon$ is chosen sufficiently small, we find
\[
\begin{aligned}
\sqrt{x+\varepsilon i}&=\sqrt{x}+\frac{\varepsilon i}{2\sqrt x}+\mathcal O(\varepsilon^2)=i\sqrt{-x}+\frac{\varepsilon}{2\sqrt{-x}}+\mathcal O(\varepsilon^2),\\  
\sqrt{x+\varepsilon z_j}&=i\sqrt{-x}-\frac{i\varepsilon}2 \frac{z_j}{\sqrt{-x}}+\mathcal O(\varepsilon^2),\ 1\leq j\leq i-1. 
\end{aligned}
\]
Once again, rescaling by $\frac{1}2$ the infinitesimal parameter $\varepsilon$, and neglecting terms of order higher than $1$ w.r.t.\ $\varepsilon$ in the above expressions, we may rewrite Expression~\eqref{eq-sq-2} as
\[
\left[\left(i\sqrt{-x}+\varepsilon,i\sqrt{-x}+\varepsilon \frac{z_1}{\sqrt{-x}},\dots,i\sqrt{-x}+\varepsilon\frac{z_{i-1}}{\sqrt{-x}},e^{i\frac{t}2},\sqrt{w_1},\dots,\sqrt{w_{n-i-1}},i\sqrt{-x_1},\dots,i\sqrt{-x_k},\sqrt{x_{k+2}},\dots,\sqrt{x_m}\right)\right],
\]
and it is easy to see that the morphism~\eqref{eq-square} maps $C_{i,0}^+\times C_{n-i,m+1}^+$ diffeomorphically to $C_{i,0}^+\times C_{n-i,k+1,l}^+$, and the induced morphism is precisely given by the product of the morphism~\eqref{eq-square} from $C_{n-i,m+1}^+$ to $C_{n-i,k+1,l}^+$ with the obvious diffeomorphism of $C_{i,0}^+$ given by
\[
C_{i,0}^+\ni \left[\left(z_1,\dots,z_i\right)\right]\mapsto \left[\left(\frac{z_1}{\sqrt{x_{k+1}-x}},\dots,\frac{z_i}{\sqrt{x_{k+1}-x}}\right)\right]\in C_{i,0}^+,
\]
where $x$ denotes the first point on the real axis in lexicographical order, and $x_{k+1}$ is the special point on real axis.

For the situation $|B|\neq 0$, we need to distinguish between two cases, namely $i)$ $B$ contains $k+1$ or $ii)$ $b$ does not contain $k+1$ (in which case, either the minimum of $B$ is greater or equal than $k+2$ or the maximum of $B$ is less or equal than $k$).

We first consider the case, where $B$ contains $k+1$, and we assume $A=[i]$ and $B=[p,q]=\{p,\dots,q\}$, where $1\leq p\leq k+1\leq q\leq m$; we further write $B=\{k\}\sqcup B_1\sqcup B_2$, where $B_1=[p,k]$ and $B_2=[k+2,q]$ (of course, $B_1$ and/or $B_2$ may be empty).
The interior of the corresponding boundary stratum of $C_{n,m}^+$, resp.\ $C_{n,k,l}^+$, is $C_{i,B}^+\times C_{n-i,([m]\smallsetminus B)\sqcup \{\bullet\}}^+$, resp.\ $C_{i,B_1,B_2}^+\times C_{n-i,[k]\smallsetminus B_1,[l]\smallsetminus B_2}^+$, and corresponding local coordinates are given by 
\[
\begin{aligned}
C_{i,B}^+\times C_{n-i,([m]\smallsetminus B)\sqcup \{\bullet\}}^+\ni &\left(\left(e^{it_1},z_1,\dots,z_{i-1},x_1,\dots,x_{k-p+1},0,x_{k-p+3},\dots,x_{q-p+1}\right),\right.\\
&\left.\left(e^{it_2},w_1,\dots,w_{n-i-1},x_1',\dots,x_{p-1}',0,x_{p+1}',\dots,x_{m-2k-q}'\right)\right),\\
C_{i,B_1,B_2}^+\times C_{n-i,[k]\smallsetminus B_1,[l]\smallsetminus B_2}^+\ni &\left(\left(e^{it_1},z_1,\dots,z_{i-1},i x_1,\dots,i x_{k-p+1},y_1,\dots,y_{q-k-1}\right),\right.\\
&\left.\left(e^{it_2},w_1,\dots,w_{n-i-1},i x_1',\dots,i x_{p-1}',y_1',\dots,y_{l-q+k+1}'\right)\right),
\end{aligned}
\]
where $t_i$, $i=1,2$, is in $(0,\pi)$, resp.\ $(0,\frac{\pi}2)$, all points in $\mathbb H^+$, resp.\ $Q^{+,+}$, are distinct in the first, resp.\ second, expression.
In the first, resp.\ second, expression, the $x_i$ and $x_i'$ are lexicographically ordered, resp.\ $x_1>\cdots x_{k-p+1}>0$, $x_1'>\cdots x_{p-1}'>0$, $0<y_1<\cdots y_{q-k-1}$ and $0<y_1'<\cdots< y'_{l-q+k+1}$.

Choosing a positive number $\varepsilon$ sufficiently small as before, we may write local coordinates of $\mathcal C_{n,m}^+$, resp.\ $\mathcal C_{n,k,l}^+$, near the interior of the boundary stratum $\mathcal C_{i,B}^+\times\mathcal C_{n-i,([m]\smallsetminus B)\sqcup \{\bullet\}}^+$, resp.\ $\mathcal C_{i,B_1,B_2}^+\times \mathcal C_{n-i,[k]\smallsetminus B_1,[l]\smallsetminus B_2}^+$, namely
\[
\begin{aligned}
&\left[\left(\varepsilon e^{i t_1},\varepsilon z_1,\dots,\varepsilon z_{i-1},e^{it_2},w_1,\dots,w_{n-i-1},\right.\right.\\
&\left.\left.x_1',\dots,x_{p-1}',\varepsilon x_1,\dots,\varepsilon x_{k-p+1},0,\varepsilon x_{k-p+3},\dots,\varepsilon x_{q-p+1},x_{p+1}',\dots,x_{m-2k-q}'\right)\right],\ \text{resp.}\\
&\left[\left(\varepsilon e^{i t_1},\varepsilon z_1,\dots,\varepsilon z_{i-1},e^{i t_2},w_1,\dots,w_{n-i-1},\right.\right.\\
&\left.\left.ix_1',\dots,ix_{p-1}',\varepsilon ix_1,\dots,\varepsilon ix_{k-p+1},0,\varepsilon y_1,\dots,\varepsilon y_{q-k+1},y_1',\dots,y'_{l-q+k+1}\right)\right].
\end{aligned}
\]
We now apply the morphism~\eqref{eq-square} to the first of the two previous expressions, getting
\[
\begin{aligned}
&\left[\left(\sqrt{\varepsilon e^{i t_1}},\sqrt{\varepsilon z_1},\dots,\sqrt{\varepsilon z_{i-1}},\sqrt{e^{it_2}},\sqrt{w_1},\dots,\sqrt{w_{n-i-1}},\right.\right.\\
&\left.\left.i\sqrt{-x_1'},\dots,i\sqrt{-x_{p-1}'},i\sqrt{-\varepsilon x_1},\dots,i\sqrt{-\varepsilon x_{k-p+1}},\sqrt{\varepsilon x_{k-p+3}},\dots,\sqrt{\varepsilon x_{q-p+1}},\sqrt{x_{p+1}'},\dots,\sqrt{x_{m-2k-q}'}\right)\right]=\\
&\left[\left(\sqrt{\varepsilon} e^{i \frac{t_1}2},\sqrt{\varepsilon}\sqrt{z_1},\dots,\sqrt{\varepsilon}\sqrt{z_{i-1}},e^{i\frac{t_2}2},\sqrt{w_1},\dots,\sqrt{w_{n-i-1}},\right.\right.\\
&\left.\left.i\sqrt{-x_1'},\dots,i\sqrt{-x_{p-1}'},\sqrt{\varepsilon}i\sqrt{-x_1},\dots,\sqrt{\varepsilon}i\sqrt{-x_{k-p+1}},\sqrt{\varepsilon}\sqrt{x_{k-p+3}},\dots,\sqrt{\varepsilon}\sqrt{x_{q-p+1}},\sqrt{x_{p+1}'},\dots,\sqrt{x_{m-2k-q}'}\right)\right],
\end{aligned}
\]
from which we read immediately that the morphism~\eqref{eq-square} maps diffeomorphically $C_{i,B}^+\times C_{n-i,([m]\smallsetminus B)\sqcup\{\bullet\}}^+$ to $C_{i,B_1,B_2}^+\times C_{n-i,[k]\smallsetminus B_1,[l]\smallsetminus B_2}^+$.

Finally, we consider the case $|B|\neq 0$, such that the maximum of $B=[j]$ with $j\leq k$.
The interior of the corresponding boundary stratum of $C_{n,m}^+$, resp.\ $C_{n,k,l}^+$, is $C_{i,j}^+\times C_{n-i,m-j+1}^+$, resp.\ $C_{i,j}^+\times C_{n-i,k-j+1,l}^+$, and corresponding local coordinates are given by 
\[
\begin{aligned}
C_{i,j}^+\times C_{n-i,m-j+1}^+\ni &\left(\left(i,z_1,\dots,z_{i-1},x_1,\dots,x_j\right),\left(e^{it},w_1,\dots,w_{n-i-1},x_1',\dots,x_{k-j+1}',0,x_{k-j+2}',\dots,x_{m-j+1}'\right)\right),\ \text{resp.}\\
C_{i,j}^+\times C_{n-i,k-j+1,l}^+\ni &\left(\left(i,z_1,\dots,z_{i-1},x_1,\dots,x_j\right),\left(e^{it},w_1,\dots,w_{n-i-1},i x_1',\dots,i x_{k-j+1}',y_1',\dots,y_l'\right)\right),
\end{aligned}
\]
where $t$ is in $(0,\pi)$, resp.\ $(0,\frac{\pi}2)$, all points in $\mathbb H^+$ and $Q^{+,+}$ are distinct in both expressions.
In the first, resp.\ second, expression, the $x_i$ and $x_i'$ are lexicographically ordered, resp.\ $x_1'>\cdots x_{k-j+1}'>0$ and $0<y_1'<\cdots< y'_l$.

Choosing a positive number $\varepsilon$ sufficiently small, we now write local coordinates of $\mathcal C_{n,m}^+$, resp.\ $\mathcal C_{n,k,l}^+$, near the interior of the boundary stratum $\mathcal C_{i,j}^+\times\mathcal C_{n-i,m-j+1}^+$, resp.\ $\mathcal C_{i,j}^+\times \mathcal C_{n-i,k-j+1,l}^+$, namely
\[
\begin{aligned}
&\left[\left(x_1'+\varepsilon i,x_1'+\varepsilon z_1,\dots,x_1'+\varepsilon z_{i-1},e^{it},w_1,\dots,w_{n-i-1},x_1',x_1'+\varepsilon x_1,\dots,x_1'+\varepsilon x_j,x_2',\dots,x_{k-j+1}',0,\right.\right.\\
&\left.\left. x_{k-j+2}',\dots,x_{m-j+1}'\right)\right],\ \text{resp.}\\
&\left[\left(ix_1'+\varepsilon,ix_1'-i\varepsilon z_1,\dots,ix_1'-i\varepsilon z_{i-1},e^{it},w_1,\dots,w_{n-i-1},ix_1',ix_1-i'\varepsilon x_1,\dots,ix_1'-i\varepsilon x_j,ix_2',\dots,ix_{k-j+1}',y_1',\dots,y_l'\right)\right].
\end{aligned}
\]
If we apply the morphism~\eqref{eq-square} to the first of the two previous expressions, we get
\[
\begin{aligned}
&\left[\left(\sqrt{x_1'+\varepsilon i},\sqrt{x_1'+\varepsilon z_1},\dots,\sqrt{x_1'+\varepsilon z_{i-1}},\sqrt{e^{it}},\sqrt{w_1},\dots,\sqrt{w_{n-i-1}},\right.\right.\\
&\left.\left.\sqrt{x_1'},\sqrt{x_1'+\varepsilon x_1},\dots,\sqrt{x_1'+\varepsilon x_j},\sqrt{x_2'},\dots,\sqrt{x_{k-j+1}'},\sqrt{x_{k-j+2}'},\dots,\sqrt{x_{m-j+1}'}\right)\right].
\end{aligned}
\]
Once again, we find
\[
\begin{aligned}
\sqrt{x_1'+\varepsilon i}&=i\sqrt{-x_1'}+\frac{\varepsilon}2\frac{1}{\sqrt{-x_1'}}+\mathcal O(\varepsilon^2),\ &\ \sqrt{x_1'+\varepsilon iz_e}&=i\sqrt{-x_1'}-i\frac{\varepsilon}2\frac{z_e}{\sqrt{-x_1'}}+\mathcal O(\varepsilon^2),\\ 
\sqrt{x_1'+\varepsilon x_e}&=i\sqrt{-x_1'}-i\frac{\varepsilon}2\frac{x_e}{\sqrt{-x_1'}}+\mathcal O(\varepsilon^2), 
\end{aligned}
\]
and using the same arguments as in the previous computations, we see that the morphism~\eqref{eq-square} maps $C_{i,j}^+\times C_{n-i,m-j+1}^+$ diffeomorphically to $C_{i,j}^+\times C_{n-i,k-j+1,l}^+$.
\end{proof}

\subsection{The choice of propagators}\label{ss-2-2}
We now discuss the propagators needed for the computations in the framework of (bi)quantization.
In particular, we discuss in detail the $4$-colored propagators: we will mainly work here with the $4$-colored propagators as introduced originally in~\cite{CFb}, and used extensively in~\cite{CT}.
The point is that we will view the biquantization techniques in~\cite{CT} in the framework of the $2$-brane formality of~\cite{CFFR}.
In~\cite{CFFR}, the authors preferred to work with the $4$-colored propagators on $\mathcal C_{2,1}^+$, in order to use the (simpler) compactified configuration spaces $\mathcal C_{n,m}^+$ of Kontsevich's type: in order to tie in with the computations in~\cite{CT}, we want to establish a more precise relationship than the one sketched in~\cite{CFFR} about the $4$-colored propagators in~\cite{CT} and in~\cite{CFFR}.

\subsubsection{The Kontsevich propagator}\label{sss-2-2-1}
We consider a pair $(z_1,z_2)$ of distinct points in $\mathbb H^+$, and we associate to it a closed $1$-form by the formula
\begin{equation}\label{eq-prop-K}
\omega(z_1,z_2)=\frac{1}{2\pi}\left[\mathrm d\ \mathrm{arg}(z_1-z_2)-\mathrm d\ \mathrm{arg}(\overline z_1-z_2)\right]=\frac{1}{2\pi}\left[\mathrm d\ \mathrm{arg}(z_1-z_2)+\mathrm d\ \mathrm{arg}(z_1-\overline z_2)\right].
\end{equation}
In Formula~\eqref{eq-prop-K}, the function $\mathrm{arg}(z)$ denotes the Euclidean angle of the complex number $z$: it can be made into a smooth function by restricting its domain of definition on $\mathbb C\smallsetminus i\mathbb R^-$.
In particular, the restriction of $\mathrm{arg}(z)$ on $\mathbb H^+$ is a smooth function.
However, we want to consider $\omega(z_1,z_2)$ as a closed $1$-form on $(\mathbb H^+\times\mathbb H^+)\smallsetminus \Delta$, $\Delta$ being the diagonal in $\mathbb H^+\times\mathbb H^+$: as such, $\omega(z_1,z_2)$ is the sum of a closed form on $(\mathbb H^+\times\mathbb H^+)\smallsetminus \Delta$ and of an exact $1$-form, where the corresponding function is $\mathrm{arg}(z_1-\overline z_2)/2\pi$.
We observe, for the sake of later computations (see~\cite{VdB} for a very nice application of this idea), that the closed $1$-form can be made into a truly exact $1$-form by restricting the domain of definition to 
\[
\{(z_1,z_2)\in(\mathbb H^+\times\mathbb H^+)\smallsetminus \Delta:\ \mathrm{Re}(z_1)=\mathrm{Re}(z_2)\Rightarrow \mathrm{Im}(z_1)>\mathrm{Im}(z_2)\}.
\]

It is not difficult to prove that the $1$-form~\eqref{eq-prop-K} descends to $C_{2,0}^+$; a bit more involved is the proof that it extends to a smooth $1$-form $\omega$ on the compactified configuration space $\mathcal C_{2,0}^+$.
The function $\eta(z_1,z_2)=\mathrm{arg}(z_1-\overline z_2)/2\pi$ also descends to $C_{2,0}^+$ and extends to a smooth function on $\mathcal C_{2,0}^+$. 
\begin{Lem}\label{l-prop-K}
The closed $1$-form~\eqref{eq-prop-K} determines a smooth, closed $1$-form $\omega$ on $\mathcal C_{2,0}^+$, which further enjoys the following properties:
\begin{itemize}
\item[$i)$] \[
\omega\vert_{\mathcal C_2\times\mathcal C_{1,0}^+}=\mathrm d\varphi,
\]
where $\mathrm d\varphi$ denotes (improperly) the normalized volume form of $\mathcal C_2\cong S^1$;
\item[$ii)$] 
\[
\omega\vert_{\mathcal C_{1,0}^+\times\mathcal C_{1,1}^+}=0,
\]
where $\mathcal C_{1,0}^+\times\mathcal C_{1,1}^+$ denotes the boundary stratum of $\mathcal C_{2,0}^+$ corresponding to the approach of the first argument $z_1$ to $\mathbb R$.
\end{itemize} 
The function $\eta(z_1,z_2)$ determines a smooth function $\eta$ on $\mathcal C_{2,0}^+$, which restricts on the boundary stratum $\mathcal C_2\times\mathcal C_{1,0}^+$ to the constant function $\pi/2$; observe that $\mathcal C_{1,0}^+\cong \{i\}$.
\end{Lem}
The $1$-form $\omega$ is usually called Kontsevich's angle form~\cite[Subsection 6.2]{K}: it will be useful, for certain computations, to recall that Kontsevich's angle function is the sum of a closed $1$-form and of an exact $1$-form, constructed by means of the function $\eta$.

We finally observe that the natural involution $(z_1,z_2)\overset{\tau}\mapsto (z_2,z_1)$ of $(\mathbb H^+\times\mathbb H^+)\smallsetminus \Delta$ yields an involution $\tau$ of $\mathcal C_{2,0}^+$: we may then consider two Kontsevich's angle forms $\omega^{\pm}$ defined through
\[
\omega^+=\omega,\ \omega^-=\tau^*(\omega).
\]
The angle forms $\omega^\pm$ have been first introduced in~\cite{CFb,CF}: they have opposite boundary conditions when one of their arguments approaches $\mathbb R$, as can be easily deduced from Lemma~\ref{l-prop-K}.

\subsubsection{The $4$-colored propagators on $\mathcal C_{2,1}^+$}\label{sss-2-2-2}
We consider a triple $(z_1,z_2,x)$ in $(\mathbb H^+\times\mathbb H^+)\smallsetminus \Delta\times \mathbb R$.

There is a natural smooth projection from $(\mathbb H^+\times\mathbb H^+)\smallsetminus \Delta\times \mathbb R$ to $(\mathbb H^+\times\mathbb H^+)\smallsetminus \Delta$, thus we may consider the pull-back $\omega^{+,+}$ of the closed $1$-form $\omega^+(z_1,z_2)$ to $(\mathbb H^+\times\mathbb H^+)\smallsetminus \Delta\times \mathbb R$.
We set $\omega^{-,-}(z_1,z_2,x)$ to be the pull-back of $\omega^-$ w.r.t.\ the very same projection.

We recall the complex square root discussed in Subsubsection~\ref{sss-2-1-3}: as already remarked, it is a biholomorphism from $\mathbb H^+$ to $Q^{+,+}$, and we associate to a triple $(z_1,z_2,x)$ in $(\mathbb H^+\times\mathbb H^+)\smallsetminus \Delta\times \mathbb R$ a pair $(\sqrt{z_1-x},\sqrt{z_2-x})$ in $(Q^{+,+}\times Q^{+,+})\smallsetminus \Delta$ (compare with the morphism of Proposition~\ref{p-square}).
We then set
\[
\begin{aligned}
\omega^{+,-}(z_1,z_2,x)&=\frac{1}{2\pi}\left[\mathrm d\ \mathrm{arg}\!\left(\sqrt{z_1-x}-\sqrt{z_2-x}\right)+\mathrm d\ \mathrm{arg}\!\left(\overline{\sqrt{z_1-x}}-\sqrt{z_2-x}\right)-\right.\\
&\phantom{=\frac{1}{2\pi}[}\left.-\mathrm d\ \mathrm{arg}\!\left(\overline{\sqrt{z_1-x}}+\sqrt{z_2-x}\right)-\mathrm d\ \mathrm{arg}\!\left(\sqrt{z_1-x}+\sqrt{z_2-x}\right)\right]\\
\omega^{-,+}(z_1,z_2,x)&=\frac{1}{2\pi}\left[\mathrm d\ \mathrm{arg}\!\left(\sqrt{z_1-x}-\sqrt{z_2-x}\right)-\mathrm d\ \mathrm{arg}\!\left(\overline{\sqrt{z_1-x}}-\sqrt{z_2-x}\right)+\right.\\
&\phantom{=\frac{1}{2\pi}[}\left.+\mathrm d\ \mathrm{arg}\!\left(\overline{\sqrt{z_1-x}}+\sqrt{z_2-x}\right)-\mathrm d\ \mathrm{arg}\!\left(\sqrt{z_1-x}+\sqrt{z_2-x}\right)\right].
\end{aligned}
\]
The two $1$-forms $\omega^{+,-}$ and $\omega^{-,+}$ are smooth and obviously closed on $(\mathbb H^+\times\mathbb H^+)\smallsetminus \Delta\times \mathbb R$.

We need to characterize more explicitly the compactified configuration space $\mathcal C_{2,1}^+$ (for whose more precise description we refer to~\cite[Section 5]{CFFR}): here, we content ourselves to describe all boundary strata of codimension $1$, which we depict as follows
\bigskip
\begin{center}
\resizebox{0.75 \textwidth}{!}{\input{I-cube_Kosz_1.pstex_t}}\\
\text{Figure 3 - The boundary strata of codimension $1$ of $\mathcal C_{2,1}^+$} \\
\end{center}
\bigskip
We observe that the boundary stratum $\alpha$ corresponds to $\mathcal C_2\times\mathcal C_{1,1}^+$, the boundary strata $\beta$ and $\gamma$ to two copies of $\mathcal C_{2,0}^+\times\mathcal C_{0,2}^+$, the boundary strata $\delta$ and $\varepsilon$ to two copies of $\mathcal C_{1,1}^+\times\mathcal C_{1,1}^+$, and the boundary strata $\eta$, $\theta$, $\zeta$ and $\xi$ to four copies of $\mathcal C_{1,0}^+\times\mathcal C_{1,2}^+$.
When it is clear from the context, we will omit to write the projections $\pi_i$, $i=1,2$, from the these spaces to the each of the factors.
We finally recall, once again, that the function $\mathrm{arg}(z)$ is well-defined and smooth on $\mathbb H^+$: in particular, the function $\eta$ from Lemma~\ref{l-prop-K}, Subsubsection~\ref{sss-2-2-1}, yields a smooth function (denoted again by $\eta$) on $\mathcal C_{1,1}^+$, when the second argument approaches $\mathbb R$.
In more down-to-earth terms, $\eta=\mathrm{arg}(z-x)/2\pi$.

It is not difficult to prove that the $4$ $1$-forms $\omega^{+,+}$, $\omega^{+,-}$, $\omega^{-,+}$ and $\omega^{-,-}$ descend to smooth, closed $1$-forms on $C_{2,1}^+$.
In fact, as the following Lemma shows (for whose proof we refer to~\cite[Lemma 5.4]{CFFR}), these in turn extend to smooth, closed $1$-forms on the compactified configuration space $\mathcal C_{2,1}^+$.
\begin{Lem}\label{l-CF}
The $1$-forms $\omega^{+,+}$, $\omega^{+,-}$, $\omega^{-,+}$ and $\omega^{-,-}$ determine smooth, closed $1$-forms on the compactified configuration space $\mathcal C_{2,1}^+$, which enjoy the following properties:
\begin{itemize}
\item[$i)$] \[
\omega^{+,+}\vert_\alpha=\mathrm d\varphi,\ \omega^{+,-}\vert_\alpha=\mathrm d\varphi-\mathrm d\eta,\  \omega^{-,+}\vert_\alpha=\mathrm d\varphi-\mathrm d\eta,\ \omega^{-,-}\vert_\alpha=\mathrm d\varphi,
\] 
where $\mathrm d\varphi$ is the normalized volume form of $\mathcal C_2\cong S^1$.
\item[$ii)$] \[
\begin{aligned}
\omega^{+,+}\vert_\beta&=\omega^+,\ & \omega^{+,-}\vert_\beta&=\omega^+,\ & \omega^{-,+}\vert_\beta&=\omega^-,\ & \omega^{-,-}\vert_\beta&=\omega^-\quad \text{and}\\
\omega^{+,+}\vert_\gamma&=\omega^+,\ & \omega^{+,-}\vert_\gamma&=\omega^-,\ & \omega^{-,+}\vert_\gamma&=\omega^+,\ & \omega^{-,-}\vert_\gamma&=\omega^-,
\end{aligned}
\]
where $\omega^\pm$ have to be understood on $\mathcal C_{2,0}^+$.
\item[$iii)$] \[
\begin{aligned}
&\omega^{+,+}\vert_\delta=\omega^{+,-}\vert_\delta=\omega^{-,+}\vert_\delta=0,\\
&\omega^{+,-}\vert_\varepsilon=\omega^{-,+}\vert_\varepsilon=\omega^{-,-}\vert_\varepsilon=0.
\end{aligned}
\] 
\item[$iv)$] \[
\begin{aligned}
&\omega^{+,-}\vert_\eta=\omega^{-,-}\vert_\eta=0,\ & &\omega^{+,+}\vert_\theta=\omega^{-,+}\vert_\theta=0,\\
&\omega^{-,+}\vert_\zeta=\omega^{-,-}\vert_\zeta=0,\ & &\omega^{+,+}\vert_\xi=\omega^{+,-}\vert_\xi=0.
\end{aligned}
\]
\end{itemize}
\end{Lem}

\subsubsection{The $4$-colored propagators on $\mathcal C_{2,0,0}^+$}\label{sss-2-2-3}
We now define on $(Q^{+,+}\times Q^{+,+})\smallsetminus \Delta$ $4$ closed, smooth $1$-forms, which, by an (apparent) abuse of notation, are denoted by $\omega^{\pm,\pm}$: namely, we set
\[
\begin{aligned}
\omega^{+,+}&=\frac{1}{2\pi}\left[\mathrm d\ \mathrm{arg}(z_1-z_2)-\mathrm d\ \mathrm{arg}(\overline z_1-z_2)-\mathrm d\ \mathrm{arg}(\overline z_1+z_2)+\mathrm d\ \mathrm{arg}(z_1+z_2)\right],\\
\omega^{+,-}&=\frac{1}{2\pi}\left[\mathrm d\ \mathrm{arg}(z_1-z_2)+\mathrm d\ \mathrm{arg}(\overline z_1-z_2)-\mathrm d\ \mathrm{arg}(\overline z_1+z_2)-\mathrm d\ \mathrm{arg}(z_1+z_2)\right],\\
\omega^{-,+}&=\frac{1}{2\pi}\left[\mathrm d\ \mathrm{arg}(z_1-z_2)-\mathrm d\ \mathrm{arg}(\overline z_1-z_2)-\mathrm d\ \mathrm{arg}(\overline z_1+z_2)+\mathrm d\ \mathrm{arg}(z_1+z_2)\right],\\
\omega^{+,+}&=\frac{1}{2\pi}\left[\mathrm d\ \mathrm{arg}(z_1-z_2)+\mathrm d\ \mathrm{arg}(\overline z_1-z_2)+\mathrm d\ \mathrm{arg}(\overline z_1+z_2)+\mathrm d\ \mathrm{arg}(z_1+z_2)\right],
\end{aligned}
\]
for an element $(z_1,z_2)$ of $(Q^{+,+}\times Q^{+,+})\smallsetminus \Delta$.

We first observe that the last three summands in the previous $1$-forms are exact $1$-forms: namely, as has been previously remarked, the function $\mathrm{arg}(z)$ is smooth and well-defined on $\mathbb C\smallsetminus (i\mathbb R^-\sqcup\{0\})$, hence the three functions appearing in the last three summands of the previous formul\ae\ are well-defined and smooth on $(Q^{+,+}\times Q^{+,+})\smallsetminus \Delta$.

It is not difficult to prove that the closed $1$-forms $\omega^{\pm,\pm}$ descends to smooth, closed $1$-forms on the open configuration space $C_{2,0,0}^+$, and that these in turn determine smooth, closed $1$-forms $\omega^{\pm,\pm}$ on the compactified configuration space $\mathcal C_{2,0,0}^+$.

Because of the results of Subsubsection~\ref{sss-2-1-3}, we already know that there is a diffeomorphism between $\mathcal C_{2,0,0}^+$ and $\mathcal C_{2,1}^+$, which smoothly extends to the compactified configuration spaces the diffeomorphism
\[
C_{2,1}^+\ni [(z_1,z_2,x)]\mapsto [(\sqrt{z_1-x},\sqrt{z_2-x})]\in C_{2,0,0}^+.
\]
We leave it to the reader to reinterpret on $\mathcal C_{2,0,0}^+$ the boundary strata of codimension $1$ of $\mathcal C_{2,1}^+$.

It is not difficult to prove that the pull-backs w.r.t.\ the morphism~\eqref{eq-square} from $\mathcal C_{2,1}^+$ to $\mathcal C_{2,0,0}^+$ of the $1$-forms $\omega^{\pm,\pm}$ on $\mathcal C_{2,0,0}^+$ are exactly the $1$ forms $\omega^{\pm,\pm}$ introduced in Subsubsection~\ref{sss-2-2-2}: {\em e.g.} for $\omega^{+,+}$, we have the obvious identity
\[
\omega^{+, +}=\frac1{2\pi} \left[\mathrm{d}\ \mathrm{arg} (z_1^2-z_2^2) -\mathrm{d}\ \mathrm{arg}(\overline{z}_1^2-z_2^2)\right],
\]
whence the claim follows.
Similar arguments work for the other cases.

According to the boundary stratification of $\mathcal C_{2,0,0}^+$, we have the following variant of Lemma~\ref{l-CF}.
\begin{Lem}\label{l-CF-q}
The $1$-forms $\omega^{+,+}$, $\omega^{+,-}$, $\omega^{-,+}$ and $\omega^{-,-}$ determine smooth, closed $1$-forms on the compactified configuration space $\mathcal C_{2,0,0}^+$, which enjoy the following properties:
\begin{itemize}
\item[$i)$] \[
\omega^{+,+}\vert_\alpha=\mathrm d\varphi+\mathrm d\eta,\ \omega^{+,-}\vert_\alpha=\mathrm d\varphi-\mathrm d\eta,\  \omega^{-,+}\vert_\alpha=\mathrm d\varphi-\mathrm d\eta,\ \omega^{-,-}\vert_\alpha=\mathrm d\varphi+\mathrm d \eta,
\]
where $\mathrm d\varphi$ is the normalized volume form of $\mathcal C_2\cong S^1$, and $\eta=\mathrm{arg}(z)/2\pi$ is a well-defined, smooth function on $\mathcal C_{1,0,0}^+$.
\item[$ii)$] \[
\begin{aligned}
\omega^{+,+}\vert_\beta&=\omega^+,\ & \omega^{+,-}\vert_\beta&=\omega^+,\ & \omega^{-,+}\vert_\beta&=\omega^-,\ & \omega^{-,-}\vert_\beta&=\omega^-\quad \text{and}\\
\omega^{+,+}\vert_\gamma&=\omega^+,\ & \omega^{+,-}\vert_\gamma&=\omega^-,\ & \omega^{-,+}\vert_\gamma&=\omega^+,\ & \omega^{-,-}\vert_\gamma&=\omega^-,
\end{aligned}]
\] 
where $\omega^\pm$ have to be understood on $\mathcal C_{2,0}^+$.
\item[$iii)$] \[
\begin{aligned}
&\omega^{+,+}\vert_\delta=\omega^{+,-}\vert_\delta=\omega^{-,+}\vert_\delta=0,\\
&\omega^{+,-}\vert_\varepsilon=\omega^{-,+}\vert_\varepsilon=\omega^{-,-}\vert_\varepsilon=0.
\end{aligned}
\]
\item[$iv)$] \[
\begin{aligned}
&\omega^{+,-}\vert_\eta=\omega^{-,-}\vert_\eta=0,\ & &\omega^{+,+}\vert_\theta=\omega^{-,+}\vert_\theta=0,\\
&\omega^{-,+}\vert_\zeta=\omega^{-,-}\vert_\zeta=0,\ & &\omega^{+,+}\vert_\xi=\omega^{+,-}\vert_\xi=0.
\end{aligned}
\]
\end{itemize}
\end{Lem}
We observe that the $1$-forms $\omega^{\pm,\pm}$, be they defined either on $\mathcal C_{2,1}^+$ or on $\mathcal C_{2,0,0}^+$, satisfy the same boundary conditions $ii)$, $iii)$ and $iv)$; on the other hand, the behavior of the $4$-colored propagators on the boundary strata $\mathcal C_2\times\mathcal C_{1,1}^+$ and $\mathcal C_2\times\mathcal C_{1,0,0}^+$ are quite different.
This can be traced back to the proof of Proposition~\ref{p-square}, when analyzing the shape of the morphism~\eqref{eq-square} on the boundary stratum $\mathcal C_2\times \mathcal C_{1,1}^+$.
Still, we have to be careful about these (seemingly) different boundary conditions for the $4$-colored propagators $\omega^{\pm,\pm}$: namely, the fact that the $4$-colored propagators, quite opposite to Kontsevich's angle form, can be written as a sum of a regular and of a singular term (the $1$-form living on $\mathcal C_{1,1}^+$ or $\mathcal C_{1,0,0}^+$ and on $\mathcal C_2$ respectively) produces a significant change in the application of Stokes' Theorem, which is the fundamental tool for proving the $2$-brane Formality Theorem, from which biquantization follows.

\subsection{Formality Theorems}\label{ss-2-3}
In this Subsection, we recall the $2$-brane Formality Theorem of~\cite{CFFR}, from which we will derive the biquantization techniques we apply later on.
Although the main computations of this Subsection are already contained in~\cite{CFFR}, we review them in some detail because of the following reasons: first, the $2$-brane Formality Theorem has been proved using superpropagators along the same patterns of~\cite{CF}, and superpropagators are better suited for keeping track of all different colors of propagators w.r.t.\ the treatment in~\cite{CT}, and second, because we deserve here a more careful treatment than in~\cite{CFFR} of the $1$-loop correction arising because of the aforementioned regular term in the $4$-colored propagators.
We thus profit of the space here to correct a slight mistake in~\cite[Subsection 7.1]{CFFR} (in the sense that the computations therein are correct, but a subtle point has been missed regarding the multidifferential operator associated to the $1$-loop correction, which we illustrate here in detail) and, more importantly, to correct a more serious mistake in~\cite{CT}, where the regular part of the restriction to the boundary stratum $\mathcal C_2\times \mathcal C_{1,1}^+$ or $\mathcal C_2\times \mathcal C_{1,0,0}^+$ is missing completely.
The correction term arising from the presence of the regular part is responsible for a quantum shift, which will be illustrated explicitly in Section~\ref{s-3}, which is predicted by representation-theoretic arguments and was otherwise absent. 

We also prove a version of~\cite[Lemmata 7.3.1.1, 7.3.3.1]{K} for the $4$-colored propagators: such vanishing lemmata are central in some computations in~\cite{CT} regarding the Harish--Chandra homomorphism.
The main idea of the proof is, once again, Stokes' Theorem, but of course here we have to be a bit more careful and slightly change the final argument.

\subsubsection{Admissible graphs}\label{sss-2-3-1}
Before entering into the technicalities of the $1$-brane and $2$-brane Formality Theorems, we need to spend some words on admissible graphs.

For a pair of non-negative integers $(n,m)$, such that $2n+m-2\geq 0$, we consider the set $\mathcal G_{n,m}$ of admissible graphs of type $(n,m)$: the integer $n$, resp.\ $m$, refers to the number of vertices of the first, resp.\ second type, {\em i.e.} vertices in $\mathbb H^+$, resp.\ on $\mathbb R$.
An admissible graph $\Gamma$ of type $(n,m)$ in the framework of the $1$-brane Formality Theorem~\cite{CF,CFFR} is an oriented graph, which may admit double edges, {\em i.e.} given any two vertices $(v_1,v_2)$, there can be more than one edge connecting $v_1$ to $v_2$, and edges departing from vertices of the second type; it does not possess short loops, {\em i.e.} there can no edge in $\Gamma$ with coincident initial and final point.
The presence of multiple edges and edges departing from $\mathbb R$ is in opposition to the definition of admissible graphs of type $(n,m)$ as in~\cite{K}.

Further, for a triple of non-negative integers $(n,k,l)$, such that $2n+k+l-1\geq 0$, we consider the set $\mathcal G_{n,k,l}$ of admissible graphs of type $(n,k,l)$, where $n$ is the number of vertices of the first type ({\em i.e.} in $Q^{+,+}$), $k$, resp.\ $l$, is the number of vertices of the first type on $i\mathbb R^+$, resp.\ $\mathbb R^+$. 
A general element $\Gamma$ of $\mathcal G_{n,k,l}$ is an oriented graphs with $n$, resp.\ $k+l$, vertices of the first, resp.\ second type, which may admit multiple edges, edges departing from $i\mathbb R^+\sqcup\{0\}\sqcup\mathbb R^+$ and even short loops.

We observe that we may also equivalently consider, for $m=k+l+1$, the set $\mathcal G_{n,m}$ of admissible graphs of type $(n,m)$, consisting of oriented graphs with $n$, resp. $m$, vertices of the first, resp.\ second type ({\em i.e.} lying in $\mathbb H^+$ and on $\mathbb R$ respectively), such that one vertex of the first type is marked and which admit multiple edges, edges departing from $\mathbb R$ and short loops: the notation is abused, but it will be clear from the context if we allow elements of $\mathcal G_{n,m}$ to possess or not short loops, which is the only additional feature that the admissible graphs for the $2$-brane Formality Theorem admit w.r.t.\ the ones in the $1$-brane Formality Theorem.
The algebraic counterpart of the geometric results of Subsubsection~\ref{sss-2-1-3} is the fact that we may freely pass from $\mathcal G_{n,m}$ to $\mathcal G_{n,k,l}$, for $m=k+l+1$, by noting that the vertex labeled by $k+1$ on $\mathbb R$ corresponds to the origin $\{0\}$.

\subsubsection{Superpropagators}\label{sss-2-3-2}
We now pick an admissible graph $\Gamma$ of type $(n,m)$ for the $1$-brane Formality Theorem of~\cite{CF}. 
As $\Gamma$ is of type $(n,m)$, its vertices correspond to a point of $\mathcal C_{n,m}^+$, and an edge $e$ determines a natural projection $\pi_e:\mathcal C_{n,m}^+\to \mathcal C_{2,0}^+$.

If we pick an admissible graph $\Gamma$ of $\mathcal G_{n,m}$ in the framework of the $2$-brane Formality Theorem, then the vertices of $\Gamma$ still define a configuration of points in $\mathcal G_{n,m}$.
An edge $e$ defines, as in the previous situation, either a natural projection $\pi_e:\mathcal C_{n,m}^+\to\mathcal C_{2,1}^+$, if $e=(v_e^i,v_e^f)$, $v_e^i\neq v_e^f$, or $\pi_e:\mathcal C_{n,m}^+\to\mathcal C_{1,1}^+$, if $e$ is a short loop.
The point in $\mathbb R$ in either $\mathcal C_{2,1}^+$ or $\mathcal C_{1,1}^+$ is the marked point of $\mathcal C_{n,m}^+$.
If, equivalently, we consider the corresponding admissible graph $\Gamma$ of type $(n,k,l)$, then an edge $e$ of $\Gamma$ determines either a projection $\pi_e:\mathcal C_{n,k,l}^+\to\mathcal C_{2,0,0}^+$ or $\pi_e:\mathcal C_{n,k,l}^+\to\mathcal C_{1,0,0}^+$.

We now consider the vector space $X=\mathbb K^d$ and two linear (or affine) subspaces $U_i$, $i=1,2$, for which we assume there is a direct sum decomposition
\begin{equation}\label{eq-orth-split}
X=(U_1\cap U_2)\overset{\perp}\oplus (U_1^\perp\cap U_2)\overset{\perp}\oplus (U_1\cap U_2^\perp)\overset{\perp}\oplus (U_1+U_2)^\perp,
\end{equation} 
w.r.t.\ a chosen inner product over $X$.
Clearly, we have
\begin{equation*}
U_1=(U_1\cap U_2)\overset{\perp}\oplus (U_1\cap U_2^\perp),\ U_2=(U_1\cap U_2)\overset{\perp}\oplus (U_1^\perp\cap U_2).
\end{equation*}

We choose linear coordinates $\{x_i\}$ on $X$ which are adapted to the orthogonal decomposition~\eqref{eq-orth-split}, {\em i.e.} there are two non-disjoint subsets $I_i$, $i=1,2$, of $[d]$, such that 
\[
[d]=\left(I_1\cap I_2\right)\sqcup\left(I_1\cap I_2^c\right)\sqcup\left(I_1^c\cap I_2\right)\sqcup\left(I_1^c\cap I_2^c\right),
\]
w.r.t.\ which $\{x_i\}$ is a set of linear coordinates on $U_1\cap U_2$, $U_1\cap U_2^\perp$, $U_1^\perp\cap U_2$ or $(U_1+U_2)^\perp$, if the index $i$ belongs to $I_1\cap I_2$, $I_1\cap I_2^c$, $I_1^c\cap I_2$ or $I_1^c\cap I_2^c$ respectively. 
Accordingly, for $I$ either one of the previous subsets of $[d]$, and $e$ an edge of admissible graph $\Gamma$ of type $(n,m)$, we set
\[
\tau_e^I=\sum_{i\in I}\iota_{\mathrm d x_i}^{(v_e^i)}\partial_{x_i}^{(v_e^f)}\in\mathrm{End}\!\left(T_\mathrm{poly}(X)^{\otimes(m+n)}\right),\ T_\mathrm{poly}(X)=\mathrm S(X^*)\otimes\wedge^\bullet X,
\]
and $\partial_{x_i}^{(v)}$ denotes the action of the differential operator on the copy of $T_\mathrm{poly}(X)$ sitting at the $v$-th position, and similarly for $\iota_{\mathrm d x_i}^{(v)}$.
We observe that $\tau_e^I$ is well-defined and has degree $-1$ w.r.t.\ the natural grading on $T_\mathrm{poly}(X)$.

We now set
\begin{align*}
A&=\mathrm S(U_1^*)\otimes \wedge (X/U_1)=\mathrm S(U_1^*)\otimes \wedge (U_1^\perp\cap U_2)\otimes \wedge (U_1+U_2)^\perp,\\
B&=\mathrm S(U_2^*)\otimes \wedge (X/U_2)=\mathrm S(U_2^*)\otimes \wedge (U_1\cap U_2^\perp)\otimes \wedge (U_1+U_2)^\perp,\\
K&=\mathrm S((U_1\cap U_2)^*)\otimes \wedge (U_1+U_2)^\perp.
\end{align*}
It is clear that $A$ and $B$ both admit a (trivial) structure of $A_\infty$-algebra, and $K$ is naturally an $A$-$B$-bimodule.

With respect to the previously introduced notation, the relevant superpropagators are then given by
\begin{align}
\label{eq-A-form}\omega^A_e&=\pi_e^*(\omega^+)\otimes\left(\tau^{I_1\cap I_2}_e+\tau^{I_1\cap I_2^c}_e\right)+\pi_e^*(\omega^-)\otimes\left(\tau^{I_1^c\cap I_2}_e+\tau^{I_1^c\cap I_2^c}_e\right),\\
\label{eq-B-form}\omega^B_e&=\pi_e^*(\omega^+)\otimes\left(\tau^{I_1\cap I_2}_e+\tau^{I_1^c\cap I_2}_e\right)+\pi_e^*(\omega^-)\otimes\left(\tau^{I_1\cap I_2^c}_e+\tau^{I_1^c\cap I_2^c}_e\right),\\
\label{eq-K-form}\omega^K_e&=\pi_e^*(\omega^{+,+})\otimes \tau^{I_1\cap I_2}_e+\pi_e^*(\omega^{+,-})\otimes \tau^{I_1\cap I_2^c}_e+\pi_e^*(\omega^{-,+})\otimes \tau^{I_1^c\cap I_2}_e+\pi_e^*(\omega^{-,-})\otimes \tau^{I_1^c\cap I_2^c}_e,
\end{align} 
for an edge $e=(v_e^i,v_e^f)$, $v_e^i\neq v_e^f$, of an admissible graph $\Gamma$ of type $(n,m)$.

We observe that the superpropagators~\eqref{eq-A-form},~\eqref{eq-B-form} and\eqref{eq-K-form} are closed $1$-forms on $\mathcal C_{2,0}^+$ and $\mathcal C_{2,1}^+$ with values in $\mathrm{End}\!\left(T_\mathrm{poly}(X)^{\otimes(m+n)}\right)$ (of course, $A$, $B$ and $K$ may be viewed as subalgebras of $T_\mathrm{poly}(X)$).
Equivalently, we may regard the superpropagator~\eqref{eq-K-form} as a closed $1$-form on $\mathcal C_{2,0,0}^+$ with values in $\mathrm{End}\!\left(T_\mathrm{poly}(X)^{\otimes(m+n)}\right)$.

Lemma~\ref{l-prop-K}, Subsubsection~\ref{sss-2-2-1}, implies the following useful boundary conditions for the superpropagators~\eqref{eq-A-form} and~\eqref{eq-B-form}:
\begin{itemize}
\item[$i)$] their restrictions to the boundary stratum $\mathcal C_{2,0}^+\times\mathcal C_{1,0}^+$ equal
\[
\omega_e^A\vert_{\mathcal C_{2,0}^+\times\mathcal C_{1,0}^+}=\omega_e^B\vert_{\mathcal C_{2,0}^+\times\mathcal C_{1,0}^+}=\mathrm d\varphi\otimes\tau_e^{[d]}; 
\]  
\item[$ii)$] their restrictions to the boundary stratum $\mathcal C_{1,0}^+\times\mathcal C_{1,1}^+$ corresponding to the approach of the first, resp.\ second, argument to $\mathbb R$ equal
\[
\begin{aligned}
\omega_e^A\vert_{\mathcal C_{1,0}^+\times\mathcal C_{1,1}^+}&=\pi_e^*(\omega^-)\otimes\left(\tau^{I_1^c\cap I_2}_e+\tau^{I_1^c\cap I_2^c}_e\right),\ &\text{resp.}\quad \omega_e^A\vert_{\mathcal C_{1,0}^+\times\mathcal C_{1,1}^+}&=\pi_e^*(\omega^+)\otimes\left(\tau^{I_1\cap I_2}_e+\tau^{I_1\cap I_2^c}_e\right),\\
\omega_e^B\vert_{\mathcal C_{1,0}^+\times\mathcal C_{1,1}^+}&=\pi_e^*(\omega^-)\otimes\left(\tau^{I_1\cap I_2^c}_e+\tau^{I_1^c\cap I_2^c}_e\right),\ &\text{resp.}\quad \omega_e^A\vert_{\mathcal C_{1,0}^+\times\mathcal C_{1,1}^+}&=\pi_e^*(\omega^+)\otimes\left(\tau^{I_1\cap I_2}_e+\tau^{I_1^c\cap I_2}_e\right).
\end{aligned}
\]
\end{itemize}
In particular, we see why admissible graphs appearing in the $1$-brane Formality Theorem may admit edges departing from $\mathbb R$, see for more details~\cite{CFb,CF}.

We now concentrate on the boundary conditions for the superpropagator~\eqref{eq-K-form} on $\mathcal C_{2,1}^+$: Lemma~\ref{l-CF} yields
\begin{itemize}
\item[$i)$] the restriction of the superpropagator~\eqref{eq-K-form} to the boundary stratum $\alpha$ of $\mathcal C_{2,1}^+$ equals
\[
\omega_e^K\vert_\alpha=\mathrm d\varphi\otimes \tau_e^{[d]}-\mathrm d\eta\otimes\left(\tau_e^{I_1\cap I_2^c}+\tau_e^{I_1^c\cap I_2}\right);
\]
\item[$ii)$] the restriction of the superpropagator~\eqref{eq-K-form} to the boundary strata $\beta$ and $\gamma$ equals
\[
\omega_e^K\vert_\beta=\omega_e^A,\ \omega_e^K\vert_\gamma=\omega_e^B;
\]
\item[$iii)$] the restriction of the superpropagator~\eqref{eq-K-form} to the boundary strata $\delta$ and $\varepsilon$ equals
\[
\omega_e^K\vert_\delta=\pi_e^*(\omega^{-,-})\otimes \tau_e^{I_1^c\cap I_2^c},\ \omega_e^K\vert_\varepsilon=\pi_e^*(\omega^{+,+})\otimes \tau_e^{I_1\cap I_2};
\]
\item[$iv)$] the restriction of the superpropagator~\eqref{eq-K-form} to the boundary strata $\\eta$, $\theta$, $\zeta$ and $\xi$ equals
\[
\begin{aligned}
\omega_e^K\vert_\eta&=\pi_e^*(\omega^{+,+})\otimes \tau_e^{I_1\cap I_2^c}+\pi_e^*(\omega^{-,+})\otimes \tau_e^{I_1^c\cap I_2},\ &\ \omega_e^K\vert_\theta=\pi_e^*(\omega^{+,-})\otimes \tau_e^{I_1\cap I_2^c}+\pi_e^*(\omega^{-,-})\otimes \tau_e^{I_1^c\cap I_2^c},\\
\omega_e^K\vert_\zeta&=\pi_e^*(\omega^{+,+})\otimes \tau_e^{I_1\cap I_2}+\pi_e^*(\omega^{+,-})\otimes \tau_e^{I_1\cap I_2^c},\ &\ \omega_e^K\vert_\xi=\pi_e^*(\omega^{-,+})\otimes \tau_e^{I_1^c\cap I_2}+\pi_e^*(\omega^{-,-})\otimes \tau_e^{I_1^c\cap I_2^c}.
\end{aligned}
\]
\end{itemize}
If we choose the superpropagator~\eqref{eq-K-form} on $\mathcal C_{2,0,0}^+$, it satisfies the same boundary conditions, with the exception of the first one, which takes the form
\[
\omega_e^K\vert_\alpha=\mathrm d\varphi\otimes \tau_e^{[d]}+\mathrm d\eta\otimes\left(\tau_e^{I_1\cap I_2}+\tau_e^{I_1^c\cap I_2^c}-\tau_e^{I_1\cap I_2^c}-\tau_e^{I_1^c\cap I_2}\right).
\]
For the sake of simplicity, we write $\tau_e^+=\tau_e^{I_1\cap I_2}+\tau_e^{I_1^c\cap I_2^c}$ and $\tau_e^-=\tau_e^{I_1^c\cap I_2}+\tau_e^{I_1\cap I_2^c}$.

We observe that the boundary conditions of type $iii)$ and $iv)$ explain why the admissible graphs appearing in the $2$-brane Formality Theorem admit edges departing from $\mathbb R$; when considering such admissible graphs in $\mathbb Q^{+,+}\sqcup i\mathbb R^+\sqcup\mathbb R^+\sqcup\{0\}$, we observe that the boundary conditions $iii)$ imply that such graphs admit edges departing from or arriving at the origin.

We now deal with the so-called superloop propagator: its origin will be explained carefully in the proof of the $2$-brane Formality Theorem, which will come later on.
For the time being, we content ourselves by noting that the superloop propagator appear only first in the $2$-brane Formality Theorem as a consequence of the boundary condition $i)$ satisfied by the superpropagator~\eqref{eq-K-form}, more precisely it arises because of the ``regular term'' containing the form $\mathrm d\eta$.

With the same notation as before, the superloop propagator associated to a short loop $e$ of an admissible graph $\Gamma$ of type $(n,m)$ is defined as the closed $1$-form on $\mathcal C_{2,1}^+$ with values in $\mathrm{End}\!\left(T_\mathrm{poly}(X)^{(m+n)}\right)$ 
\[
\omega_e^K=\frac{1}2\pi_e^*(\mathrm d\eta)\otimes (\mathrm{div}_{(v)}^+-\mathrm{div}^-_{(v)}),\ e=(v,v),
\]
where 
\[
\mathrm{div}^+_{(v)}=\sum_{k\in (I_1\cap I_2)\sqcup (I_1^c\cap I_2^c)}\iota_{\mathrm d x_k}^{(v)}\partial_{x_k}^{(v)},\quad \mathrm{div}^-_{(v)}=\sum_{k\in (I_1^c\cap I_2)\sqcup (I_1\cap I_2^c)}\iota_{\mathrm d x_k}^{(v)}\partial_{x_k}^{(v)}.
\]
We observe that the superloop propagator is exact: this fact will be used in all subsequent computations.
Notice that the superloop propagator on $\mathcal C_{2,0,0}^+$ is defined by the same formula without the rescaling by $1/2$ (because of the morphism~\eqref{eq-square} from $\mathcal C_{2,1}^+$ to $\mathcal C_{2,0,0}^+$).

\subsubsection{The formality morphisms}\label{sss-2-3-3}
We consider $X$, $U_1$ and $U_2$ as before, to which we associate the graded vector spaces $A$, $B$ and $K$.
Using the superpropagators~\eqref{eq-A-form},~\eqref{eq-B-form} and~\eqref{eq-K-form}, and keeping in mind the notation in the previous Subsubsections, we set
\begin{align}
\label{eq-c-A}\mathcal O^A_\Gamma(\gamma_1|\cdots|\gamma_n|a_1|\cdots|a_m)&=\mu_{n+m}^B\left(\int_{\mathcal C_{n,m}^+}\prod_{e\in E(\Gamma)}\omega^A_e(\gamma_1|\cdots|\gamma_n|a_1|\cdots|a_m)\right),\\
\label{eq-c-B}\mathcal O^B_\Gamma(\gamma_1|\cdots|\gamma_n|a_1|\cdots|a_m)&=\mu_{n+m}^B\left(\int_{\mathcal C_{n,m}^+}\prod_{e\in E(\Gamma)}\omega^B_e(\gamma_1|\cdots|\gamma_n|b_1|\cdots|b_m)\right),\\
\label{eq-c-K}\mathcal O_\Gamma^K(\gamma_1|\cdots|\gamma_n|a_1|\cdots|a_k|k|b_1|\cdots|b_l)&=\mu_{m+n}^K\left(\int_{\mathcal C_{n,m}^+}\prod_{e\in E(\Gamma)}\omega^K_e(\gamma_1|\cdots|\gamma_n|a_1|\cdots|a_k|k|b_1|\cdots|b_l)\right),
\end{align}
where $\gamma_i$, $i=1,\dots,n$, are elements of $T_\mathrm{poly}(X)$, $a_i$ and $b_i$ are elements of $A$ and $B$ respectively, $k$ is an element of $K$; $E(\Gamma)$ is the set of edges of an admissible graph $\Gamma$ of type $(n,m)$; $\mu^A$, $\mu^B$ and $\mu^K$ denotes the multiplication map on $T_\mathrm{poly}(X)$, followed by the projection onto $A$, $B$ and $K$ respectively.

Since $\Gamma$ may have multiple edges, there is a combinatorial subtlety to be taken into account: in all previous formul\ae, whenever there are multiple edges between two vertices $(v^i,v_f)$, for $v^i\neq v_f$, we must divide by the factorial of the number of such edges.
We observe that short loops cannot be multiple edges, as the superpropagator for a short loop squares obviously to $0$.

We also observe that the product on formul\ae~\eqref{eq-c-A},~\eqref{eq-c-B} and~\eqref{eq-c-K} are well-defined and do not depend on the order of the factors: namely, the total degree of any superpropagator appearing in these formul\ae\ is $0$, as the $1$-form piece has (form) degree $1$, while the multidifferential operator piece has degree $-1$.

Using the multidifferential operators defined in~\eqref{eq-c-A},~\eqref{eq-c-B} and~\eqref{eq-c-K}, we set
\begin{align}
\label{eq-A-mor}\mathcal U_A^n(\gamma_1|\cdots|\gamma_n)(a_1|\dots|a_m)&=(-1)^{\left(\sum_{i=1}^n|\gamma_i|-1\right)m}\sum_{\Gamma\in\mathcal{G}_{n,m}}\mathcal{O}_{\Gamma}^A(\gamma_1|\cdots|\gamma_n|a_1|\cdots|a_m),\\
\label{eq-B-mor}\mathcal U_B^n(\gamma_1|\cdots|\gamma_n)(a_1|\dots|a_m)&=(-1)^{\left(\sum_{i=1}^n|\gamma_i|-1\right)m}\sum_{\Gamma\in\mathcal{G}_{n,m}}\mathcal{O}_{\Gamma}^B(\gamma_1|\cdots|\gamma_n|b_1|\cdots|b_m),\\
\label{eq-K-mor}\mathcal U_K^n(\gamma_1|\cdots|\gamma_n)(a_1|\cdots|a_k|k|b_1|\cdots|b_l)&=(-1)^{\left(\sum_{i=1}^n|\gamma_i|-1\right)m}\sum_{\Gamma\in\mathcal{G}_{n,m}}\mathcal{O}_{\Gamma}^K(\gamma_1|\cdots|\gamma_n|a_1|\cdots|a_k|k|b_1|\cdots|b_l),\\
\label{eq-A_inf-bimod}\mathrm d_K^{k,l}(a_1|\cdots|a_k|k|b_1|\cdots|b_l)&=\sum_{\Gamma\in\mathcal G_{0,m}}\mathcal O_\Gamma^K(a_1|\cdots|a_k|k|b_1|\cdots|b_l),
\end{align}
with the above notation.

Some observations are necessary here.
The morphisms~\eqref{eq-A-mor} and~\eqref{eq-B-mor} and~\eqref{eq-K-mor} are multilinear maps from $T_\mathrm{poly}(X)$ to the multidifferential operators on $A$ and $B$ respectively; the morphisms~\eqref{eq-K-mor} and~\eqref{eq-A_inf-bimod} are multilinear maps from $T_\mathrm{poly}(X)$ to the multidifferential operators from $A^{\otimes k}\otimes K\otimes B^{\otimes l}$ to $K$.
All multidifferential operators appearing in the previous formul\ae\ are non-trivial only if the number of edges of the admissible graphs of type $(n,m)$ equals $2n+m-2$: since to each edge of an admissible graph is associated a contraction operator (which lowers degrees by $1$), it follows immediately that the morphisms~\ref{eq-A-mor},~\eqref{eq-B-mor},~\eqref{eq-K-mor} have degree $2-n$, and that the morphism~\eqref{eq-A_inf-bimod} has degree $2-m$.

We refer to~\cite[Section 3]{CFFR} for a short introduction to $A_\infty$-categories in the present framework (see~\cite{Kel,Lef-Has} for more details on $A_\infty$-categories and related issues), which is needed for the statement of the main theorem ($1+2$-brane Formality Theorem) of the present Section.
We only recall that $T_\mathrm{poly}(X)$ has a structure of dg (short for differential graded) Lie algebra with trivial differential and Schouten--Nijenhuis bracket (extending the natural Lie bracket on polynomial vector fields on $X$); similarly, the Hochschild cochain complex of an $A_\infty$-category $\mathcal A$ (roughly, an abelian category, whose spaces spaces of morphisms admit the structure of $A_\infty$-algebras and $A_\infty$-bimodules) is also a dg Lie algebra with Hochschild differential (the $A_\infty$-structure itself) and Gerstenhaber bracket (which is well-defined an any sort of Hochschild cochain complex). 
\begin{Thm}\label{t-form}
We may regard $A$, $B$ as $A_\infty$-algebras, whose only non-trivial Taylor component is given by the corresponding natural (graded) commutative products: then, the morphisms~\eqref{eq-A_inf-bimod} fit into the Taylor components of a non-trivial $A_\infty$ $A$-$B$-bimodule structure over $K$, which restricts to the natural $A$ left- and $B$-right module structures on $K$.

Furthermore, the morphisms~\eqref{eq-A-mor},~\eqref{eq-B-mor} and~\eqref{eq-K-mor} fit into the Taylor components of an $L_\infty$-morphism $\mathcal U$ from $T_\mathrm{poly}(X)$ to the (completed) Hochschild cochain complex of the $A_\infty$-category $\mathcal A$ with two objects $U_i$, $=1,2$, and spaces of morphisms given by 
\[
\mathrm{Hom}_{\mathcal A}(U_1,U_1)=A,\ \mathrm{Hom}_{\mathcal A}(U_1,U_1)=B,\ \mathrm{Hom}_{\mathcal A}(U_1,U_2)=K,\ \mathrm{Hom}_{\mathcal A}(U_1,U_1)=\{0\},
\]
with the respective $A_\infty$-structures.
Finally, the $L_\infty$-morphism $\mathcal U$ extends to an $L_\infty$-quasi-isomorphism by suitably completing the graded vector spaces $A$, $B$, $K$.
\end{Thm}
\begin{proof}
The first claim has been proved in detail in~\cite[Proposition 6.5]{CFFR}, to which we refer.

The second claim splits into three claims, namely $\mathcal U$ consists of three morphisms $\mathcal U_A$, $\mathcal U_B$ and $\mathcal U_K$, where $\mathcal U_A$, resp.\ $\mathcal U_B$, is a pre-$L_\infty$-morphism from $T_\mathrm{poly}(A)$, resp.\ $T_\mathrm{poly}(B)$, to the (completed) Hochschild cochain complex of $A$, resp.\ $B$, and $\mathcal U_K$ is a collection of maps from $T_\mathrm{poly}(X)$ to the mixed component $\mathrm C^\bullet(A,B,K)$ of the (completed) Hochschild cochain complex of the above $A_\infty$-category $\mathcal A$.
Here, we have used the (non-canonical) identification of dg Lie algebras $T_\mathrm{poly}(X)=T_\mathrm{poly}(A)=T_\mathrm{poly}(B)$.

The fact that $\mathcal U_A$ and $\mathcal U_B$ are $L_\infty$-morphisms has been proved in detail in~\cite{CF}; they extend to $L_\infty$-quasi-isomorphisms by suitably completing $A$ and $B$.

The fact that the morphism $\mathcal U_K$ satisfies the required $L_\infty$-identities has been proved in detail in~\cite[Theorem 7.2]{CFFR}: we profit nonetheless for discussing an incorrect issue in the proof regarding the superloop propagator.
The superloop propagator, which has been defined above, is manifestly different from the one considered in~\cite[Subsection 7.1]{CFFR}: the point is that the actual superloop propagator is the correct one.
We may repeat the proof of~\cite[Theorem 7.2]{CFFR} {\em verbatim} until the discussion of boundary strata of codimension $1$ of the form $\mathcal C_A\times C_{([n]\smallsetminus A)\sqcup\{\bullet\},m}^+$, where $|A|=2$: the following discussion on how the corresponding integral contribution looks like is precisely the same, {\em i.e.} the only situation that matter arise when there are at least one and at most two edges connecting the two vertices labeled by $A$, {\em i.e.} pictorially
\bigskip
\begin{center}
\resizebox{0.6 \textwidth}{!}{\input{b-i-contr.pstex_t}}\\
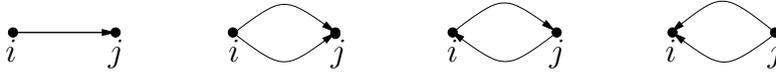
\text{Figure 4 - The four possible loop-free subgraphs $\Gamma_A$ yielding non-trivial boundary contributions of type $i)$} \\
\end{center}
\bigskip

We are interested only in the contributions from the last three subgraphs (which we denote collectively by $\Gamma_A$).
Taking into account the fact that the second graph $\Gamma_A$ has $2$ multiple edges (thus recalling the normalization factor $1/2$), its contribution equals
\[
\int_{\mathcal C_2}\omega_{\Gamma_A}^K=-\pi_e^*(\mathrm d\eta)\otimes \tau^{[d]}_e\tau^-_e=\frac{1}2\pi_e^*(\mathrm d\eta)\otimes \tau^{[d]}_e(\tau^+_e-\tau_e^-),
\]
where $\pi_e$ is here the projection with respect to the ``phantom'' short loop arising from the contraction of the vertices of the subgraph $\Gamma_A$.
The novelty with respect to the corresponding computations in the proof of~\cite[Theorem 7.2]{CFFR} lies in the re-writing of the second term in the previous chain of equalities; of course, we have used the obvious fact that $(\tau_e^{[d]})^2=0$.
The fourth graph in Figure 4 yields a similar contribution.
The third graph, on the other hand, yields the contribution
\[
\int_{\mathcal C_2}\omega_{\Gamma_A}^K=-\pi_{e}^*(\mathrm d\eta)\otimes\tau^{[d]}_{e_1}\tau^-_{e_2}-\pi_{e}^*(\mathrm d\eta)\otimes\tau^{[d]}_{e_2}\tau^-_{e_1}=\frac{1}2\pi_{e}^*(\mathrm d\eta)\otimes\tau^{[d]}_{e_1}(\tau^+_{e_2}-\tau_{e_2}^-)+\frac{1}2\pi_{e}^*(\mathrm d\eta)\otimes\tau^{[d]}_{e_2}(\tau^+_{e_1}-\tau_{e_1}^-),
\] 
where $e_1=(i,j)$, $e_2=(j,i)$, and $e$ is (improperly) the ``phantom'' short loop arising from the contraction of the two vertices $i$, $j$.
Here, we have used the obvious fact that $\tau_{e_1}^{[d]}\tau_{e_2}^{[d]}=-\tau_{e_2}^{[d]}\tau_{e_1}^{[d]}$.

The factor $1/2$ before the function $\eta$ on $\mathcal C_{1,1}^+$ (which we have tacitly omitted) is compatible with the fact that the pull-back of $\eta$ from $\mathcal C_{1,0,0}^+$ to $\mathcal C_{1,1}^+$ is precisely the rescaled function $\eta$ on $\mathcal C_{1,1}^+$.

Therefore, the same arguments as in the corresponding part of the proof of~\cite[Theorem 7.2]{CFFR} show that the right compensation for the contributions coming from the last three graphs in Figure 4 is given precisely by the superloop propagator $\omega^K_e$, which differ from the superloop propagator chosen in the proof of~\cite[Theorem 7.2]{CFFR} in its multidifferential operator part: the trick is to prove that we may rewrite the multidifferential operator parts of the contributions coming from the last three graphs in Figure 4 using the difference $\tau_e^+-\tau_e^-$, which is exactly the term appearing if we do the computations using the compactified configuration spaces $\mathcal C_{n,k,l}^+$ instead of $\mathcal C_{n,m}^+$. 
\end{proof}

\subsubsection{Biquantization as a consequence of Theorem~\ref{t-form}}\label{sss-2-3-4}
We consider now the particular situation $X=\mathfrak g^*$, for $\mathfrak g$ a finite-dimensional Lie algebra over $\mathbb K$, and for a given Lie subalgebra $\mathfrak h$ thereof, we set $U_1=X$ and $U_2=\mathfrak h^\perp$, the annihilator of $\mathfrak h$ in $\mathfrak g$.
We observe that, later on, we will consider $U_2$ to be the affine space $\lambda+\mathfrak h^\perp$, where $\lambda$ is a character of $\mathfrak h$: the results of the previous Subsubsection still hold true in this situation.

For the sake of explicit computations, we choose a complementary subspace of $\mathfrak h$ in $\mathfrak g$, {\em i.e.} we choose a subspace $\mathfrak p$ of $\mathfrak g$, such that $\mathfrak g=\mathfrak h\oplus\mathfrak p$.
We observe that, in general, $\mathfrak p$ is not $\mathfrak h$-invariant with respect to the restriction of the adjoint representation.
Still, in the case of symmetric pairs $(\mathfrak k,\mathfrak p)$, $\mathfrak p$ is a $\mathfrak k$-module.

We thus apply~\cite[Theorem 7.2]{CFFR} to this situation (we only observe that, in this framework, we do not consider completed algebras, as in~\cite{CFFR}: still, Theorem 7.2 holds true, the only difference is that we have to drop the property of the $L_\infty$-morphism to be an $L_\infty$-quasi-isomorphism): we may view the Poisson structure on $X$ as a Maurer--Cartan (shortly, form now on, MC) element of $T_\mathrm{poly}(X)$, and its image with respect to the $L_\infty$-morphism from~\cite[Theorem 7.2]{CFFR} is a MC element in the Hochschild cochain complex (with mixed component completed) of the $A_\infty$-category $\mathrm{Cat}_\infty(A,B,K)$, with objects $U_i$, $i=1,2$.

Using the previous prescriptions, we have
\[
A=\mathrm S(\mathfrak g),\ B=\mathrm S(\mathfrak p)\otimes\wedge \mathfrak h^*,\ K=\mathrm S(\mathfrak p);
\]
a bit improperly, we sometimes write $\mathfrak p=\mathfrak g/\mathfrak h$ (as it is an identification only of vector spaces, obviously not of $\mathfrak h$-modules).

Since $\mathfrak g$ is a Lie algebra, $X=\mathfrak g^*$ is a Poisson manifold with linear Poisson bivector $\pi$, and $U_1$ and $U_2$ are coisotropic submanifolds thereof.
The linear Poisson structure on $X$ determines a Maurer--Cartan element of $T_\mathrm{poly}(X)$: for a choice of a formal parameter $\hbar$, the image of $\hbar\pi$ with respect to the $L_\infty$-morphism $\mathfrak U$ from Theorem~\ref{t-form} is a Maurer--Cartan element $\mathcal U(\hbar\pi)$ in the (completed) Hochschild cochain complex of the $A_\infty$-category $\mathcal A$, which is a concept needing some unraveling.

The Hochschild cochain complex of $\mathcal A$ splits into three terms, namely the Hochschild cochain complex of $A$, the one of $B$ and a graded vector space which contains $A$, $B$ and $K$: general elements of the mixed term $\mathrm C^\bullet(A,B,K)$ are multilinear maps from $A^{\otimes k}\otimes K\otimes B^{\otimes l}$ to $K$.
From the general theory of Hochschild cochain complexes it is known that Maurer--Cartan elements of Hochschild cochain complexes correspond to $A_\infty$-structures on the underlying graded vector spaces: in our situation, a Maurer--Cartan element is precisely a structure of $A_\infty$-algebra on both $A$ and $B$, and a corresponding structure of $A_\infty$-$A$-$B$-bimodule on $K$, or, equivalently, to an $A_\infty$-structure on the category $\mathcal A$.

Now we consider $A_\hbar=A[\![\hbar]\!]$, and similarly for $B_\hbar$ and $K_\hbar$: it is clear that the structure of $A_\infty$-category on $\mathcal A$ extends to $\hbar$-linearly to $\mathcal A_\hbar$, whose objects are the same objects of $\mathcal A$, but whose morphism spaces are replaced by $A_\hbar$, $B_\hbar$ and $K_\hbar$ endowed with the $\hbar$-linearly extended $A_\infty$-structure $\mu$.
We have a natural Hochschild differential $\mathrm d_\mathrm H$ on the Hochschild cochain complex of $\mathcal A_\hbar$, given by the adjoint representation of $\mu$ with respect to the Gerstenhaber bracket. 
The element $\mathcal U(\hbar)$ satisfies the Maurer--Cartan equation 
\[
\mathrm d_\mathrm H\mathcal U(\hbar\pi)+\frac{1}2\left[\mathcal U(\hbar\pi),\mathcal U(\hbar\pi)\right]=\frac{1}2\left[\mu+\mathcal U(\hbar\pi),\mu+\mathcal U(\hbar\pi)\right]=0,
\]
{\em i.e.} $\mu+\mathcal U(\hbar\pi)$ is a Maurer--Cartan element for $\mathcal A_\hbar$, which deforms (with respect to the formal parameter $\hbar$) the ``classical'' $A_\infty$-structure on $\mathcal A$.

Since $A_\hbar$ is concentrated in degree $0$ by construction, and $\mu_A$ (the component of $\mu$ in the Hochschild cochain complex of $A$) is the obvious $\hbar$-linear commutative, associative product on $A_\hbar$, then $\mu_A+\mathcal U_A(\hbar\pi)$ is an associative product $\star_{A_\hbar}$ on $A_\hbar$, which deforms non-trivially $\mu_A$: $(A_\hbar,\star_{A_\hbar})$ is a deformation quantization of $(A,\mu_A)$ in the sense of~\cite{K}.

The image of $\pi$ in $T_\mathrm{poly}(A)$ with respect to the dg Lie algebra isomorphism $T_\mathrm{poly}(X)\cong T_\mathrm{poly}(A)$ (depending on a choice of $\mathfrak p$) is a Maurer--Cartan element in $T_\mathrm{poly}(A)$, which is a sum of a three polyvector fields, $\pi_0$, $\pi_1$ and $\pi_2$, where $\pi_i$ is an $i$-th polyvector field of polynomial degree $2-i$.
We observe that $A=\mathrm S(\mathfrak g)$ and $B=\mathrm S(\mathfrak p\oplus\mathfrak h^*[-1])$, thus it makes sense to speak about polynomial degree for elements of $A$ and $B$; $[\bullet]$ is the degree-shifting functor (hence, the polynomial grading of $B$ does not coincide with the internal grading coming from the functor $[-1]$), {\em e.g.} the internal degree of $\pi_i$ is $1$, $i=0,1,2$.
As a Maurer--Cartan element of $T_\mathrm{poly}(A)$, $\pi$ defines a $P_\infty$-structure on $B$, in other words, $\pi$ defines a Poisson algebra structure on $B$ up to homotopy: for example, $\pi_0$ is a homological vector field over $B$, whose cohomology identifies with the Chevalley--Eilenberg cohomology of the $\mathfrak h$-module $\mathrm S(\mathfrak g/\mathfrak h)$, which in turn inherits from $\pi_1$ (which is a bivector field of internal degree $1$) a structure of graded Poisson algebra.
We notice that, in degree $0$, this corresponds to the well-known fact that Poisson reduction endows the commutative algebra $\mathrm S(\mathfrak g/\mathfrak h)^\mathfrak h$ with a Poisson structure coming from the natural one on $A=\mathrm S(\mathfrak g)$.
%We will return to this issue in a subsequent paper~\cite{CFR}, where we deal with the extension to higher degree cohomology of the relative Duflo conjecture; we will touch the subject later on in the present paper in degree $0$ cohomology.
The Maurer--Cartan element $\mu_B+\mathcal U_B(\hbar\pi)$ is an $A_\infty$-structure on $B_\hbar$, deforming the obvious $A_\infty$-structure on $B$: thus, a $P_\infty$-algebra structure on $B$ produces an $A_\infty$-structure {\em via} the graded version of deformation quantization, see also~\cite{CF}.
We observe that the $A_\infty$-structure $\mu_B+\mathcal U_B(\hbar\pi)$ is the sum of (possibly) infinitely many components of different internal degree: in particular, the component of internal degree $2$ is an element of $B_\hbar$ of degree $2$, the curvature of the $A_\infty$-structure.
If it non-trivial, then we cannot talk about the cohomology of $A_\infty$-algebra $(B_\hbar,\mu_B+\mathcal U_B(\hbar\pi))$, and some problems may arise: luckily, in the present framework, the curvature vanishes, see {\em e.g.}~\cite{CFb,CT} and later on.
We finally observe that the term of order $1$ with respect to $\hbar$ of the $A_\infty$-structure on $B_\hbar$ is precisely the ($\hbar$-shifted) $P_\infty$-structure on $B_\hbar$: thus, if we select its vector field piece, we get the $\hbar$-shifted Chevalley--Eilenberg differential on $B_\hbar$.

Finally, the mixed component $\mu_K+\mathcal U_K(\hbar\pi)$ determines a deformation of the $A_\infty$-$A$-$B$-bimodule $\mu_K$ structure on $K$: we do not spend here much words, because we will deal with $\mu_K+\mathcal U_K(\hbar\pi)$ in the rest of the paper, at least in degree $0$.
We only observe that, through $\mu_K$, we may re-prove classical Koszul duality between $A$ and $B$ (both are graded quadratic algebras), and its deformation quantization permits to extend the Koszul duality to the deformed case $(A_\hbar,\star_\hbar)$ and $(B_\hbar,\mu_B+\mathcal U_B(\hbar\pi))$.
%We will deal with this issue regarding equivalences of categories of representations in a subsequent paper~\cite{FR}
%The deformation quantization of Koszul duality seems to be also central in a seemingly fruitful strategy to prove the relative Duflo conjecture~\cite{CFR}.

Biquantization as in~\cite{CT} is the specialization to degree $0$ of the data presented above.
In particular, $(A_\hbar,\star_{A_\hbar})$ is an $A_\infty$-algebra concentrated in degree $0$, hence its cohomology equals itself; the piece of $B_\hbar$ of degree $0$ equals $\mathrm S(\mathfrak p)\cong\mathrm S(\mathfrak g/\mathfrak h)$ endowed with a differential $\mathrm d_{B_\hbar}^0$ and with an associative product $\star_{B_\hbar}$ up to homotopy.
Finally, $K_\hbar$ is also concentrated in degree $0$, hence its cohomology with respect to $\mathrm d_{K_\hbar}^{0,0}$ (the $(0,0)$-component of $\mu_K+\mathcal U_K(\hbar\pi)$) equals itself, hence $K_\hbar$ becomes with respect to $d_{K_\hbar}^{1,0}=\star_L$ a left $(A_\hbar,\star_{A_\hbar})$- and with respect to $\mathrm d_{K_\hbar}^{0,1}=\star_R$ a right $(\mathrm H^0(B_\hbar),\star_{B_\hbar})$-module (the latter also because of the vanishing of the curvature of the $A_\infty$-structure $\mu_B+\mathcal U(\hbar\pi)$).

Later on, still in the framework of finite-dimensional Lie algebras and Lie subalgebras thereof, we will consider the more general framework, where both $A_\hbar$ and $B_\hbar$ are $A_\infty$-algebras with no curvature, and $K_\hbar$ is a graded $A_\infty$-$A_\hbar$-$B_\hbar$-bimodule, hence the $0$-th cohomologies $\mathrm H^0(A_\hbar)$, $\mathrm H^0(B_\hbar)$ become associative algebras and $\mathrm H^0(K_\hbar)$ is an $\mathrm H^0(A_\hbar)$-$\mathrm H^0(B_\hbar)$-bimodule.

\subsubsection{Symmetries of the $4$-colored propagators}\label{sss-2-3-5}
For later purposes, we now exhibit certain symmetries of the $2$-colored and $4$-colored propagators, which we now discuss in some detail.

The complex upper half-plane $\mathbb H^+$ has two obvious symmetries, namely the reflection with respect to the imaginary axis $i\mathbb R$, given by $z\overset{\sigma}\mapsto -\overline z$, and the inversion with respect to the unit half-circle, given by $z\overset{\tau}\mapsto 1/\overline z$: both maps extend to $\mathbb H^+\sqcup \mathbb R$, and they define two orientation-reversing involutions $\sigma$ and $\tau$ of it.
Equivalently, $\mathcal Q^{+,+}\sqcup i\mathbb R^+\sqcup\{0\}\sqcup\mathbb R^+$ admits two orientation-reversing involutions $\sigma$ and $\tau$, where $z\overset{\sigma}\mapsto i\overline{z}$ and $z\overset{\tau}\mapsto \frac{1}{\overline z}$.

It is not difficult to prove that $\sigma$ and $\tau$ descend both to involutions of $C_{n,m}^+$ and $C_{n,k,l}^+$, and that, using the same techniques as in the proof of Proposition~\ref{p-square}, Subsubsection~\ref{sss-2-1-3}, $\sigma$ and $\tau$ extend to involutions of the compactified configuration spaces $\mathcal C_{n,m}^+$ and $\mathcal C_{n,k,l}^+$.
We observe that $\sigma$ and $\tau$ are orientation-preserving, resp.\ -reversing, if and only if $n+m-1$ is even, resp.\ odd.

We then have the following technical Lemma about the behavior of the $2$-colored and $4$-colored propagators with respect to the action of $\sigma$ and $\tau$. 
\begin{Lem}\label{l-symm-4}
The $2$-colored and $4$-colored propagators behave as follows with respect to the involutions $\sigma$ and $\tau$ on the respective compactified configuration spaces $\mathcal C_{2,0}^+$ and $\mathcal C_{2,1}^+$:
\[
\begin{aligned}
\sigma^*(\omega^+)&=-\omega^+,\ & \sigma^*(\omega^-)&=-\omega^-,\ & \tau^*(\omega^+)&=-\omega^++2\pi_1^*(\mathrm d\eta),\ & \tau^*(\omega^-)&=-\omega^-+2\pi_2^*(\mathrm d\eta),\\
\sigma^*(\omega^{+,+})&=-\omega^{+,+}, & \sigma^*(\omega^{+,-})&=-\omega^{-,+}, & \sigma^*(\omega^{-,+})&=-\omega^{+,-}, & \sigma^*(\omega^{-,-})&=-\omega^{-,-},\\
\tau^*(\omega^{+,+})&=-\omega^{+,+}+2 \pi_1^*(\mathrm d\eta), & \tau^*(\omega^{+,-})&=-\omega^{+,-}, & \tau^*(\omega^{-,+})&=-\omega^{-,+}, & \tau^*(\omega^{-,-})&=-\omega^{-,-}+2 \pi_2^*(\mathrm d\eta),
\end{aligned}
\]
where now $\pi_i$, $i=1,2$, denotes the two natural projections from $\mathcal C_{2,0}^+$ onto $\mathcal C_{1,0}^+$ or from $\mathcal C_{2,1}^+$ to $\mathcal C_{1,1}^+$.
Similar formul\ae\ hold true for the $4$-colored propagators on $\mathcal C_{2,0,0}$, keeping in track a rescaling before the exact $1$-form $\eta$.
\end{Lem}

\subsubsection{Kontsevich's Vanishing Lemmata}\label{sss-2-3-6}
We now need a Vanishing Lemma for the $4$-colored propagators, reminiscent of the Vanishing Lemmata in~\cite[Subsubsubsection 7.3.3.1]{K}.
We observe that Kontsevich's Vanishing Lemmata in~\cite[Subsubsubsection 7.3.3.1]{K} are key ingredients in the proof of the globalization of its $L_\infty$-Formality-quasi-isomorphism: in this sense, the Vanishing Lemma we are going to state and prove here (the main application being for later computations regarding the generalized Harish-Chandra homomorphism) play also a central {\em r\^ole} in the globalization of the $2$-brane $L_\infty$-Formality-quasi-isomorphism of~\cite{CFFR}, but do not indulge here on this point, referring to upcoming work for more details.

We consider the three natural projections $\pi_{ij}$, $i\leq i<\leq 3$, from $\mathcal C_{3,0}^+$ onto $\mathcal C_{2,0}^+$, which smoothly extend the projections $[(z_1,z_2,z_3)]\overset{\pi_{ij}}\to [(z_i,z_j)]$ to the corresponding compactified configuration spaces.
The typical fiber of the projection $\pi_{ij}$ is $2$-dimensional, and it is not difficult to verify that it is a smooth manifold with corners (hence, it admits a natural stratification, whose description, at least in codimension $1$, will be made explicit later on).
We improperly denote by the same symbol the natural projection $\pi_{ij}$, $i\leq i<\leq 3$, from $\mathcal C_{3,1}^+$ onto $\mathcal C_{2,1}^+$, which this times extends the projection $[(z_1,z_2,z_3,x)]\overset{\pi_{ij}}\mapsto[(z_i,z_j,x)]$: again, its fiber is a smooth $2$-dimensional manifold with corners.

For any two smooth $1$-forms $\eta_i$, $i=1,2$, on $\mathcal C_{2,0}^+$ or $\mathcal C_{2,1}^+$, we define a smooth function on $\mathcal C_{2,0}^+$ or $\mathcal C_{2,1}^+$ {\em via} the integral
\begin{equation}\label{eq-int-van}
\Omega(\eta_1,\eta_2)=\pi_{13,*}\!\left(\pi_{12}^*(\eta_1)\wedge\pi_{23}^*(\eta_2)\right)=\int_{z_2\in\mathcal Q^{+,+}\smallsetminus\{z_1,z_3\}}\eta_1(z_1,z_2)\wedge\eta_2(z_2,z_3),
\end{equation}
where $\pi_{ij,*}$ denotes integration along the fiber of the projection $\pi_{ij}$.
\begin{Lem}\label{l-vanish-4}
The function $\Omega(\eta_1,\eta_2)$ vanishes, whenever $\eta_1=\eta_2$ is either one of the $2$-colored propagators or either one of the $4$-colored propagators.
\end{Lem}
\begin{proof}
The claim for the $2$-colored propagators $\omega^+$ and $\omega^-$ is precisely the content of the vanishing lemmata in~\cite[Subsubsubsection 7.3.3.1]{K}, to which we refer for a proof.
We observe that the proof below for the $4$-colored propagators applies with minor changes ({\em e.g.} the final argument involves the involution $\sigma$ and not $\tau$) applies to the statement for $2$-colored propagators.

For symmetry reasons, it suffices to prove the claim for the $4$-colored propagators $\omega^{+,+}$ and $\omega^{+,-}$.

We first prove the claim for $\eta_1=\eta_2=\omega^{+,+}$: the idea is to show first that $\Omega(\eta_1,\eta_2)$ is a constant function, and then to use Lemma~\ref{l-symm-4} to prove that, for a well-suited choice of arguments, $\Omega(\eta_1,\eta_2)$ equals minus itself.

We therefore compute the exterior derivative of $\Omega(\eta_1,\eta_2)$: we make use of generalized Stokes' Theorem, and the fact that $\eta_1=\eta_2$ is a closed $1$-form, yields
\[
\mathrm d\Omega(\eta_1,\eta_2)=\pi_{13,*}^\partial(\pi_{12}^*(\eta_1)\pi_{23}^*(\eta_2)),
\]
where $\pi_{13,*}^\partial$ denotes integration along the codimension $1$-boundary strata of the fiber of $\pi_{13}$: there are five such boundary strata, which correspond to $i)$ the point labeled by $z_2$ approaching either $\mathbb R$ on the left or on the right of the marked point on $\mathbb R$ or the marked point itself, or to $ii)$ the point labeled by $z_2$ approaching either the point labeled by $z_1$ or $z_2$.

In the three situations in $i)$, the corresponding contribution vanishes in view of Lemma~\ref{l-CF}, $iii)$ and $iv)$.

In both situations described in $ii)$, the boundary fibration is trivial, namely $\mathcal C_2\times \mathcal C_{2,1}^+$: an easy computation in local coordinates for $\mathcal C_{3,1}^+$ near the boundary strata in $i)$ shows that there are no orientation signs appearing, and we finally get, using Lemma~\ref{l-CF}, $i)$,
\[
\mathrm d\Omega(\eta_1,\eta_2)=\int_{\mathcal C_2^+}\mathrm d\varphi\ \omega^{+,+}+\int_{\mathcal C_2^+}\omega^{+,+}\mathrm d\varphi=0.
\]

Therefore, $\Omega(\eta_1,\eta_2)$ is a constant function on $\mathcal C_{2,1}^+$, whose value is completely determined by a choice of a point in $\mathcal C_{2,1}^+$, and a natural choice is $(i,2i,0)$.
Since $i\mathbb R$ is the fixed point set of $\sigma$, $\sigma$ preserves $\Omega(\eta_1,\eta_2)$, whence
\[
\Omega(\eta_1,\eta_2)=\sigma^*(\Omega(\eta_1,\eta_2)=-\pi_{13,*}(\sigma^*(\pi_{12}^*(\eta_1))\wedge \sigma^*(\pi_{23}^*(\eta_1))=-\Omega(\eta_1,\eta_2),
\]
where the minus sign in the second equality arises because $s$ is orientation-reversing on the first quadrant, while the third equality is a consequence of Lemma~\ref{l-symm-4}.

The very same arguments can be applied in the situation $\eta_1=\eta_2=\omega^{-,-}$.

We consider now the case $\eta_1=\eta_2=\omega^{+,-}$.
The computation of the exterior derivative of $\Omega(\eta_1,\eta_2)$ in this case is similar to the previous one: we only observe that the boundary condition for boundary strata of type $i)$ let appear a regular term $\mathrm d\eta$, whose contribution to integration is trivial.
The arguments for dealing with the other strata are, once again, a consequence of Lemma~\ref{l-CF}, $iii)$ and $iv)$.

Therefore, $\Omega(\eta_1,\eta_2)$ is uniquely determined by a given point in $\mathcal C_{2,1}^+$: quite differently from the previous case, we choose a pair of points lying on the unit circle.
Recalling that the unit circle is the fixed point {\em locus} of the involution $\tau$, we get in this case
\[
\Omega(\eta_1,\eta_2)=\tau^*(\Omega(\eta_1,\eta_2)=-\pi_{13,*}(\tau^*(\pi_{12}^*(\eta_1))\wedge \tau^*(\pi_{23}^*(\eta_1))=-\Omega(\eta_1,\eta_2),
\]
by the very same arguments as in the previous case, because $t$ is orientation-reversing and because of Lemma~\ref{l-symm-4}.
\end{proof}
As an application of Lemma~\ref{l-vanish-4}, we briefly sketch the vanishing of the curvature of the $A_\infty$-algebra $B_\hbar$, with the previously introduced notations.
As already mentioned, the curvature of $B_\hbar$ is the piece of degree $2$ in $B_\hbar$ of the Maurer--Cartan element $\mu_B+\mathcal U_B(\hbar\pi)$: more explicitly, the curvature is given by the formal power series 
\[
\mathcal U(\hbar\pi)_0=\sum_{n\geq 1}\frac{1}{n!}\mathcal U_B^n(\underset{n}{\underbrace{\hbar\pi|\cdots|\hbar\pi}})=\sum_{n\geq 1}\frac{1}{n!}\sum_{\Gamma\in\mathcal{G}_{n,0}}\mathcal{O}_{\Gamma}^B(\underset{n}{\underbrace{\hbar\pi|\cdots|\hbar\pi}})=\sum_{n\geq 1}\frac{1}{n!}\sum_{\Gamma\in\mathcal{G}_{n,0}}\mu_{n}^B\left(\int_{\mathcal C_{n,0}^+}\prod_{e\in E(\Gamma)}\omega^B_e(\underset{n}{\underbrace{\hbar\pi|\cdots|\hbar\pi}})\right).
\]
For an admissible graph $\Gamma$ in $\mathcal G_{n,0}$, the rightmost integral is non-trivial only if the degree of the integrand equals $2n-2$.
To each vertex of the first type of $\Gamma$ is associated a copy of the linear Poisson bivector $\hbar\pi$, hence from each vertex depart exactly two arrows: each arrow, by definition, corresponds to a derivation and a contraction.
In particular, the differential operator $\mathcal O_\Gamma^B$ has degree $2n-2$: since $\hbar\pi$ is linear, the polynomial degree of the object on the rightmost part of the previous chain of equalities is $n-(2n-2)=-n+2$, whence $1\leq n\leq 2$.
When $n=1$, $\mu_1^B(\hbar\pi)$ vanishes because of the coisotropy of $U_2$.
For $n=2$, there is only one possible admissible graph $\Gamma$ of type $(2,0)$, namely the loop graph connecting the two vertices of the first type: the corresponding operator $\mathcal O_\Gamma^B$ vanishes because of Lemma~\ref{l-vanish-4}.

\subsection{Quantum reduction algebras}\label{ss-2-4}
The present Subsection presents the results of~\cite[Section 2]{CT} using a slightly different perspective, coherent with the approach to biquantization we have introduced in the previous Subsection: however, we think it useful to review many results mainly because of the notation, which will be then used extensively in the rest of the paper.

Thus, we will mostly concentrate on the dg vector space $B_\hbar$, where we denote by $\mathrm d_{B_\hbar}$ its differential ({\em i.e.} the piece of degree $1$ of its $A_\infty$-structure): we have already observed that $\mathrm d_{B_\hbar}=\hbar\mathrm d_\mathrm{CE}+\mathcal O(\hbar^2)$, where $\mathrm d_\mathrm{CE}=\pi_0$ is the Chevalley--Eilenberg differential on $B_\hbar$.
The admissible graphs $\Gamma$ appearing in $\mathrm d_{B_\hbar}$ are of type $(n,1)$, possibly with multiple edges and admitting edges departing from $\mathbb R$: each edge $e$ of an element $\Gamma$ of $\mathcal C_{n,1}$ is the sum of two colored edges, namely $e^+$ and $e^-$ according to Formula~\eqref{eq-B-form}, Subsubsection~\ref{sss-2-3-2}.
The properties of the superpropagator~\eqref{eq-B-form} imply that an edge $e$ of $\Gamma$ arriving to the only vertex of the second type has color $+$, while an edge departing from it has color $-$.

One of the main tools we will use throughout the present Subsection is Lemma~\ref{l-vanish-4}, Subsubsection~\ref{sss-2-3-6}.
Namely, we assume $v$ is a vertex of the first type of $\Gamma$ of type $(n,1)$ with two edges at it of the form $e_1=(\bullet_1,v)$ and $e_2=(v,\bullet_2)$, where $\bullet_i$, $i=1,2$, denotes some other vertex (notice that now we allow $\bullet_1=\bullet_2$): then the configuration at $v$ is either $(e_1^+,e_1^-)$ or $(e_1^-,e_2^+)$.
Of course, since to each edge of the first type of $\Gamma$ is associated a copy of the linear Poisson bivector $\hbar\pi$, we assume that from $v$ as above departs a third edge, which does not join any other vertex of $\Gamma$: this edge has color $-$.
This ``phantom'' edge is the ``edge to $\infty$'', using the terminology of~\cite{CT}: dimensional arguments imply that each admissible graph $\Gamma$ of type $(n,1)$ admit precisely one vertex of the first type with a phantom edge.

\subsubsection{Symmetric pairs and, more generally, Lie subalgebras of trivial extension class}\label{sss-2-4-1}
We consider here the case of a symmetric pair $\mathfrak g=\mathfrak k\oplus\mathfrak p$, or, more generally, of a Lie subalgebra $\mathfrak h\subseteq\mathfrak g$ admitting an $\mathfrak h$-invariant complementary subspace $\mathfrak p$: we observe that this is equivalent to the triviality of the extension class $\alpha$ of the short exact sequence of $\mathfrak h$-modules $\mathfrak h\hookrightarrow\mathfrak g\twoheadrightarrow\mathfrak g/\mathfrak h$.

By definition, a symmetric pair is a pair $(\mathfrak g,\sigma)$, where $\mathfrak g$ is a finite-dimensional Lie algebra over $\mathbb K$ and $\sigma$ is an involutive Lie algebra automorphism: thus, $\mathfrak g=\mathfrak k\oplus\mathfrak p$ is the direct sum of the $+1$-eigenspace $\mathfrak k$ and the $-1$-eigenspace $\mathfrak p$ of $\sigma$.
In particular, we have the Cartan relations
\begin{equation}\label{eq-cartan}
[\mathfrak k,\mathfrak k]\subseteq \mathfrak k,\ [\mathfrak k,\mathfrak p]\subseteq \mathfrak k,\ [\mathfrak p,\mathfrak p]\subseteq \mathfrak k. 
\end{equation}
In particular, $\mathfrak k$ is a Lie subalgebra of $\mathfrak g$ and $\mathfrak p$ is a $\mathfrak k$-module.
The graded algebra $B$ in the case of a symmetric pair equals $B=\mathrm S(\mathfrak k)\otimes\wedge^\bullet(\mathfrak p)$.

Of course, if the above extension class $\alpha$ vanishes, then $\mathfrak g$ admits simply a decomposition $\mathfrak g=\mathfrak h\oplus\mathfrak p$, with relations $[\mathfrak h,\mathfrak h]\subseteq\mathfrak h$ and $[\mathfrak h,\mathfrak p]\subseteq\mathfrak p$.

The claim is that in the case of a symmetric pair $(\mathfrak g,\sigma)$ or of a pair $(\mathfrak g,\mathfrak h)$ with trivial extension class the quantized differential $\mathrm d_{B_\hbar}$ equals simply the ($\hbar$-shifted) Chevalley--Eilenberg differential $\hbar\mathrm d_\mathrm{CE}$ on $B_\hbar=\mathrm C^\bullet(\mathfrak h,\mathrm S(\mathfrak g/\mathfrak h))[\![\hbar]\!]$.

Namely, we consider an admissible graph $\Gamma$ of type $(n,1)$, $n\geq 2$, appearing in Formula~\eqref{eq-B-mor}, Subsubsection~\ref{sss-2-3-3}, and we know from the above considerations that $\Gamma$ possesses a vertex $v$ of the first type with a phantom edge $e_\mathrm{gh}$ and two edges $e_1=(\bullet_1,v)$ and $e_2=(v,\bullet_2)$.
The configuration at the vertex is either $(e_1^+,e_2^-,e_\mathrm{gh}^-)$ or $(e_1^-,e_2^+,e_\mathrm{gh}^-)$: in Lie algebraic terms, these two configurations correspond to $[\mathfrak h,\mathfrak p]\subseteq\mathfrak h$ and $[\mathfrak h,\mathfrak h]\subseteq\mathfrak p$ respectively, which is a contradiction to the Cartan relations, for $\mathfrak g$ a symmetric pair, or to the fact that $\mathfrak k$ has a $\mathfrak k$-invariant complement.
This implies that only the contributions of order $1$ with respect to $\hbar$ matter, whence the claim.

\subsubsection{Cohomology of degree $0$}\label{sss-2-4-2}
The content of the present Subsubsection presents some arguments for dealing with the classification of admissible graphs appearing in the computation of the differential $\mathrm d_{B_\hbar}$, which will also appear in other contexts related to biquantization.

We consider here the case of a Lie subalgebra $\mathfrak h$ of $\mathfrak g$ with no assumptions on the extension class of $\mathfrak h\subseteq \mathfrak g$: thus, we only assume to have picked out some complementary subspace $\mathfrak p$ for explicit computations.
\begin{Prop}\label{p-quant-0}
The admissible graphs of type of $(n,1)$, $n\geq 1$, appearing in the restriction of the differential $\mathrm d_{B_\hbar}$ to $\mathrm S(\mathfrak p)$ are either of type Bernoulli ({\em i.e.} connected graphs with one root and one phantom edge), or of type wheel ({\em i.e.}  connected graphs whose edges between the vertices of the first type form an oriented loop with one phantom edge), or of mixed type ({\em i.e.} a Bernoulli-type graphs attached to a wheel-type graph), see also Figure 5.
\end{Prop}
\begin{proof}
We consider, for $n\geq 1$, an admissible graph of type $(n,1)$, and we denote by $p$ the number of edges of $\Gamma$ arriving at the vertex of the second type.

Degree reasons imply that such a graph admits one phantom edge, hence the actual edges connecting vertices of $\Gamma$ are $2n-1$.
Dimensional reasons imply also that $p\geq 1$: namely, if $p$ were $0$, since $n\geq 1$ by assumption, there would be a $1$-dimensional submanifold of $\mathcal C_{n,1}^+$ over which there is nothing to integrate (a subset of $\mathbb R$), hence integration would yield $0$.

The admissible graph $\Gamma$ is connected in the sense that, if we remove from it the edges to the vertex of the second type, we obtain a connected graph in the strict sense of the world, as $p\geq 1$.
The connectedness of $\Gamma$ is also a consequence of dimensional reasons: if $\Gamma=\Gamma_1\sqcup\Gamma_2$, then either $\Gamma_1$ or $\Gamma_2$ would have a phantom edge. 
W.l.o.g.\ we assume $\Gamma_1$ has a phantom edge, hence $\Gamma_2$ does not, which means that all its edges provide differential forms to be integrated: if $n_i$ is the number of vertices of the first type of $\Gamma_i$, $i=1,2$, then the integral over $\mathcal C_{n_2,1}^+$ of the differential form associated to $\Gamma_2$ vanishes, as its degree is $2 n_2$, while the dimension of $\mathcal C_{n_2,1}^+$ is $2n_2-1$.

The boundary conditions for the superpropagator~\eqref{eq-B-form} imply immediately that edges arriving at the vertex of second type are all colored by $+$, whence $p\leq (n-1)+1=n$: in fact, from the vertex of the first type from which departs the phantom edge departs another edge, which may or may not arrive at the vertex of the second type. 

On the other hand, there are $2n-1-p$ edges from vertices of the second type arriving at (distinct) vertices of the first type: this implies that the polynomial degree of the differential operator on $B_\hbar$ associated to $\Gamma$ is $n-(2n-1-p)=-n+1+p\geq 0$, whence $p\geq n-1$: it follows then immediately $p=n-1$ or $p=n$. 

We first consider the case $p=n$: from every vertex of the first type of $\Gamma$ departs one edge to the vertex of the second type, $\Gamma$ is connected, has a phantom edge and there is a single vertex of the first type, which is the final point of none of its edges.
Such an admissible graph is obviously of Bernoulli type.

Then, assume $p=n-1$: the only vertex from which does not depart an edge to the vertex of the second type may or may not be the vertex from which departs the phantom edge.
In the first case, the connectedness of $\Gamma$ implies that it is of wheel-type, while in the second case, it must be a wheel type, from which departs an edge hitting the root of a Bernoulli-type graph.
Here, the root of a Bernoulli-type graph is the only vertex of the second type, which is the final point of none of its edges.
\end{proof}
Here is a pictorial representation of the three types of graphs appearing in $\mathrm d_{B_\hbar}$ according to the previous Proposition:
\bigskip
\begin{center}
\resizebox{0.6 \textwidth}{!}{\input{d-quant-gr.pstex_t}}\\
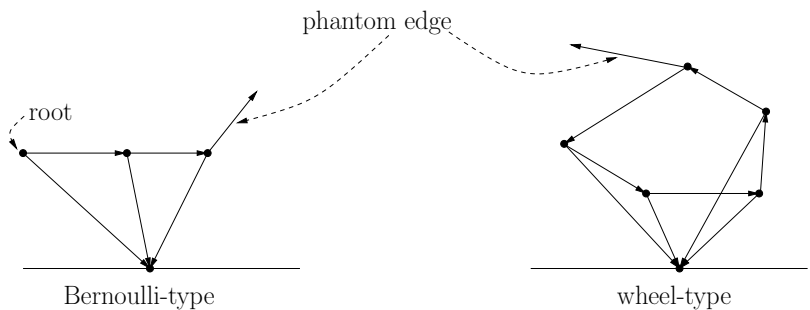
\text{Figure 5 - Bernoulli-type and wheel-type graphs} \\
\end{center}
\bigskip

Although for most of the computations in this framework we do not really need it, we want to understand the differential $\mathrm d_{B_\hbar}$ on the whole of $B_\hbar$: as the next proposition shows, the results of Proposition~\ref{p-quant-0} with a slight addition suffice.
\begin{Prop}\label{p-quant-d}
The admissible graphs of type $(n,1)$ appearing in $\mathrm d_{B_\hbar}$ are either the admissible graphs of Proposition~\ref{p-quant-0} or connected graphs which are brackets of two Bernoulli-type graphs ({\em i.e.} there is an edge departing from the vertex of the second type to a vertex of the first type, whose two outgoing edges arrive at the roots of two Bernoulli-type graphs, see Figure 6).
\end{Prop}
\begin{proof}
We first observe that, since we are considering $\mathrm d_{B_\hbar}$ on the whole of $B_\hbar$, admissible graphs $\Gamma$ can have edges departing from the vertex of the second type.

We first assume that that all edges departing from the vertex of the second type are phantom edges: in this case, we are reduced to the very same analysis as in the proof of Proposition~\ref{p-quant-0}.

We now assume that the admissible graph $\Gamma$ of type $(n,1)$ possesses $k$ edges departing from the vertex of the second type and arriving to $k$ vertices of the first type (these vertices are distinct because of the linearity of the Poisson bivector $\hbar\pi$).
Because of degree reasons, $\Gamma$ has $k+1$ phantom edges departing from vertices of the first type.
Dimensional reasons imply further that no vertex of the first type may have two phantom edges.
We denote by $p$ the number of edges departing from vertices of the first type and arriving to the vertex of the second type: the very same arguments as in the proof of Proposition~\ref{p-quant-0} imply that $p\geq n$.
On the other hand, the polynomial degree of the differential operator associated to $\Gamma$ equals $(n-k)-(2n-(k+1)-p)=-n+1+p\geq 0$, whence either $p=n-1$ or $p=n$.

We first consider the case $p=n$: in this situation, from every vertex of the first type departs exactly one edge to the vertex of the second type.
Thus, $\Gamma$ is the disjoint union of exactly one Bernoulli-type graph and of either wheel-like graphs or Bernoulli-wheel-type graphs or Bernoulli-type graphs, whose root is the endpoint of an edge departing from the vertex of the second type (at least such a graph appears here, as $k\geq 1$ by assumption).
Since $\Gamma$ is a disjoint union of subgraphs, the corresponding integral factors into integrals over compactified configuration spaces of the form $\mathcal C_{m,1}^+$, for $1\leq m\leq n$.
For a subgraph of the last type as above, the corresponding integral vanishes, because the degree of the integrand is $2m$, while the dimension of the fiber is $2m-1$.

We consider now the case $p=n-1$: the corresponding differential operator is translation-invariant.
If $\Gamma$ is the disjoint union of connected subgraphs, at least one of which is a Bernoulli-type graph, whose root is the final point of an edge departing from the vertex of the second type, the last argument in the previous paragraph yields triviality of the corresponding differential operator.
A subgraph like the one we have analyzed always appear, if there is a vertex of the first type of $\Gamma$, from which departs one edge to the vertex of the second type and which is the endpoint of an edge issued from the said vertex.
As $p=n-1$, there can be exactly one edge departing from the vertex of the second type and arriving at the only vertex of the first type, from which no edge depart to $\mathbb R$: these two edges meet two distinct roots (again, because $p=n-1$) of Bernoulli-type graphs.
Therefore, $\Gamma$ is the disjoint union of subgraphs as in Proposition~\ref{p-quant-0} and of a connected subgraph of the said type.
\end{proof}
Here is a pictorial representation of the new type of connected graphs appearing in $\mathrm d_{B_\hbar}$ in higher degrees according to the previous Proposition:
\bigskip
\begin{center}
\resizebox{0.35 \textwidth}{!}{\input{d-quant-gen.pstex_t}}\\
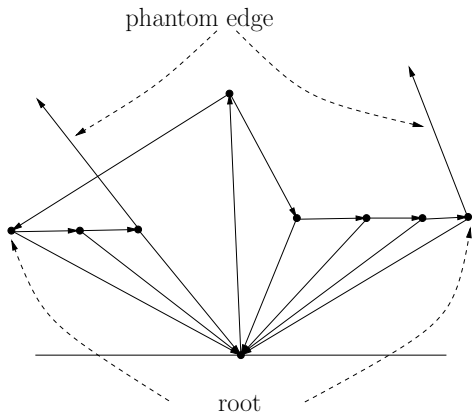
\text{Figure 6 - A connected graph for $\mathrm d_{B_\hbar}$ appearing in higher degrees} \\
\end{center}
\bigskip

We consider now a general admissible graph $\Gamma$ of type $(n,1)$ appearing in $\mathrm d_{B_\hbar}$: we may consider the involution $\sigma$ of $\mathcal C_{n,1}^+$ from Subsubsection~\ref{sss-2-3-5}.
It preserves, resp.\ reverses, orientation if $n$ is even, resp.\ odd; the pull-back of integrand in Formula~\eqref{eq-c-B} with respect to $\sigma$ equals itself up to a global $-1$-sign, hence if $n$ is even,  the contributions to $\mathrm d_{B_\hbar}$ from admissible graphs of type $(n,1)$ with $n$ even are trivial.
We may therefore write $\mathrm d_{B_\hbar}$ as
\[
\mathrm d_{B_\hbar}=\hbar\mathrm d_\mathrm{CE}+\hbar^3 \mathrm d_3+\hbar^5 \mathrm d_5+\cdots,
\]
where only odd powers of $\hbar$ appear.
Therefore, if $f$ is a general element of $B_\hbar$, then $f$ is $\mathrm d_{B_\hbar}$-closed if and only if both its even and odd part with respect to $\hbar$ are $\mathrm d_{B_\hbar}$-closed.

\subsubsection{The generalized Iwasawa decomposition}\label{sss-2-4-3}
In this Subsubsection we review in some detail the generalized Iwasawa decomposition of a symmetric pair $\mathfrak g=\mathfrak k\oplus\mathfrak p$: we want to stress that we do not present here any new results, but simply need to fix notation and conventions in view of later applications, namely the Harish-Chandra homomorphism in the framework of deformation quantization.
A detailed discussion of the Iwasawa decomposition for semisimple symmetric pairs can be found in~\cite[Chapter 1, Section 13]{Dix}; the main reference to the generalized Iwasawa decomposition for a general symmetric pair is~\cite[Subsections 1.1, 1.2 and 1.3]{T1}.

We consider a general symmetric pair $(\mathfrak g,\sigma)$ with Cartan decomposition $\mathfrak g=\mathfrak k\oplus\mathfrak p$.
For an element $\xi$ of $\mathfrak k^\perp$, we set $\mathfrak g(\xi)=\left\{x\in\mathfrak g:\ \mathrm{ad}^*(x)(\xi)=0\right\}$: it is a Lie subalgebra of $\mathfrak g$ and the Lie algebra involution $\sigma$ restricts to $\mathfrak g(\xi)$, whose Cartan decomposition we denote by $\mathfrak g(\xi)=\mathfrak k(\xi)\oplus\mathfrak p(\xi)$.

An element $\xi$ of $\mathfrak k^\perp$ is said to be regular, if the dimension of $\mathfrak g(\xi)$ is minimal among the subalgebras $\mathfrak g(\eta)$, $\eta$ in $\mathfrak k^\perp$.
The set of regular elements of $\mathfrak k^\perp$ is a Zarisky-open subset of $\mathfrak k^\perp=\mathfrak p^*$, and for every regular element $\xi$ of $\mathfrak k^\perp$, $[\mathfrak k(\xi),\mathfrak p(\xi)]=0$, see~\cite[Chapter 1.11]{Dix} for more details on regular linear functionals on a Lie algebra and~\cite[Definition 1.1.2.1 and Lemma 1.2.2.1]{T1} for a similar discussion in the present situation.

For $\xi$ in $\mathfrak k^\perp$ regular as above, we set $\mathfrak a(\xi)=\mathfrak p(\xi)\oplus[\mathfrak p(\xi),\mathfrak p(\xi)]$: the Cartan relations together with the previous remark imply that it is a nilpotent Lie subalgebra of $\mathfrak g$.
We denote by $\mathfrak s_\xi$ a maximal torus of $\mathfrak a(\xi)$, which is additionally preserved by $\sigma$: according to~\cite[Subsubsection 1.2.3]{T1}, there exists only one maximal torus $\mathfrak s_\xi$ of $\mathfrak a(\xi)$, which is exactly the set of semisimple elements of $\mathfrak p(\xi)$.
It is moreover central in $\mathfrak g(\xi)$.
The regular element $\xi$ is said to be generic (or, following the terminology of~\cite[Definition 1.2.3.1]{T1}, very regular), if the dimension of $\mathfrak s_\xi$ is maximal among all $\mathfrak s_\eta$, for $\eta$ regular.

According to~\cite[Lemma 1.3.1.1]{T1}, for $\xi$ very regular in $\mathfrak k^\perp$, $\mathfrak g$ admits a root decomposition $\mathfrak g=\mathfrak g_0\oplus\bigoplus_{\alpha\in\Delta}\mathfrak g_\alpha$, where $\mathfrak g_0$ is the centralizer of $\mathfrak s_\xi$ in $\mathfrak g$, and we set $\mathfrak g_\alpha=\left\{x\in\mathfrak g:\ \mathrm{ad}(t)(x)=\alpha(t)x,\ t\in\mathfrak s_\xi\right\}$, for $\alpha$ in $\mathfrak s_\xi^*$.
An element $\alpha$ of $\mathfrak s_\xi$ is said to be a root, if $\mathfrak g_\alpha\neq\{0\}$; the set of all roots is denoted by $\Delta$.
Moreover, $\mathfrak g_0$ is $\sigma$-stable and inherits the structure of symmetric pair, whence $\mathfrak g_0=\mathfrak k_0\oplus\mathfrak p_0$ its Cartan decomposition.
Finally, for a root $\alpha$ in $\Delta$, $\sigma(\mathfrak g_\alpha)=\mathfrak g_{-\alpha}$, whence $\mathfrak g$ admits a triangular decomposition $\mathfrak g=\mathfrak n_-\oplus\mathfrak g_0\oplus\mathfrak n_+$, where $\mathfrak n_\pm$ denotes the direct sum of all root spaces associated to positive/negative roots.

From the triangular decomposition follows the generalized Iwasawa decomposition $\mathfrak g=\mathfrak k\oplus\mathfrak p_0\oplus \mathfrak n_+$: in fact, see also~\cite[Proposition 1.13.11]{Dix}, $\mathfrak g=\mathfrak k+\mathfrak p_0+\mathfrak n_+$ directly from the previous triangular decomposition.
Namely, a general element $x$ of $\mathfrak g$ can be (uniquely) written as $x=x_-+k_0+p_0+x_+$, with $x_\pm$ in $\mathfrak n_\pm$, $k_0$ in $\mathfrak k_0$ and $p_0$ in $\mathfrak p_0$.
Then, $x_-=(x_-+\sigma(x_-))-\sigma(x_-)$, and obviously $x_-+\sigma(x_-)$ belongs to $\mathfrak k$, while $\sigma(x_-)$ is in $\mathfrak n_+$, whence the first claim follows.
Now, we assume $k+p_0+x_+=0$: then, $\sigma(k+p_0+x_+)=k-p_0+\sigma(x_+)=0$, whence $2p_0+x_+-\sigma(x_+)=0$.
From the triangular decomposition of $\mathfrak g$ follows automatically $p_0=x_+=0$, and thus also $k=0$.
The relations for the triangular decomposition of $\mathfrak g$ imply that $\mathfrak k_0\oplus \mathfrak n_+$ is a Lie subalgebra of $\mathfrak g$, and thus is $(\mathfrak k_0\oplus\mathfrak n_+)^\perp$ a coisotropic submanifold of $\mathfrak g^*$.
\begin{Prop}\label{p-iwasawa}
For a general symmetric pair $(\mathfrak g,\xi)$ and a very regular element $\xi$ of $\mathfrak k^\perp$, the reduction space $\mathrm H^0(B_\hbar)$ for the coisotropic submanifold $(\mathfrak k_0\oplus\mathfrak n_+)^\perp$ of $\mathfrak g^*$ equals $\mathrm S(\mathfrak p_0)^{\mathfrak k_0}[\![\hbar]\!]$. 
\end{Prop}
\begin{proof}
We perform a slight change in the proof, namely we consider the triangular decomposition $\mathfrak g=\mathfrak n_-\oplus\mathfrak k_0\oplus\mathfrak p_0\oplus\mathfrak n_+$: we observe that $\mathfrak n_-\oplus\mathfrak p_0$ is, in general, not a module for the Lie subalgebra $\mathfrak k_0\oplus\mathfrak n_+$.
According to previous arguments, the differential $\mathrm d_{B_\hbar}$ can be written as a formal power series $\hbar\mathrm d_\mathrm{CE}+\hbar^3\mathrm d_3+\cdots$, where only odd powers of $\hbar$ appear, and a general element $f$ satisfying $\mathrm d_{B_\hbar}(f)=0$ can be written as $f=f_0+\hbar^2 f_2+\cdots$, where only even powers of $\hbar$ appear.
In particular, $f_0$ is $\mathrm d_\mathrm{CE}$-closed: according to~\cite[???]{T1}, $f_0$ belongs then to $\mathrm S(\mathfrak p_0)^{\mathfrak k_0}$.

Therefore, to prove the claim, it suffices to prove that the operators $\mathrm d_{2n+1}$, $n\geq 1$ act trivially on $\mathrm S(\mathfrak p_0)^{\mathfrak k_0}$.
The differential operator $\mathrm d_{2n+1}$ is the sum of differential operators associated to admissible graphs in $\mathcal G_{2n+1,1}$, see Subsubsection~\ref{sss-2-3-1}.
In view of Proposition~\ref{p-quant-0}, 
an admissible graph $\Gamma$ yielding a (possibly) non-trivial contribution to $\mathrm d_{2n+1}$ is either of type Bernoulli or of type wheel or of mixed type.

We first consider $\Gamma$ of type Bernoulli in $\mathcal G_{2n+1,1}$: by Proposition~\ref{p-quant-0}, $\Gamma$ has a root ({\em i.e.} a vertex of the first type with no incoming edges), a phantom edge ({\em i.e.} an edge with no final point) and from all vertices of the first type there is exactly on outgoing edge to the unique vertex of the first type.
To every edge $e$ of $\Gamma$ corresponds a superpropagator, which in turn yields a coloring $e=e^++e^-$, where the color $\pm$ corresponds to linear coordinates on $\mathfrak n_-\oplus\mathfrak p_0$ and $\mathfrak k_0\oplus\mathfrak n_+$ respectively.
If the vertex of the second type is associated to $f_0$ in $\mathrm S(\mathfrak p_0)^{\mathfrak k_0}$, all edges pointing to the said vertex are colored by $+$, as they correspond to derivative with respect to $\mathfrak p_0$.
The phantom edge is obviously colored by $-$, as it corresponds to an element of $\mathfrak k_0^*\oplus\mathfrak n^*_+$.
We consider the vertex $v$ of the first type from which departs the phantom edge: the edge departing from $v$ to the only vertex of the second type is colored by $+$, thus in view of Lemma~\ref{l-vanish-4}, the only incoming edge to $v$ with initial point a distinct vertex of the first type is colored by $-$.
Therefore, all edges connecting two distinct vertices of the first type are colored by $-$: to be even more precise, they correspond to derivatives and contractions with respect to $\mathfrak n_+$.
Finally, the root of $\Gamma$ carries an element of $\mathfrak n_+$, which vanishes upon restriction.

We then consider an admissible graph $\Gamma$ of type wheel.
Therefore, $\Gamma$ has a phantom edge with initial point $v$, a vertex of the first type, and no root, and from each vertex of the first type different from $v$ depart exactly one edge to the only vertex of the second type. 
The phantom edge is colored by $-$, and it corresponds to an element of $\mathfrak k_0^*\oplus\mathfrak n_+^*$; as the differential operator corresponding to $\Gamma$ acts on $\mathrm S(\mathfrak p_0)^{\mathfrak k_0}$, all edges pointing to the only vertex of the second type are colored by $+$.
If we first assume that the phantom edge of $\Gamma$ corresponds to an element of $\mathfrak k_0^*$, the relations $[\mathfrak k_0,\mathfrak n_\pm]\subseteq \mathfrak n_\pm$ and $[\mathfrak k_0,\mathfrak p_0]\subseteq \mathfrak p_0$ imply that the vertex $v$ presents automatically a configuration of the form $(e_1,e_2,e_\mathrm{gh})$, $e_1=(\bullet,v)$, $e_2=(v,\bullet)$, with the same coloring, {\em i.e.} $(e_1^+,e_2^+,e_\mathrm{gh}^-)$ or $(e_1^-,e_2^-,e_\mathrm{gh}^-)$, thus the corresponding integral weight vanishes because of Lemma~\ref{l-vanish-4}.
If the phantom edge of $\Gamma$ corresponds to an element of $\mathfrak n_+^*$, we first assume that the edge $e_2$ departing from $v$ corresponds to derivation and contraction with respect to $\mathfrak k_0$ or $\mathfrak p_0$ ({\em i.e.} the color of $e_2$ is $-$ and $+$ respectively): the relations $[\mathfrak k_0,\mathfrak n_+]\subseteq \mathfrak n_+$ and $[\mathfrak p_0,\mathfrak n_+]\subseteq \mathfrak n_+$, together with the Cartan relations for the small symmetric pair $\mathfrak g_0$ imply that the edge $e_1$ carries simultaneously derivation and contraction with respect to $\mathfrak n_+$ and $\mathfrak k_0$ or $\mathfrak p_0$, which is impossible.
If $e_2$ corresponds to either $\mathfrak n_-$ or $\mathfrak n_-$, the relation $[\mathfrak n_\pm,\mathfrak p_0]\subseteq \mathfrak n_\pm$ implies that to $v$ corresponds a configuration of the form either $(e_1^+,e_2^+,e_\mathrm{gh}^-)$ or $(e_1^+,e_2^+,e_\mathrm{gh}^-)$, which implies triviality of the corresponding integral weight. 

We finally consider an admissible graph $\Gamma$ in $\mathcal G_{2n+1,1}$ of mixed type Bernoulli-wheel.
In this case, $\Gamma$ has a unique vertex $v$ of the first type from which departs the phantom edge, and there is a vertex $w$, from which departs an edge from a wheel-like graph $\Gamma_1$ to the root of a Bernoulli-like graph $\Gamma_2$.
The coloring of $\Gamma_2$ can be deduced {\em via} the same arguments used previously for a Bernoulli type graph: in particular, the color of the edge with the root of $\Gamma_2$ as final point is $-$: more precisely, it corresponds to derivation and contraction with respect to $\mathfrak n_+^*$.
The Cartan relations for the small symmetric pair $\mathfrak g_0$ together with $[\mathfrak k_0,\mathfrak n_\pm]\subseteq \mathfrak n_\pm$ and $[\mathfrak p_0,\mathfrak n_\pm]\subseteq \mathfrak \mathfrak n_\pm$ imply that the internal edges of the wheel-like graph $\Gamma_1$ have the same color, either $+$ or $-$; the corresponding edges correspond to derivation and contraction with respect to either $\mathfrak n_-$ and $\mathfrak n_+$.
The corresponding differential operators are of the form
\[
\text{either}\ \mathrm{tr}_{\mathfrak n_-}\!\left(\mathrm{ad}(X_1)\cdots\mathrm{ad}(X_p)\mathrm{ad}(Y)\right)\ \text{or}\ \mathrm{tr}_{\mathfrak n_+}\!\left(\mathrm{ad}(X_1)\cdots\mathrm{ad}(X_p)\mathrm{ad}(Y)\right),
\]
where $X_i$ are general elements of $\mathfrak p_0$ and $Y$ is an element either of $\mathfrak n_-$ or $\mathfrak n_+$.
Such operators are clearly trivial due to the nilpotence of both $\mathfrak n_\pm$.
\end{proof}

\subsubsection{Polarizations}\label{sss-2-4-4}
This short Subsubsection also serves the purpose of fixing notation and conventions for certain issues, which will be dealt later on by means of deformation quantization.

For more details on polarizations of Lie algebras, we refer to~\cite[Chapter 1.12]{Dix}.
A general element $\xi$ in $\mathfrak g^*$ defines a skew-symmetric bilinear form on a finite-dimensional Lie algebra $\mathfrak g$ {\em via} the assignment $B_\xi(x_1,x_2)=\langle\xi,[x_1,x_2]\rangle$.
It is clear that $B_\xi$ restricts to a non-degenerate skew-symmetric bilinear form on $\mathfrak g/\mathfrak g(\xi)$, whence $\dim \mathfrak g+\dim\mathfrak g(\xi)$ is even.

A polarization $\mathfrak b$ of $\xi$ as above is a Lie subalgebra of $\mathfrak g$, which is isotropic with respect to $B_\xi$ ({\em i.e.} $\xi$ defines a character for $\mathfrak b$) and of maximal dimension ({\em i.e.} the dimension of $\mathfrak b$ equals $(\dim\mathfrak g+\dim\mathfrak g(\xi))/2$).
The isotropy condition on $\mathfrak b$ makes it automatically an algebraic subalgebra of $\mathfrak g$.
A polarization $\mathfrak b$ of $\xi$ satisfies Pukanszky's condition, if $\mathfrak b=\mathfrak g(\xi)+\mathfrak b_u$, where $\mathfrak b_u$ denotes the unipotent radical of $\mathfrak b$, see also~\cite[Subsubsection 1.5.2]{T1} for equivalent characterizations of Pukanszky's condition.
\begin{Prop}\label{p-polar}
For a finite-dimensional Lie algebra $\mathfrak g$ and a general element $\xi$ of $\mathfrak g^*$, the reduction space $\mathrm H^0(B_\hbar)$ associated to the coisotropic submanifold $\xi+\mathfrak b^\perp$ equals $\mathbb K[\![\hbar]\!]$.
\end{Prop}
\begin{proof}
The proof makes use of~\cite[Lemma 1.5.2.2]{T1} concerning equivalent characterizations of Pukanszky's condition: in fact, we use the equivalence between the above characterization of Pukanszky's condition and the one stating that $\mathrm{Ad}^*(B)(\xi)=\xi+\mathfrak b^\perp$, for $B$ an algebraic, connected subgroup of $G$ (an algebraic group with Lie algebra $\mathfrak g$) with Lie algebra $\mathfrak b$.
According to~\cite[Subsubsection 1.5.2]{T1}, $\mathrm{Ad}^*(B)(\xi)$ is a Zarisky open subset of $\xi+\mathfrak b^\perp$.

We consider an element $f=f_0+\mathcal O(\hbar)$ of $B_\hbar^0$, $f_i$ in $\mathbb K[\xi+\mathfrak b^\perp]$.
If $f$ is $\mathrm d_{B_\hbar}$-closed, then $f_0$ is $B$-invariant, which, by the above, means that $f_0$ is constant.
By recurrence, the higher order identities reduce to $f_i$ $B$-invariant, $i\geq 1$, whence the claim follows.
\end{proof}

\subsection{Products on quantum reduction algebras}\label{ss-2-5}
After having discussed in some detail the relevant quantum reduction algebras which we will encounter in the sequel, we are now interested in a detailed discussion of the existence of associative products on quantum reduction algebras.

We consider the dg vector space $B_\hbar$ in its full generality for a finite-dimensional Lie algebra $\mathfrak g$ and a Lie subalgebra $\mathfrak h$ thereof, to which we associate a dg algebra $B$, whose deformation quantization $B_\hbar$ is a flat $A_\infty$-algebra: in particular, this means that the graded vector space $\mathrm H^\bullet(B_\hbar)$ is endowed with an associate product.
More precisely, the $A_\infty$-structure $\mu_B+\mathcal U_B(\hbar\pi)$, where $\pi$ denotes here the $P_\infty$-structure on $B$ Fourier-dual to the Poisson bivector on $X=\mathfrak g^*$, consists of infinitely many Taylor components 
\[
\mathcal U_B(\hbar\pi)^n:B_\hbar^{\otimes n}\to B_\hbar[2-n],\ n\geq 1,
\]
which satisfy an infinite series of quadratic identities between them.
For our purposes, we need only know that $\mathcal U_B(\hbar\pi)^1=\mathrm d_{B_\hbar}$, and that $\mu_B+\mathcal U_B(\hbar\pi)^2$ defines a $\mathbb K$-bilinear pairing of degree $0$ on $B_\hbar$, which is compatible with $\mathrm d_{B_\hbar}$ and which is associative up to a the homotopy $\mathcal U_B(\hbar\pi)^3$ with respect to $\mathrm d_{B_\hbar}$.
Therefore, $\mu_B+\mathcal U_B(\hbar\pi)^2$ descends to the quantum reduction space $\mathrm H^\bullet(B_\hbar)$ to an associative product, which we denote for simplicity by $\star_{B_\hbar}$: its restriction to $\mathrm H^0(B_\hbar)$ defines an obvious deformation of the commutative product on $B^0$.

We refer to~\cite[Section 3]{CT} for a careful description of the deformed product $\star_{B_\hbar}$ on $\mathrm H^0(B_\hbar)$ in the case of a symmetric pair $(\mathfrak g,\sigma)$ (with Cartan decomposition $\mathfrak g=\mathfrak k\oplus\mathfrak p$) and of a character $\chi$ of the Lie subalgebra $\mathfrak k$.
We only observe that in~\cite[Section 3]{CT}, the authors use the notation $\star_{\mathrm{CF},\lambda}=\star_{B_\hbar}$.

\section{Applications of biquantization in Lie theory for symmetric pairs}\label{s-3}
In the present Section, we discuss the first relevant applications of biquantization, as discussed in detail in the previous Section, to concrete problems in Lie theory in the framework of a symmetric space $(\mathfrak g,\sigma)$ with standard Cartan decomposition $\mathfrak g=\mathfrak k\oplus\mathfrak p$.

In fact, the results presented here are the revisited versions of the results of~\cite[Section 4]{CT} taking into account the more precise and correct approach to biquantization presented in the previous Section.

\subsection{A comparison between the quantum deformed product and Rouvi\`ere's product}\label{ss-3-1}
We consider the quantum reduction algebra $(\mathrm H^0(B_\hbar),\star_{B_\hbar})$ for a symmetric pair $(\mathfrak g,\sigma)$ with the standard Cartan decomposition.
On the other hand, we consider the quantum deformed algebra $(A_\hbar,\star_{A_\hbar})$ and the $A_\infty$-$A_\hbar$-$B_\hbar$-bimodule $K_\hbar$: unraveling the $A_\infty$-identities for $K_\hbar$, we see that $K_\hbar$, which is concentrated in degree $0$, becomes actually an $(A_\hbar,\star_{A_\hbar})$-$(\mathrm H^0(B_\hbar),\star_{B_\hbar})$-bimodule, and we denote by $\star_L$ and $\star_R$ the respective left $(A_\hbar,\star_{A_\hbar})$- and right $(\mathrm H^0(B_\hbar),\star_{B_\hbar})$-action on $K_\hbar$.

More explicitly, in the present framework, we have $A_\hbar=\mathrm S(\mathfrak g)[\![\hbar]\!]$, $B_\hbar^0=K_\hbar=\mathrm S(\mathfrak p)[\![\hbar]\!]$.
We observe that the linearity of the Poisson bivector $\pi$ on $X=\mathfrak g^*$ has an interesting by-product, as already remarked in~\cite[Subsubsection 8.3.1]{K}: although there are infinitely many bidifferential operators appearing in the deformed product $\star_{A_\hbar}$, the action of the infinite series $\star_{A_\hbar}$, for $\hbar=1$, on $\mathrm S(\mathfrak g)\otimes\mathrm S(\mathfrak g)$ is well-defined, as can be proved by inspecting the integral weights and counting degrees.
Similar arguments hold true also for the components of the $A_\infty$-structure $\mu_B+\mathcal U_B(\hbar\pi)$ on $B_\hbar$ and for the $A_\infty$-$A_\hbar$-$B_\hbar$-bimodule structure on $K_\hbar$.
As a consequence, we may consider the $A_\infty$-algebras $A_\hbar$ and $B_\hbar$ and the $A_\infty$-$A_\hbar$-$B_\hbar$-bimodule $K_\hbar$ as polynomial deformations with respect to $\hbar$, and in particular we may safely consider the value of the parameter $\hbar=1$.

Therefore, we consider the associative algebras $(A,\star_A)$, the $A_\infty$-algebra $(B,\mu_B+\mathcal U_B(\pi))$, with corresponding quantum reduction algebra $(\mathrm H^0(B),\star_B)$, and the $(A,\star_A)$-$(\mathrm H^0(B),\star_B)$-bimodule $(K,\star_L,\star_R)$.
The deformed differential $\mathrm d_B$ on $B^0$ is now a differential operator from $\mathrm S(\mathfrak p)$ to $\mathrm S(\mathfrak p)\otimes \mathfrak k^*$ of infinite order, whose action is well-defined.
In a similar way, the pairing $\mu_B+\mathcal U_B(\pi)^2$ is a well-defined bidifferential operator on $\mathrm S(\mathfrak p)$ of infinite order, and so $\star_L$ and $\star_R$ (we should more precise on $\star_R$, as we should speak of the infinite-order bidifferential operator on $\mathrm S(\mathfrak p)$ coming from the $A_\infty$-$A$-$B$-bimodule structure on $K$).

The deformed product $\star_A$ on $A=\mathrm S(\mathfrak g)$ has been explicitly characterized in~\cite{K,BCKT}.
More precisely, we consider the following function on $\mathfrak g$, 
\begin{equation}\label{eq-duf}
q(x)=\underset{\mathfrak g}\det \!\left(\frac{\sinh\!\left(\frac{\mathrm{ad}(x)}2\right)}{\frac{\mathrm{ad}(x)}2}\right),
\end{equation}
which is analytic in a neighborhood of $0$.
It can be expanded in a power series of the polynomials $c_n(x)=\mathrm{tr}_\mathfrak g(\mathrm{ad}(x)^n)$, $n\geq 1$.
Alternatively, it may be viewed as an element of the completed symmetric algebra $\widehat{\mathrm S}(\mathfrak g^*)$ of the dual of $\mathfrak g$, and as such, as an invertible, $\mathfrak g$-invariant, infinite-order differential operator with constant coefficients acting on $A$.
Similar arguments hold true also if we consider its square root $\sqrt{q}$: we denote by $\partial_{\sqrt{q}}$ the corresponding invertible, $\mathfrak g$-invariant, infinite-order differential operator on $A$. 
We further denote by $\beta$ the Poincar\'e--Birkhoff--Witt (shortly, PBW) isomorphism from $\mathrm S(\mathfrak g)$ to $\mathrm U(\mathfrak g)$, {\em i.e.} $\beta$ is the symmetrization morphism
\[
\mathrm S(\mathfrak g)\ni x_1\cdots x_n\mapsto \frac{1}{n!}\sum_{\sigma\in\mathfrak S_n}x_{\sigma(1)}\cdots x_{\sigma(n)}\in\mathrm U(\mathfrak g),\ x_j\in \mathfrak g,\ j=1,\dots,n.
\]
Then, the deformed product $\star_A$ on $\mathrm S(\mathfrak g)$ is related with the product in $\mathrm U(\mathfrak g)$ {\em via}
\begin{equation}\label{eq-DK}
\beta\!\left(\partial_{\sqrt{q}}(f_1)\star_A\partial_{\sqrt{q}}(f_2)\right)=\beta(\partial_{\sqrt{q}}(f_1))\cdot\beta(\partial_{\sqrt{q}}(f_2)),\ f_i\in A,\ i=1,2, 
\end{equation}
and $\cdot$ denotes here the product in $\mathrm U(\mathfrak g)$.

The way the operator $\partial_{\sqrt{q}}$ arises in the framework of deformation quantization has been elucidated in detail in~\cite[Subsubsections 8.3.1, 8.3.2 and 8.3.3]{K}, combining the results therein with~\cite{Sh}.
We also refer to~\cite[Part II]{BCKT} for a complete overview of the applications of deformation quantization as in~\cite{K} in Lie theory.

The motivation for the following computations lies in the comparison in Identity~\eqref{eq-DK} between the UEA $\mathrm U(\mathfrak g)$ and the quantum deformed algebra $(A,\star_A)$, which are related precisely by the ``strange'' automorphism $\partial_{\sqrt{q}}$ of the symmetric algebra $\mathrm S(\mathfrak g)$, which appears also in Duflo's Theorem: namely, the composition $\beta\circ \partial_{\sqrt{q}}$ defines an algebra isomorphism between $\mathrm S(\mathfrak g)^\mathfrak g$ and the center of $\mathrm U(\mathfrak g)$.
The main point is that Kontsevich's deformed product $\star_A$ contains bidifferential operators, which are represented in terms of wheel-like graphs: such graphs are precisely responsible for the appearance of the ``strange'' automorphism $\partial_{\sqrt{q}}$.

Quite similarly, in the case of a symmetric pair $(\mathfrak g,\sigma)$, we may consider the associative algebra $(\mathrm H^0(B),\star_B)$, where $\mathrm H^0(B)=\mathrm S(\mathfrak p)^\mathfrak k$.
It is worth observing that Poisson reduction methods yield a Poisson structure on $\mathrm S(\mathfrak p)^\mathfrak k$ simply by restriction, and the product $\star_B$ defines a deformation quantization of $\mathrm S(\mathfrak p)^\mathfrak k$ in the sense of Kontsevich. 
On the other hand, for any choice of a character $\chi$ of $\mathfrak k$ ({\em i.e.} a $1$-dimensional $\mathfrak k$-representation on $\mathbb K$), the PBW isomorphism $\beta$ induces a direct sum decomposition $\mathrm U(\mathfrak g)=\beta(\mathrm S(\mathfrak p))\oplus \mathrm U(\mathfrak g)\cdot \mathfrak k^{-\chi}$, where $\mathfrak k^{-\chi}$ denotes the affine subspace of $\mathrm U(\mathfrak g)$ spanned by elements of the form $x-\chi(x)$, $x$ in $\mathfrak k$.
Then, Rouvi\`ere defines also a ``deformation quantization'' of the Poisson algebra $\mathrm S(\mathfrak p)^\mathfrak k$ {\em via} the formula
\begin{equation}\label{eq-rouv}
\beta(f_1\# f_2)=\beta(f_1)\cdot\beta(f_2)\ \text{modulo $\mathrm U(\mathfrak g)\cdot \mathfrak k^{-\chi}$,}\ f_i\in\mathrm S(\mathfrak p)^\mathfrak k,\ i=1,2.
\end{equation}

The PBW isomorphism (of vector spaces) is obviously $\mathfrak g$-invariant, hence it is automatically $\mathfrak k$-invariant: therefore, $\beta$ restricts to a $\mathfrak k$-invariant isomorphism of vector spaces from $\mathrm S(\mathfrak p)$ to $\beta(\mathrm S(\mathfrak p))\subseteq \mathrm U(\mathfrak g)$.
In particular, from the above decomposition $\mathrm U(\mathfrak g)=\beta(\mathrm S(\mathfrak p))\oplus \mathrm U(\mathfrak g)\cdot \mathfrak k^{-\chi}$, it follows immediately that the right-hand side of Identity~\eqref{eq-rouv} defines a unique bilinear pairing $\#$ on $\mathrm S(\mathfrak p)$, which restricts to an associative product on $\mathfrak k$-invariant elements.
%\begin{Rem}\label{r-comm}
%We observe another important issue of Identity~\eqref{eq-rouv}: to bring a product in $\mathrm U(\mathfrak g)$ of the form $\beta(f_1)\cdot\beta(f_2)$ to an element of the form $\beta(g)$, for some $g$ in $\mathrm S(\mathfrak p)$ working inside $\mathrm U(\mathfrak g)$, we permute factors in $\mathrm U(\mathfrak g)$, getting by definition Lie brackets.
%In the case at hand, the Cartan relations for the symmetric pair $(\mathfrak g,\sigma)$ imply that $\beta(f_1)\cdot\beta(f_2)$ and $\beta(g)$ differ by an element which is actually in $\mathrm U(\mathfrak g)\cdot \mathfrak k_1^{-\chi}$, where $\mathfrak k_1=[\mathfrak p,\mathfrak p]\subseteq \mathfrak k$. 
%We observe that $\mathfrak g_1=\mathfrak k_1\oplus\mathfrak p$ is a symmetric subpair of $\mathfrak g$.
%This apparently harmless observation will prove quite useful in later computations.
%\end{Rem}
We are now going to compare the products $\star_B$ and $\#$ {\em via} biquantization techniques.
As one could naturally guess from Identity~\eqref{eq-DK}, the two products on $\mathrm S(\mathfrak p)^\mathfrak k$ do not coincide, but are related to each other in a similar fashion, {\em i.e.} through a ``relative'' counterpart of Duflo's ``strange'' automorphism.
The novelty of the approach through biquantization is the fact that we use it to compare $\star_B$ on $\mathrm H^0(B)$ with $\star_A$ on $A$; the rest of the proof, {\em i.e.} the comparison of different automorphisms of $\mathrm S(\mathfrak p)$ similar in shape to Duflo's ``strange'' automorphism, is really similar to the proof presented in~\cite[Subsection 3.1]{K} with due modifications. 

Of course, the upcoming discussion can be generalized to the framework of some Lie subalgebra $\mathfrak h$ of any finite-dimensional Lie algebra $\mathfrak g$ over $\mathbb K$: we will discuss generalizations of the results presented here elsewhere, in particular in relationship with equivalences of categories of representations of Lie algebras and corresponding subalgebras and with the relative Duflo conjecture.

\subsubsection{A version of Duflo's ``strange'' automorphism for symmetric pairs}\label{sss-3-1-1}
Using the previous notation and conventions, we define the following operator
\begin{equation}\label{eq-A-wheel}
\mathcal A(f)=f\star_L 1,\ f\in A=\mathrm S(\mathfrak g)
\end{equation}
from $A$ to $K=\mathrm S(\mathfrak p)$.
Of course, we could have first defined the operator $\mathcal A_\hbar$ from $A_\hbar$ to $K_\hbar$, for $\hbar$ a formal parameter.
By its very construction, $\mathcal A_\hbar$ is a deformation of the surjective projection from $A$ to $K$, which we denote by $\pi$.
By its very construction, $\mathcal A_\hbar=\pi\circ \mathbb A_\hbar$, where $\pi$ is extended $\hbar$-linearly to $A_\hbar$, while $\mathbb A_\hbar$ is a formal series of differential operators on $A$, where $\mathbb A_0=\mathrm{id}$, and $\mathbb A_n$ (the coefficient of degree $n$ with respect to $\hbar$) is a differential operator of order $n$.
In particular, it is clear that $\mathbb A_\hbar$ is invertible.
By the same arguments as before, we may safely set $\hbar=1$, and thus we get an invertible differential operator $\mathbb A$ on $A$ of infinite order.

We consider $(A,\star_A)$ as a left $(A,\star_A)$-module: then, $\mathcal A$ is a surjective morphism of $(A,\star_A)$-modules from $(A,\star_A)$ to $(K,\star_L)$, whence $K\cong A/I$, where $I=\mathrm{Ker}(\mathcal A)$.

By the very definition of $\mathcal A$, $I=A\star_A \mathbb A^{-1}(\mathfrak k)$: in fact, $A\star_A \mathbb A^{-1}(\mathfrak k)\subseteq I$ by its very construction.
To prove the opposite inclusion, we re-introduce momentarily the formal parameter $\hbar$.
For $\hbar=0$, the ideal $I=\langle\mathfrak k\rangle$ of $A$ is finitely-generated.
In fact, as $I$ is the two-sided ideal of $A$, viewed here as a commutative algebra, generated by the ideal $\mathfrak k$, viewed here as an ideal of linear functions on $X=\mathfrak g^*$.
In particular, there is a surjective morphism $A^{\oplus \dim\mathfrak k}\to I\to 0$ of $A$-modules: by the previous argument, there is a morphism of left $A_\hbar$-modules $A\star_{A_\hbar} \mathbb A_\hbar^{-1}(\mathfrak k)\to I_\hbar$, which is a formal $\hbar$-deformation of the surjective morphism $A^{\oplus \dim\mathfrak h}\to I$, whence the surjectivity of the deformed morphism follows.
Therefore, $I_\hbar=A_\hbar\star_{A_\hbar} \mathbb A_\hbar^{-1}(\mathfrak k)$, whence the claim follows by setting safely $\hbar=1$.

It remains to compute $\mathbb A^{-1}(\mathfrak k)$.

Writing $\mathbb A=\mathrm{id}+\sum_{n\geq 1}\mathbb A_n$, $\mathbb A_n$, $n\geq 1$, has order $n$ by construction; furthermore, $\mathbb A_n$, $n\geq 1$, has no constant term.
Namely, $\mathbb A_n$ is specified by differential operators associated to admissible graphs in $\mathcal G_{n,1}$: recalling from the previous Section the construction of the differential operator associated to $\Gamma$ admissible of type $(n,1)$, if $\Gamma$ has no edge pointing to the only vertex of the second type, the corresponding integral weight vanishes by a dimensional argument.

The operator $\mathbb A^{-1}$ is completely determined by the power series expansion of $\mathbb A$: once again, it is of the form $\mathrm{id}+\sum_{n\geq 1} \widetilde{\mathbb A}_n$, where $\widetilde{\mathbb A}_n$ has no constant term, for $n\geq 1$, as follows by an easy computation.

We thus compute $\mathbb A(k)$, for $k$ a general element of $\mathfrak k$.
Because of degree reasons, see also~\cite[Subsubsection 8.3.1]{K}, $\mathbb A(k)=k+\mathbb A_1(k)$.
Further, it is readily checked that $\mathbb A_1$ is the sum of two differential operators, associated to the following admissible graphs of type $(1,2)$
\bigskip
\begin{center}
\resizebox{0.45 \textwidth}{!}{\begin{picture}(0,0)%
\includegraphics{def-proj.pstex}%
\end{picture}%
\setlength{\unitlength}{3947sp}%
\begingroup\makeatletter\ifx\SetFigFont\undefined%
\gdef\SetFigFont#1#2#3#4#5{%
  \reset@font\fontsize{#1}{#2pt}%
  \fontfamily{#3}\fontseries{#4}\fontshape{#5}%
  \selectfont}%
\fi\endgroup%
\begin{picture}(9511,3667)(1146,-6416)
\end{picture}
}\\
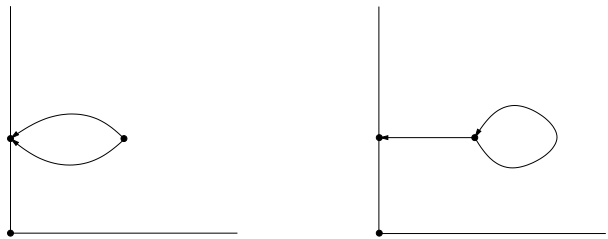
\text{Figure 7 - The only two admissible graphs of type $(1,2)$ appearing in $\mathbb A_1$}
\end{center}
\bigskip

The contribution of the first graph is trivial, because $k$ is viewed as a linear function on $X$.

We consider the second graph: we want to observe that such a graph did not appear in the computations performed in~\cite{CT}.
First of all, the corresponding integral weight is 
\begin{equation}\label{eq-loop-weight}
\int_{\mathcal C_{1,2}^+}\rho\omega^{+,-},
\end{equation}
omitting wedge products.
In fact, in the superpropagator $\omega_e$, only two of the $4$-colored propagators are non-trivial, namely $\omega^{+,+}$ and $\omega^{+,-}$ by construction; since it acts as a derivation on an element of $\mathfrak k$, by its very definition, the part with $\omega^{+,+}$ vanishes.
\begin{Lem}\label{l-loop-weight}
The integral $\int_{\mathcal C_{1,2}^+}\mathrm d\eta\wedge\omega^{+,-}$ equals $\frac{1}4$.
\end{Lem}
\begin{proof}
The integral weight associated to the previous admissible graph is $\int_{\mathcal C_{1,2}^+}\mathrm d\eta\omega^{+,-}$, where we have suppressed wedge products between the forms in the integrand.

The $1$-form $\rho$ is exact, whence 
\[
\int_{\mathcal C_{1,2}^+}\mathrm d\eta\omega^{+,-}=\int_{\partial C_{1,2}^+}\eta\omega^{+,-},
\]
where we have used notation from Subsection~\ref{ss-2-2}.
Hence, it suffices to compute all boundary contributions to evaluate the integral.

The boundary strata of $\mathcal C_{1,2}^+$ are of the type $\mathcal C_{A,B}^+\times \mathcal C_{1\smallsetminus A,[2]\smallsetminus B\sqcup \{\bullet\}}^+$, for $A$ a subset of $[1]$ and $B$ an ordered subset of $[2]$, such that $0\leq |A|\leq 1$, $0\leq |B|\leq 2$ and $|A|+|B|\leq 2$.
Dimensional arguments imply that there are only two types of such boundary strata, $\mathcal C_{0,2}^+\times \mathcal C_{1,1}^+$ and $\mathcal C_{0,3}^+\times \mathcal C_{1,0}^+$, which correspond to five different situations.

We consider the boundary stratum $\mathcal C_{0,2}^+\times\mathcal C_{1,1}^+$, which corresponds to the situation where $i_1)$ the point on the positive real axis approaches the origin, $ii_1)$ the point in the interior of the first quadrant collapses to the point on the positive real axis or $iii_1)$ the point $1$ approaches the origin.
The boundary conditions for $\omega^{-,+}$ yield triviality of the contributions $i_1)$ and $iii_1)$; the second one yields $\frac{1}4\int_{\mathcal C_{1,1}^+}\omega^+=\frac{1}4$, and we have already included orientation signs.

The boundary stratum $\mathcal C_{0,3}^+\times\mathcal C_{1,0}^+$ corresponds to the point $1$ approaching either the positive imaginary axis or the positive real axis: in the first case, the corresponding contribution vanishes by means of the boundary conditions for $\omega^{+,-}$, while in the second case, the function $\eta$ vanishes when its argument approaches the real axis.
\end{proof}
Recalling now the construction of the superpropagators in biquantization from Subsubsection~\ref{sss-2-3-2}, the differential operator corresponding to the second graph in Figure 7 is
\[
\frac{1}4\left[\mathrm{tr}_{\mathfrak p}(\mathrm{ad}_\mathfrak{k}(\bullet))-\mathrm{tr}_{\mathfrak k}(\mathrm{ad}_\mathfrak{k}(\bullet)\right)]=\delta(\bullet)-\frac{1}4\mathrm{tr}_\mathfrak g(\mathrm{ad}(\bullet)),
\]
whence $\mathbb A(k)=k+\delta(k)-\frac{1}4\mathrm{tr}_\mathfrak g(\mathrm{ad}(k))$.

Therefore, we have
\[
\mathbb A^{-1}(\mathbb A(k))=k=\mathbb A^{-1}(k)+\delta(k)-\frac{1}4\mathrm{tr}_\mathfrak g(\mathrm{ad}(k)),\ \text{whence}\ \mathbb A^{-1}(k)=k-\delta(k)+\frac{1}4\mathrm{tr}_\mathfrak g(\mathrm{ad}(k)),\ k\in\mathfrak k.
\]
We observe that we have used the fact that the terms of $\mathbb A_\hbar^{-1}$ of degree higher or equal than $1$ are differential operators without constant term.

Putting all previous arguments together, we have the identification of right $(A,\star_A)$-modules 
\[
I=A\star_A \mathfrak k^{-\delta+\frac{1}4\mathrm{tr}_\mathfrak g\circ\mathrm{ad}}.
\]

Further, we consider the restriction of $\mathcal A$ to $K$, viewed here as a subalgebra of $A$.
First of all, any admissible graph $\Gamma$ of type $(n,2)$ yielding a possibly non-trivial contribution to $\mathcal A$ has exactly $2n$ edges because of dimensional reasons.
To any edge $e$ of $\Gamma$ corresponds a superpropagator $\omega_e$, whose components are only of type $(+,+)$ or $(+,-)$: it follows immediately from the boundary conditions for both of them that $\Gamma$ has no edge departing from the only vertex on the positive real axis, and that any edge pointing to this vertex has a corresponding superpropagator with color $(+,+)$.
This excludes immediately double edges pointing to the only vertex of $\Gamma$ on the positive real axis. 
In particular, this means that any admissible graph $\Gamma$ of type $(n,2)$ yielding a non-trivial contribution to $\mathcal A$ in the present situation must satisfy the following rule: from any vertex of the first of $\Gamma$ departs at most one edge pointing to the only vertex on the positive real axis (of course, because of the presence of short loops, the initial and final point of such an edge may coincide).
We assume therefore that $p\leq n$ edges have the only vertex of the second type on the positive real axis as endpoint.
If $p<n$, there are then $2n-p$ edges, whose endpoints are both vertices of the first type (double edges and short loops are allowed). 
Since every edge is associated to a derivation, the polynomial degree associated to the vertices of the first type of $\Gamma$ is $-n+p$ (counting $n$ because of the linearity of the Poisson structure and $p-2n$ derivations), which is strictly negative, leading to a contradiction.
Therefore, form any vertex of the first type of $\Gamma$ depart exactly one edge to the only vertex of the positive real axis and one to a vertex of the first type; the polynomial degree of the corresponding differential operator on $K$ is immediately $0$. 
Because of the linearity of the Poisson structure, exactly one edge has a vertex of the first type as final point, whence admissible graphs of type $(n,2)$ contributing non-trivially to $\mathcal A$ are disjoint unions of wheel-type graphs; we observe that the $1$-wheel may in principle appear.
Further, only wheel-like graphs with $n$ even contribute (possibly) non-trivially to $\mathcal A$.
Namely, by the previous argument, from every vertex of the first type departs exactly one edge to the only vertex of the second type on the positive real axis, whose color is $(+,+)$ and whose operator-valued part corresponds to derivation and contraction with respect to $\mathfrak p$.
The Cartan relation $[\mathfrak k,\mathfrak p]\subseteq \mathfrak p$ implies that at each vertex of the first type in a wheel-like graph $\Gamma$, the edge arriving at such a vertex must have color either $(+,+)$ or $(+,-)$, while the edge departing from it on th wheel must have opposite color.
In other words, the edges of the cycle in a wheel-like graph $\Gamma$ must have alternating colours $(+,+)$ and $(+,-)$: this, in turn, excludes immediately $n$-wheels with $n$ odd.

Summarizing all previous arguments, the restriction of the operator $\mathcal A$ to $K$, which we denote (improperly) by the same symbol, defines an invertible, translation-invariant differential operator on $K$.
Its symbol, regarded as an element of the completed symmetric algebra $\widehat{\mathrm S}(\mathfrak p)$ and defined through $j_\mathcal A(x)=e^{-x}\mathcal A(e^x)$, has the explicit form 
\[
j_\mathcal A(x)=\exp\left(\sum_{n\geq 1}W_{2n}^\mathcal A\mathrm{tr}_\mathfrak p(\mathrm{ad}^{2n}(x))\right),\ x\in\mathfrak p,
\]
where $W_{2n}^\mathcal A$, $n\geq 1$, denotes the integral weight of following wheel-like graph:
\bigskip
\begin{center}
\resizebox{0.3 \textwidth}{!}{\input{wheel_A_2n.pstex_t}}\\
\text{Figure 8 - The wheel-like graph with integral weight $W_{2n}^\mathcal A$}
\end{center}
\bigskip
We observe that $j_\mathcal A$ is analytic in a neighborhood of $0$ in $\mathfrak p$.
We observe that the Cartan relations for the symmetric pair $(\mathfrak g,\sigma)$ imply immediately that, for a general element $x$ of $\mathfrak p$, $\mathrm{ad}(x)^2$ is a well-defined endomorphism of $\mathfrak p$, thus all even powers of the adjoint representation restricted to $\mathfrak p$: hence the above expression is well-defined.
Finally, all previous computations imply also the direct sum decomposition of $(A,\star_A)$:
\[
A=K\oplus (A\star_A\mathfrak k^{-\delta+\frac{1}4\mathrm{tr}_\mathfrak g\circ\mathrm{ad}}).
\]

On the other hand, $K=B^0=\mathrm S(\mathfrak p)$ by definition.
We may therefore consider the endomorphism $\mathcal B$ of $K$ defined through 
\begin{equation}\label{eq-B-wheel}
\mathcal B(f)=1\star_R f,\ f\in K=B^0=\mathrm S(\mathfrak p).
\end{equation}
We may repeat almost {\em verbatim} the previous arguments to evaluate explicitly the operator $\mathcal B$: it is an invertible, translation-invariant differential operator on $K$ of infinite order, with symbol in $\widehat{\mathrm S}(\mathfrak p)$ given by
\[
j_\mathcal B(x)=\exp\left(\sum_{n\geq 1}W_{2n}^\mathcal B\mathrm{tr}_\mathfrak p(\mathrm{ad}^{2n}(x))\right),\ x\in\mathfrak p,
\]
where $W_{2n}^\mathcal B$, $n\geq 1$, denotes the integral weight of following wheel-like graph:
\bigskip
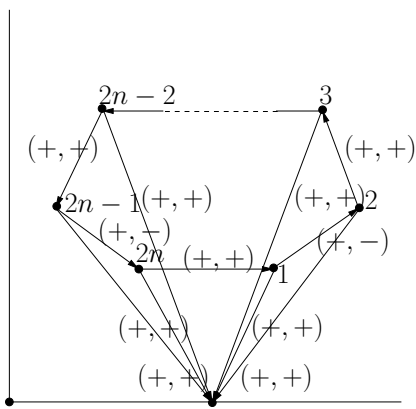
\begin{center}
\resizebox{0.3 \textwidth}{!}{\input{wheel_B_2n.pstex_t}}\\
\text{Figure 9 - The wheel-like graph with integral weight $W_{2n}^\mathcal B$}
\end{center}
\bigskip
Once again, notice that $j_\mathcal B$ is analytic in a neighborhood of $0$ in $\mathfrak p$.

At this point one could wonder whether or not the integral weights $W_{2n}^\mathcal A$ and $W_{2n}^\mathcal B$ coincide (which would imply that $j_\mathcal A=j_\mathcal B$).
First of all, we observe that in the Formul\ae\ for $j_\mathcal A$ and $j_\mathcal B$, only weights of even wheels appear because of the vanishing of the the differential operators corresponding to odd wheel-like graphs.
In fact, {\em e.g.} the $1$-wheels $\mathcal W_1^\mathcal A$ and $\mathcal W_1^\mathcal B$ are both computable and yield distinct results.

This can be seen either by computing separately both integral weights or by computing {\em e.g.} $W_1^\mathcal A$ and finding then a relationship between $W_1^\mathcal A$ and $W_1^\mathcal B$.

First of all, the integral weight $W_1^\mathcal A$ is explicitly 
\[
W_1^\mathcal A=\int_{\mathcal C_{1,2}^+}\rho\omega^{+,+}.
\]
It can be computed by the same technique used in Lemma~\ref{l-loop-weight}.
Of course, there are certain differences to be taken into account, namely, the different boundary conditions satisfied by the $4$-colored propagator $\omega^{+,+}$.
Here, the only boundary strata of $\mathcal C_{1,2}^+\cong\mathcal C_{1,1,0}^+$ which yield non-trivial contributions are $i)$ the stratum corresponding to the approach of the point in $Q^{+,+}$ to the only point on $i\mathbb R^+$ and $ii)$ the approach of the only point on $i\mathbb R^+$ to the origin.
Both integrals are readily computed, as well as their orientation signs, which then yield the desired result.

We now consider the following admissible graph of type $(2,1)$:
\bigskip
\begin{center}
\resizebox{0.2 \textwidth}{!}{\begin{picture}(0,0)%
\includegraphics{aer_wheel.pstex}%
\end{picture}%
\setlength{\unitlength}{3947sp}%
\begingroup\makeatletter\ifx\SetFigFont\undefined%
\gdef\SetFigFont#1#2#3#4#5{%
  \reset@font\fontsize{#1}{#2pt}%
  \fontfamily{#3}\fontseries{#4}\fontshape{#5}%
  \selectfont}%
\fi\endgroup%
\begin{picture}(4867,4867)(3546,-6416)
\put(5791,-4771){\makebox(0,0)[lb]{\smash{\SetFigFont{20}{24.0}{\rmdefault}{\mddefault}{\updefault}{\color[rgb]{0,0,0}$(+,+)$}%
}}}
\end{picture}
}\\
\text{Figure 10 - The ``aerial'' $1$-wheel}
\end{center}
\bigskip
It represents a $2$-form on the $3$-dimensional smooth manifold with corners $\mathcal C_{2,1}^+$, therefore, in virtue of Stokes' Theorem, 
\[
\int_{\mathcal C_{2,1}^+}\mathrm d(\mathrm d\eta\omega^{+,+})=0=\int_{\partial C_{2,1}^+}\mathrm d\eta\omega^{+,+}.
\]
The boundary strata of codimension $1$ of $\mathcal C_{2,1}^+$ have been illustrated explicitly in Subsubsection~\ref{sss-2-2-2}: either because of the boundary conditions for $\omega^{+,+}$ and $\rho$ or because of dimensional reasons, it is not difficult to verify that only three boundary strata yield non-trivial contributions, namely the strata $\alpha$, $\theta$ and $\zeta$.

The boundary stratum $\theta$, resp.\ $\zeta$, yields precisely $W_1^\mathcal B$, resp.\ $W_1^\mathcal A$; on the other hand, it is not difficult to verify that the boundary stratum $\alpha$ yields the non-trivial contribution $1/4$, as can be verified by a direct computation.
Thus, in general, we cannot expect the weights $W_n^\mathcal A$ and $W_n^\mathcal B$ to coincide.

\subsubsection{Explicit comparison of the products of Rouvi\`ere and Cattaneo--Felder}\label{sss-3-1-2}
The operator $\mathcal A$ is surjective, and its restriction to $K=\mathrm S(\mathfrak p)\subseteq A$ is an automorphism, while $\mathcal B$ is an automorphism of $K$, whence 
\[
1\star_R f=\mathcal B(f)=\mathcal A(\mathcal A^{-1}(\mathcal B(f)))=\mathcal A^{-1}(\mathcal B(f))\star_L 1,\ f\in K=\mathrm S(\mathfrak p).
\]

We now recall from Subsubsection~\ref{sss-2-4-1} that the quantum reduction algebra (at $\hbar=1$) $\mathrm H^0(B)=\mathrm S(\mathfrak p)^\mathfrak k$.
We thus consider two elements $f_i$, $i=1,2$, of $\mathrm H^0(B)$, endowed with the deformed associative product $\star_B$; then $(K,\star_R)$ becomes a right $(\mathrm H^0(B),\star_B)$-module, and the latter module structure is compatible with the left $(A,\star_A)$-module structure on $(K,\star_L)$, whence
\[
\begin{aligned}
&1\star_R (f_1\star_B f_2)=\left(\mathcal A^{-1}(\mathcal B(f_1\star_B f_2))\right)\star_L 1=\\
=& (1\star_R f_1)\star_R f_2=\left(\mathcal A^{-1}(\mathcal B(f_1))\star_L 1\right)\star_R f_2=\mathcal A^{-1}(\mathcal B(f_1))\star_L (1\star_R f_2)=\mathcal A^{-1}(\mathcal B(f_1))\star_L \left(\mathcal A^{-1}(\mathcal B(f_2))\star_L 1\right)=\\
=& \left(\mathcal A^{-1}(\mathcal B(f_1))\star_A \mathcal A^{-1}(\mathcal B(f_2))\right)\star_L 1,
\end{aligned}
\]
whence from the previous computations follows
\begin{equation}\label{eq-rouviere}
(\mathcal A^{-1}\circ \mathcal B)(f_1\star_B f_2)=(\mathcal A^{-1}\circ \mathcal B)(f_1)\star_A (\mathcal A^{-1}\circ \mathcal B)(f_2)\ \text{modulo $A\star_A \mathfrak k^{-\delta+\frac{1}4\mathrm{tr}_\mathfrak g\circ\mathrm{ad}}$}.
\end{equation}

We apply the operator $\beta\circ\partial_{\sqrt{q}}$ on both sides of Identity~\eqref{eq-rouviere} and because of Identity~\eqref{eq-DK}, we get
\[
(\beta\circ\partial_{\sqrt{q}}\circ\mathcal A^{-1}\circ\mathcal B)(f_1\star_B f_2)=(\beta\circ\partial_{\sqrt{q}}\circ\mathcal A^{-1}\circ\mathcal B)(f_1)\cdot (\beta\circ\partial_{\sqrt{q}}\circ\mathcal A^{-1}\circ\mathcal B)(f_2)\ \text{modulo $\mathrm U(\mathfrak g)\cdot \mathfrak k^{-\delta+\frac{1}4\mathrm{tr}_\mathfrak g\circ\mathrm{ad}}$},
\]
for $f_i$ in $\mathrm S(\mathfrak p)^\mathfrak k$, $i=1,2$.

We introduce the function $J$ on $\mathfrak p$ defined {\em via}
\[
j(x)=\underset{\mathfrak p}\det\!\left(\frac{\sinh(\mathrm{ad}(x))}{\mathrm{ad}(x)}\right),\ x\in\mathfrak p.
\]
This function is the determinant of the exponential map for the symmetric pair $(\mathfrak g,\sigma)$: it can be written as a formal power series of traces in $\mathfrak p$ of even powers of the restriction to $\mathfrak p$ of the adjoint representation of $\mathfrak g$.
Similarly to what has been done before, we define $\partial_{\sqrt{j}}$ as the invertible, translation-invariant, $\mathfrak k$-invariant differential operator of infinite order on $\mathrm S(\mathfrak p)$, whose symbol is exactly the square root of $j$.

We define a modified version of the previously introduced Rouvi\`ere's product {\em via}
\[
\beta\!\left(\partial_{\sqrt{j}}(f_1\# f_2)\right)=\beta(\partial_{\sqrt{j}}(f_1))\cdot\beta(\partial_{\sqrt{j}}(f_2))\ \text{modulo $\mathrm U(\mathfrak g)\cdot\mathfrak k^{-\delta+\frac{1}4\mathrm{tr}_\mathfrak g\circ\mathrm{ad}}$.} 
\]
We observe that the left-hand side of the previous Identity is well-defined, because of the $\mathfrak k$-invariance and invertibility of $\partial_{\sqrt{j}}$. 

The key point is now the following identity, which puts into relationship the functions $q$, $J$, $j_\mathcal A$ and $j_\mathcal B$: 
\begin{equation}\label{eq-DR-symb}
j_\mathcal A(x)\sqrt{j(x)}=j_\mathcal B(x)\sqrt{q(x)},\ \forall x\in\mathfrak p.
\end{equation}

The previous identity is the relative version, in the case of a symmetric pair $(\mathfrak g,\sigma)$, of the results of~\cite[Subsubsection 8.3.4]{K}.
Its proof is a consequence of~\cite[Lemma 14 and Proposition 12]{CT}: it is result which is left unaltered by the changes to biquantization which have been previously discussed.
\begin{Thm}\label{t-Rou-CF}
For a general symmetric pair $(\mathfrak g,\sigma)$, Rouvi\`ere's product $\#$ on $\mathrm S(\mathfrak p)^\mathfrak k$ coincides with the product $\star_B$ on $\mathrm H^0(B)=\mathrm S(\mathfrak p)^\mathfrak k$, {\em i.e.}
\[
\beta\!\left(\partial_{\sqrt{j}}(f_1\star_B f_2)\right)=\beta\!\left(\partial_{\sqrt{j}}(f_1)\right)\cdot \beta\!\left(\partial_{\sqrt{j}}(f_2)\right)\ \text{modulo $\mathrm U(\mathfrak g)\cdot\mathfrak k^{-\delta+\frac{1}4\mathrm{tr}_\mathfrak g\circ\mathrm{ad}}$},
\] 
for $f_i$, $i=1,2$, a general element of $\mathrm H^0(B)=\mathrm S(\mathfrak p)^\mathfrak g$.
\end{Thm}
We observe that the characters $\delta$ and $\delta-1/4\ \mathrm{tr}_\mathfrak g\circ \mathrm{ad}$ differ precisely by a character of the Lie algebra $\mathfrak g$ itself (the trace of its adjoint representation) restricted to a character of the Lie subalgebra $\mathfrak h$; as a consequence, $1/4\ \mathrm{tr}_\mathfrak g\circ\mathrm{ad}$ yields an $\mathfrak h$-equivariant map from $\mathfrak g$ to the base field $\mathbb K$ ({\em i.e.} a linear functional on $\mathfrak g$ which vanishes on the subspace $[\mathfrak h,\mathfrak g]$), thanks to which we may actual consider only the character $\delta$, instead of the sum $\delta-1/4\mathrm{tr}_\mathfrak g\circ\mathrm{ad}$.

More precisely, to the symmetric pair $\mathfrak g=\mathfrak k\oplus \mathfrak p$ we may associate a sub-symmetric pair $\mathfrak g_1=\mathfrak k_1\oplus \mathfrak p$, where $\mathfrak k_1=[\mathfrak p,\mathfrak p]$; it is clear that $\mathfrak g_1$ is an ideal of $\mathfrak g$.

The reason is that, to pass from the expression on the right-hand side of the Identity in Theorem~\ref{t-Rou-CF} to the one on the left-hand side, we produce terms which are actually in $\mathrm U(\mathfrak g)\cdot \mathfrak k_1^{-\delta+\frac{1}4\ \mathrm{tr}_\mathfrak g\circ\mathrm{ad}}$, because we reverse two elements in $\mathfrak p$, actually producing an element of $\mathfrak k_1$ (whose commutator with elements of $\mathfrak p$ remains in $\mathfrak p$), and it is obvious that the second summand in the character vanishes on $\mathfrak k_1$, being a character of $\mathfrak g$.  

If we now consider a non-trivial character $\chi$ of $\mathfrak k$, we may associate $X=U_1=\mathfrak g^*$, $U_2=\chi+\mathfrak k^\perp$.
Obviously, $A=\mathrm S(\mathfrak g)$, $B^0=K=\mathrm S(\mathfrak g)/\langle \mathfrak k^{-\chi}\rangle\cong \mathrm S(\mathfrak p)$, the last isomorphism being induced by an affine morphism.

The arguments of Subsubsection~\ref{sss-3-1-1} can be repeated almost {\em verbatim}.
The only relevant difference is that the kernel of the surjective module homomorphism $\mathcal A$ from $(A,\star_A)$ to $(K,\star_L)$ is identified with $A\star_A \mathfrak k^{-\chi-\delta+\frac{1}4\mathrm{tr}_\mathfrak g\circ \mathrm{ad}}$.
Further, the restriction of $\mathcal A$ to $B^0$ and the operator $\mathcal B$ have the same shape as previously.

Still, there is an associative product $\star_B$ on $\mathrm H^0(B)\cong\mathrm S(\mathfrak p)^\mathfrak k$: of course, now the product $\star_B$ depends explicitly on the character $\chi$.
Identity~\eqref{eq-rouviere} is consequently modified as 
\[
(\mathcal A^{-1}\circ \mathcal B)(f_1\star_B f_2)=(\mathcal A^{-1}\circ \mathcal B)(f_1)\star_A (\mathcal A^{-1}\circ \mathcal B)(f_2)\ \text{modulo $A\star_A \mathfrak k^{-\chi-\delta+\frac{1}4\mathrm{tr}_\mathfrak g\circ\mathrm{ad}}$},
\]
from which we deduce 
\[
\beta\!\left(\partial_{\sqrt{j}}(f_1\star_B f_2)\right)=\beta\!\left(\partial_{\sqrt{j}}(f_1)\right)\cdot \beta\!\left(\partial_{\sqrt{j}}(f_2)\right)\ \text{modulo $\mathrm U(\mathfrak g)\cdot\mathfrak k^{-\chi-\delta+\frac{1}4\mathrm{tr}_\mathfrak g\circ\mathrm{ad}}$},
\] 
for any two $\mathfrak k$-invariant elements of $\mathrm S(\mathfrak p)$.
Of course, once again, we may safely remove the character $1/4\ \mathrm{tr}_\mathfrak g\circ\mathrm{ad}$ in the previous identity.

\begin{Rem}\label{r-comm-corr} 
The product $\star_B$ is commutative on $\mathrm S(\mathfrak p)^\mathfrak k$ because of the symmetry of the function $E(X, Y)$, see~\cite[Lemma 11]{CT}, whence the algebra $(\mathrm U(\mathfrak g)/\mathrm U(\mathfrak g)\cdot \mathfrak k^{-\delta})^\mathfrak k$ is also commutative.  
Now the arguments of~\cite[Subsubsection 4.2.1]{CT} are no longer correct because of the contribution of
short loop terms that make the symmetry argument inefficient. 
Therefore, Theorem 5 in~\cite{CT}, which states the commutativity of $(\mathrm U(\mathfrak g)/\mathrm U(\mathfrak g)\cdot \mathfrak k^{-z\delta})^\mathfrak k$, for any real number $z$, is also no longer correct.
\end{Rem}

\subsection{Differential operators expressed {\em via} exponential coordinates in a symmetric pair}\label{ss-3-3}
We consider the triple $X=U_1=\mathfrak g^*$ and $U_2=\mathfrak k^\perp$, for a symmetric pair $(\mathfrak g,\sigma)$.
Therefore, we have two associative algebras $(A,\star_A)$, $(\mathrm H^0(B),\star_B)$, and a $(A,\star_A)$=$(\mathrm H^0(B),\star_B)$-bimodule $K$, where $A=\mathrm S(\mathfrak g)$, $\mathrm H^0(B)=\mathrm S(\mathfrak p)^\mathfrak k$ and $K=\mathrm S(\mathfrak p)$.

Through the function $q$ on $\mathfrak g$, we define the Kashiwara--Vergne density function $D(x,y)$, for $(x,y)$ a general element of $\mathfrak g\times\mathfrak g$, {\em via}
\[
D(x,y)=\frac{\sqrt{q(x)}\sqrt{q(y)}}{\sqrt{q(\mathrm{BCH}(x,y))}},
\]
where the Baker--Campbell--Hausdorff formula $\mathrm{BCH}(x,y)$ is defined by $\exp(x)\exp(y)=\exp(\mathrm{BCH}(x,y))$, for $(x,y)$ in a neighborhood of $(0,0)$.
The Kashiwara--Vergne density function has been introduced in
%~\cite{KV}
to formulate the famous Kashiwara--Vergne conjecture, a general statement for finite-dimensional Lie algebras about deformations of the Baker--Campbell--Hausdorff formula for the product of exponentials in a Lie group: the Kashiwara--Vergne conjecture leads to a proof of Duflo's Theorem about the center of $\mathrm U(\mathfrak g)$ in terms of the $\mathfrak g$-invariant symmetric algebra $\mathrm S(\mathfrak g)^\mathfrak g$.
The Kashiwara--Vergne conjecture has been proved in the general case in
%~\cite{AM}, motivated by earlier results in~\cite{T3}
using deformation quantization techniques to find a suitable deformation of the BCH formula; recently,
%in~\cite{AT},
a different approach to the Kashiwara--Vergne conjecture using Drinfel'd associators has been found, and 
%in~\cite{AT1}
a relationship between the latter approach and the former {\em via} deformation quantization has been elucidated.
We will also discuss a relative Kashiwara--Vergne conjecture in the framework of symmetric pairs later on.

For $x$, $y$ general elements of $\mathfrak p$, we denote by $e^x$ and $e^y$ the exponential of $x$ and $y$, viewed as linear functions on $X=\mathfrak g^*$.
Then, we 
\[
e^x\star_A e^y=\frac{D(x,y)}{D(P(x,y),K(x,y))}e^{P(x,y)}\star_A e^{K(x,y)},\quad P=P(x,y)\in\mathfrak p,\quad K=K(x,y)\in\mathfrak k,
\]
and the power series $P$ and $K$ are defined {\em via} the exponential map for symmetric spaces.
More precisely, it has been proved that symmetric spaces admit an exponential map, which is a diffeomorphism from a neighborhood of $0$ in $\mathfrak p$ into its image in $G/K$, where $G$ is a symmetric pair (in the sense of Lie groups) and $K$ is the fixed point set of an involution $\sigma$ of $G$ (which is a Lie group automorphism): it is simply the restriction of the exponential map of $\mathfrak g$ to right $K$-cosets in $G$.
For $x$, $y$ in a sufficiently small neighborhood of $0$ in $\mathfrak p$, such that $\exp(x)$ and $\exp(y)$ both exist, we have $\exp(x)\exp(x)=\exp(P(x,y))\exp(K(x,y))$.

An important observation at this point is that both $P$ and $K$, as previously defined, are power series in the free Lie algebra generated by $x$, $y$: in particular, the Cartan relations for $(\mathfrak g,\sigma)$ imply that $K$ is an element of $\mathfrak k_1=[\mathfrak p,\mathfrak p]$. 

For $x$, $y$ in a sufficiently small neighborhood of $0$ of $\mathfrak p$ as before, we consider the expression $e^{K(x,y)}\star_L 1=\mathcal A(e^{K(x,y)})$.
First of all, $\mathcal A(e^K)$ is a constant element of $K=\mathrm S(\mathfrak p)$.
By its very definition, $K=K(x,y)$ is an element of $\mathfrak k_1$: in the computation of a summand $K^n\star_L 1$, only the part of the differential operator acting on $K^n$ of degree $n$ survives, because either of degree reasons or of the fact that the restriction of $K$ as a linear function to $\mathfrak k^\perp$ vanishes.
Using the involution $s$ of the preceding Subsection together with the computations of Subsection~\ref{ss-3-1}, it is easy to prove that $K\star_L 1=\delta(K)-1/4\ \mathrm{tr}_\mathfrak g(\mathrm{ad}(K))=\delta(K)$, because, as observed before, $K$ belongs to $\mathfrak k_1$.
We consider, more generally, graphs appearing in the computation of $K^n\star_2 1$, $n\geq 2$.
We may actually repeat almost {\em verbatim} the arguments about the shape of the graphs appearing in $\mathcal A$, $\mathcal B$: the same arguments imply that short loops may appear only at vertices of the first type, which are linked to the only vertex of the second type on the positive real axis corresponding to $K^n$ through a single edge, while more complicated graphs are wheels.
The Cartan relations for the symmetric pair $(\mathfrak g,\sigma)$ imply that the rays of such wheels are colored by $(+,-)$, while the wheel itself has all edges either colored by $(+,+)$ or $(+,-)$.
In particular, it follows that $\mathcal A(e^K)$ consists of an infinite series of wheels, where the short loop graph is considered as the $1$-wheel: in this case, the $1$-wheel contribution appears explicitly.
\begin{Lem}\label{l-sym-K}
For $x$ a general element of $\mathfrak k$, the function $\mathcal A(e^x)=e^x\star_L 1$ satisfies
\[
\mathcal A(e^x)=\sqrt{q(x)}e^{\delta(x)-\frac{1}4\mathrm{tr}_\mathfrak g(\mathrm{ad}(x))};
\]
in particular, if $x$ belongs to the Lie subalgebra $\mathfrak k_1$, the exponent on the right-hand side of the previous equality simplifies to $\delta(x)$.
\end{Lem}
\begin{proof}
To compute $\mathcal A(e^x)$, for $x$ in $\mathfrak k$, it suffices to replace $K$ in the previous computations by $x$.

The rest of the proof follows along the same lines of the proof of~\cite[Lemma 15]{CT}, but we have to observe now that 
\[
\mathcal A(e^x)=e^x\star_L 1=\exp\!\left(\delta(x)-\frac{1}4\mathrm{tr}_\mathfrak g(\mathrm{ad}(x))+\sum_{n\geq 2} w_n^\mathfrak k \mathrm{tr}_\mathfrak k(\mathrm{ad}(x)^n)+\sum_{n\geq 2} w_n^\mathfrak p \mathrm{tr}_\mathfrak p(\mathrm{ad}(x)^n)\right),
\]
for certain integral weights $w_n^\mathfrak k$ and $w_n^\mathfrak p$.
More precisely, such weights are associated to the colored wheel-like graphs
\bigskip
\begin{center}
\resizebox{0.6 \textwidth}{!}{\input{wheel_exp.pstex_t}}\\
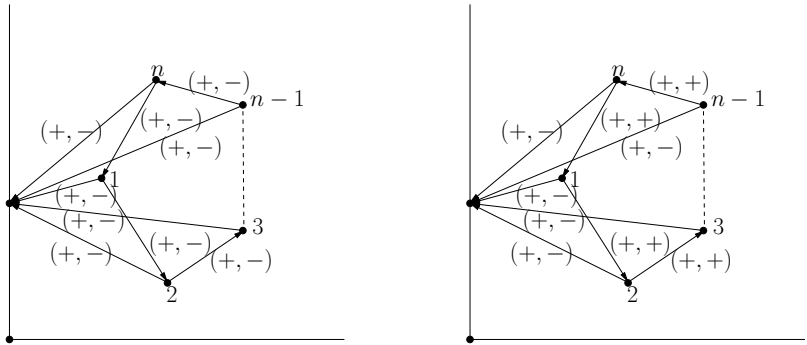
\text{Figure 11 - The possible wheel-like graphs appearing in $\mathcal A(e^x)$, for $x$ in $\mathfrak k$}
\end{center}
\bigskip
Further, also $\sqrt{q(x)}$ can be written in a similar form, with the only difference that it does not contain terms proportional to the trace of the adjoint representations of $\mathfrak k$ on itself or on $\mathfrak p$.

Therefore, the very same computations of~\cite[Lemma 15]{CT} imply the above identity.
\end{proof}
We observe that we may get rid of the factor $e^{\delta-\frac{1}4\mathrm{tr}_\mathfrak g\circ\mathrm{ad}}$, viewed as an element of the completion $\widehat{\mathrm S}(\mathfrak k)$ simply by applying the biquantization techniques to the modified triple $X=U_1=\mathfrak g$, $U_2=-\delta+1/4\ \mathrm{tr}_\mathfrak g\circ\mathrm{ad}+\mathfrak k^\perp$.
To the previous triple are associated two associative algebra $(A,\star_A)$ and $(\mathrm H^0(B),\star_B)$ and a bimodule $K$, where $A=\mathrm S(\mathfrak g)$, $B^0=K=\mathrm S(\mathfrak p)$, and $H^0(B)=\mathrm S(\mathfrak p)^\mathfrak k$.
Notice that the product $\star_B$ on $\mathrm H^0(B)$, for $B$ defined by the modified triple, does not coincide with $\star_B$ for the initial triple; similarly, the operator $\mathcal A$ on $A$ for the modified triple also does not coincide with the operator $\mathcal A$ for the initial triple, as their kernels do not obviously coincide.
Still, both their restrictions to $B^0$ coincide, thus also their symbols, and similarly for $\mathcal B$.
On the other hand, as already observed, for the above modified triple we have the identity
\[
\mathcal A(e^x)=e^x\star_L 1=\sqrt{q(x)},\ x\in\mathfrak k.
\]

We consider the following expression, using the notation from above, 
\[
\begin{aligned}
(e^x\star_A e^y)\star_L 1&=\frac{D(x,y)}{D(P,K)} \left(e^P\star_A e^K\right)\star_L 1=\frac{\sqrt{q(x)}\sqrt{q(y)}}{\sqrt{q(P)}\sqrt{q(K)}} e^P\star_L(e^K\star_L 1)=\frac{\sqrt{q(x)}\sqrt{q(y)}}{\sqrt{q(P)}}j_\mathcal A(P)e^P=\\
&=\frac{\sqrt{q(x)}\sqrt{q(y)}}{\sqrt{j(P)}}j_\mathcal B(P)e^P=\\
&=e^x\star_L(e^y\star_L1)=j_\mathcal A(y) e^x\star_L e^y=\frac{j_\mathcal A(y)}{j_\mathcal B(y)} e^x\star_L (1\star_R e^y),\ x,y\in\mathfrak p.
\end{aligned}
\]
Comparing the expression on the second line with the rightmost expression on the third line, we find
\[
e^x\star_L (1\star_R e^y)=\frac{\sqrt{q(x)}\sqrt{j(y)}j_\mathcal B(P)}{\sqrt{j(P)}} e^P
\]
We may view the expressions on both sides of the previous identity as analytic functions on a sufficiently small neighborhood of $(0,0)$ in $\mathfrak p\times \mathfrak p$.

We consider further an element $T$ of $\mathrm H^0(B)=\mathrm S(\mathfrak p)^\mathfrak k$, for which we obtain
\[
\begin{aligned}
e^x\star_L (1\star_R T)&=(e^x\star_L 1)\star_R T=j_\mathcal A(x) e^x\star_R T=\\
&=T_y\!\left(e^x\star_L (1\star_R e^y)\right)|_{y=0}=T_y\!\left(\frac{\sqrt{q(x)}\sqrt{j(y)}j_\mathcal B(P)}{\sqrt{j(P)}} e^P\right)\bigg\vert_{y=0},
\end{aligned}
\]
where we regard in the second line $T$ as a differential operator on $\mathfrak p^*$ with respect to the variable $y$.

Therefore, using Identity~\eqref{eq-DR-symb}, we get the following expression,
\[
e^x\star_R T=T_y\!\left(\frac{\sqrt{j(x)}\sqrt{j(y)}j_\mathcal B(P)}{\sqrt{j(P)}j_\mathcal B(P)} e^P\right)\bigg\vert_{y=0},
\]
whose right-hand side is, according to the arguments of~\cite[Section 6]{Rouv}. precisely the differential operator $\beta(\partial_{sqrt{j}}(T))$ expressed {\em via} exponential coordinates on the symmetric space $G/K$, for $G$, $K$ connected, simply connected Lie groups with Lie algebras $\mathfrak g$, $\mathfrak k$ respectively, up to a modification by the analytic function $\sqrt{j}/j_\mathcal B$ on $\mathfrak p$.
A deep consequence of the fact that the product $\star_B$ coincides with Rouvi\`ere's product, together with a result about the existence of polarizations compatible with the structure of symmetric pair, for which we refer to~\cite{T1,T2}, implies that the symbol $j_\mathcal B$ is constant and thus equal to $1$.

Therefore, the expression $e^x\star_L T$, for $T$ in $\mathrm S(\mathfrak p)^\mathfrak k$, viewed here as an element of $K$, truly expresses in terms of biquantization the differential operator $T$ through exponential coordinates on $G/K$.

\subsection{Deformation of the Baker--Campbell--Hausdorff formula for symmetric pairs}\label{ss-3-4}
As already briefly remarked in the previous Subsection, the problem of deforming the BCH formula and the BCH density function $D$ (see above) for a finite-dimensional Lie algebra $\mathfrak g$ is related to Duflo's Theorem {\em via} the Kashiwara--Vergne conjecture.
As the KV conjecture relies on the exponential map on Lie algebras, it seems natural to formulate a similar conjecture for a general symmetric pair, because also symmetric pairs admit an exponential map.
We refer to~\cite{Rouv} for more details on the KV conjecture for symmetric pairs.
%%%% ancora una?

The exponential map for a symmetric pair $(\mathfrak g,\sigma)$ with Cartan decomposition $\mathfrak g=\mathfrak k\oplus\mathfrak p$ is a well-defined diffeomorphism $\exp_{G/K}$ from a sufficiently small neighborhood of $0$ in $\mathfrak p$ to its image in $G/K$, for $G$, $K$ as in the previous Subsection: it is induced from the exponential map $\exp_G$ of $\mathfrak g$.
Standard manipulations, see~\cite[Section 2]{Rouv} for more details, imply that the exponential map for the symmetric pair $(\mathfrak g,\sigma)$ yields a BCH formula $\mathrm{BCH}_\mathfrak p$ {\em via}
\[
\exp_G(2\mathrm{BCH}_\mathfrak p(x,y))=\exp_G(x)\exp_G(2y)\exp_G(x),
\]
for $(x,y)$ in a sufficiently small neighborhood of $(0,0)$ in $\mathfrak p\times\mathfrak p$.
The element $\mathrm{BCH}_\mathfrak p(x,y)$ belongs to $\mathfrak p$, and is expressible in terms of even iterated brackets of $x$ and $y$.
  
The KV conjecture for symmetric pairs expresses a deformation of the function $\mathrm{BCH}_\mathfrak p$ in terms of $\mathfrak k$-adjoint vector fields on $\mathfrak p\times\mathfrak p$; we refer once again to~\cite{Rouv} for a complete introduction to this issue and for a discussion of some consequences.
There is also another claim, which expresses the corresponding variation in terms of the said $\mathfrak k$-adjoint vector fields of the density function for the symmetric pair, which is defined as
\[
D_\mathfrak p(x,y)=\frac{\sqrt{j(x)}\sqrt{j(y)}}{\sqrt{j(\mathrm{BCH}_\mathfrak p(x,y))}},
\]
again for $(x,y)$ in a sufficiently small neighborhood of $(0,0)$, where all functions make sense.
The deformation of $D_\mathfrak p$ contains traces of the adjoint $\mathfrak k$-action on $\mathfrak p$.

We are now going to illustrate how biquantization techniques yield two different deformations of $\mathrm{BCH}_\mathfrak p$ and $D_\mathfrak p$ in the sense elucidated above.
These two deformations can be characterized in terms of the $A_\infty$-structures on $A$, $B$ and $K$ through the quadratic relations between the corresponding Taylor coefficients: this is not immediately recognizable from the approach we take, which is in turn essentially motivated by the results of~\cite{T3,AM}.

We consider the triple $X=U_1=\mathfrak g^*$ and $U_2=\delta-1/4\ \mathrm{tr}_\mathfrak g\circ\mathrm{ad}+\mathfrak k^\perp$ and the corresponding algebras $A$, $\mathrm H^0(B)$ and $A$-$\mathrm H^0(B)$-bimodule $K$.

\subsubsection{First deformation}\label{sss-3-4-1}
First of all, we may consider, for $(x,y)$ in a sufficiently small neighborhood of $(0,0)$ in $\mathfrak p\times\mathfrak p$, where $\mathrm{BCH}_\mathfrak p$ and $D_\mathfrak p$, as well as $\sqrt{q}$, $\sqrt{j}$, $j_\mathcal A$ and $j_\mathcal B$ are well-defined, the function
\[
\begin{aligned}
(e^x\star_A e^y)\star_L 1&=\frac{\sqrt{q(x)}\sqrt{q(y)}}{\sqrt{j(\mathrm{BCH}_\mathfrak p(x,y))}} e^{\mathrm{BCH}_\mathfrak p(x,y)}=\\
&=e^x\star_L(e^y\star_L 1)=j_\mathcal A(y) e^x\star_L e^y,
\end{aligned}
\] 
where we have used results of the previous Subsection, setting $P(x,y)=\mathrm{BCH}_\mathfrak p(x,y)$.

On the other hand, we may also consider 
\[
e^x\star_L(1\star_R e^y)=j_\mathcal B(y) e^x\star_L e^y=e^x\star_L e^y,
\]
where we have used once again the aforementioned fact that $j_\mathcal B\equiv 1$, whence, recalling Identity~\eqref{eq-DR-symb},
\[
e^x\star_L e^y=\frac{\sqrt{q(x)}\sqrt{j(y)}}{\sqrt{j(\mathrm{BCH}_\mathfrak p(x,y))}} e^{\mathrm{BCH}_\mathfrak p(x,y)}
\]
Finally, we also have
\[
(e^x\star_L 1)\star_R e^y=j_\mathcal A(x) e^x\star_R e^y.
\]

The problem is therefore, how to relate $(e^x\star_L 1)\star_R e^y$ with $e^x\star_L (1\star_R e^y)$ in order to draw a bridge between the previous two formul\ae: the two expressions do not coincide, as both actions are compatible to each other only in cohomology.
But the $A_\infty$-nature of $B$ and of the bimodule $K$ permit to control explicitly the failure for the compatibility between $\star_L$ and $\star_R$: in facts, we find
\begin{equation}\label{eq-deform-1}
(e^x\star_L 1)\star_R e^y-e^x\star_L(1\star_R e^y)=\mathrm d_K^{1,1}(e^x,1,\mathrm d_B(e^y))=\mathrm d_K^{1,1}(e^x,1,e^y\mathrm d_B(y))
\end{equation}
where $\mathrm d_B$ is the Chevalley--Eilenberg differential on the complex $B=\mathrm S(\mathfrak p)\otimes \wedge(\mathfrak k^*)$, and $\mathrm d^{1,1}_K$ denotes the $(1,1)$-Taylor component of the $A_\infty$-bimodule structure on $K$.
More precisely,  
\[
\mathrm d_K^{1,1}(a_1|k|b_1)=\sum_{n\geq 0}\frac{1}{n!}\sum_{\Gamma\in\mathcal G_{n,3}}\mu_{n+3}^K\left(\int_{\mathcal C_{n,3}^+}\prod_{e\in E(\Gamma)}\omega^K_e(\underset{n}{\underbrace{\pi|\cdots|\pi}}|a_1|k|b_1)\right),
\]
and dimensional arguments imply that $b_1$ must be an element of $B^1$, otherwise the previous expression is trivial.

Identity~\eqref{eq-deform-1} may be derived in a slightly different way, which makes apparent the fact that it embodies the KV deformation problem sketched at the beginning of the present Subsection.
Namely, for $n\geq 1$, we consider the forgetful projection $\pi_{n,1,1}$ from $\mathcal C_{n,1,1}^+$ to $\mathcal C_{0,1,1}^+$: in this situation, we prefer to consider the compactified configuration spaces $\mathcal C_{n,1,1}^+$ to highlight the fact that we consider the functions $e^x$, $1$ and $e^y$ to be put on the positive imaginary axis, at the origin and on the positive real axis respectively.
We observe that the fiber of $\pi_{n,1,1}$ at a generic point of $\mathcal C_{0,1,1}^+$ is an orientable compact smooth manifold with corners of dimension $2n$, hence we may consider the push-forward (or integration along the fiber) $\pi_{n,1,1,*}$ with respect to $\pi_{n,1,1}$.

In particular, we may consider the expression
\[
\sum_{n\geq 0}\frac{1}{n!}\sum_{\Gamma\in\mathcal G_{n,3}}\mu_{n+3}^K\left(\pi_{n,1,1,*}\!\left(\prod_{e\in E(\Gamma)}\omega^K_e\right)(\underset{n}{\underbrace{\pi|\cdots|\pi}}|e^x|1|e^y)\right),
\]
for $(x,y)$ as above.

First of all, for $(x,y)$ as above, the previous expression is a smooth function on $\mathcal C_{0,1,1}^+$: it is a consequence of the fact that, for $n\geq 1$,  the push-forward $\pi_{n,1,1,*}$ selects the piece of the integrand of (form) degree bigger or equal than $2n$.
To the three vertices of the second type are associated functions, while to each vertex of the first type of an admissible graph $\Gamma$ of type $(n,3)$ is associated a copy of the linear Poisson bivector $\pi$: as a consequence, the form degree of each integrand must be precisely $2n$, whence the first claim.
We may further divide the previous expression by $j_\mathcal A(x)$: this ``normalization'' is due to previous computations, and it does not affect the following computations.

Because of the fact that $\pi$ is a linear Poisson bivector, we may use the arguments in~\cite[Chapter 2]{BCKT} or~\cite{Kath} to prove that the ``normalized'' function on $\mathcal C_{0,1,1}^+$ given by the previous expression can be re-written in the form
\begin{equation}\label{eq-1-def}
D_\mathfrak p^1(x,y)e^{\mathrm{BCH}_\mathfrak p^1(x,y)},
\end{equation}
where the exponent $\mathrm{BCH}_\mathfrak p^1(x,y)$ is a smooth $\mathfrak p$-valued function on $\mathcal C_{0,1,1}^+$, corresponding to the connected graphs of Lie type, whose unique root is in $\mathfrak p$, while $D_\mathfrak p^1(x,y)$ is a smooth $\mathbb K$-valued function on $\mathcal C_{0,1,1}^+$, corresponding to graphs of Lie type with roots in $\mathfrak k$ and wheel-like graphs (possibly with graphs of type Lie attached to their spokes).
Both $\mathrm{BCH}_\mathfrak p^1(x,y)$ and $D_\mathfrak p^1(x,y)$ are weighted sums over the graphs highlighted right above, where now the corresponding integral weights are smooth functions on $\mathcal C_{0,1,1}^+$.

The deformation formul\ae\ we are interested into can be therefore computed by taking the exterior derivative of~\eqref{eq-1-def} as a function on $\mathcal C_{0,1,1}^+$: we may therefore apply the generalized Stokes Theorem for integration along the fiber of $\pi_{n,1,1}$ in the first integral formula for~\eqref{eq-1-def}.
For any admissible graph $\Gamma$ of type $(n,3)$ as above, the corresponding integrand is a closed form, whence it suffices to consider the corresponding integral along the boundary strata of codimension $1$ of the generic fiber.
For $n\geq 1$, a general boundary stratum of codimension $1$ in the generic fiber corresponds either $i)$ to the collapse of points in $Q^{+,+}$ labeled by a subset $A$ of $[n]$ of cardinality $2\leq |A|\leq n$  to a single point in $Q^{+,+}$, or $ii)$ to the approach of points in $Q^{+,+}$ labeled by a subset $A$ of $[n]$ of cardinality $0\leq |A|\leq n$ to $i\mathbb R^+$, to the origin or to $\mathbb R^+$.
Notice that no boundary stratum appears, where either the point on $i\mathbb R^+$ or $\mathbb R^+$ approaches the origin (this is because we are considering the boundary strata of codimension $1$ of the generic fiber). 

Standard dimensional arguments, Kontsevich's Vanishing Lemma~\cite[Lemma]{K} and the boundary conditions for the $4$-colored propagators $\omega^{+,+}$ and $\omega^{+,-}$ imply that the only boundary strata yielding non-trivial contributions correspond to boundary strata of type $ii)$, where points in $Q^{+,+}$ labeled by $A\subseteq [n]$ approach the point on $\mathbb R^+$.
We refer to~\cite{T3} or~\cite[Chapter 2]{BCKT}
for similar computations.
 
The sum over all such contributions yields a smooth $\widehat{\mathrm S}(\mathfrak p)$-valued $1$-form on $\mathcal C_{0,1,1}^+$, depending on $(x,y)$ as above, whose integral over $\mathcal C_{0,1,1}^+$ identifies with the right-hand side of Identity~\eqref{eq-deform-1}.

The generalized Stokes Theorem, see {\em e.g.} the computations in~\cite{T3}, implies that the previous $1$-form specifies a smooth $\mathfrak k$-valued $1$-form $\omega_1(x,y)$ satisfying the identities
\[
\begin{aligned}
\mathrm d\mathrm{BCH}_\mathfrak p^1(x,y)&=\langle \left[y,\omega_1(x,y)\right],\partial_y\mathrm{BCH}_\mathfrak p^1(x,y)\rangle,\\
\mathrm d D_\mathfrak p^1(x,y)&=\langle \left[y,\omega_1(x,y)\right],\partial_y D_\mathfrak p^1(x,y)\rangle+\mathrm{tr}_\mathfrak p\!\left(\mathrm{ad}(y)\partial_y\omega_1(x,y)\right)D_\mathfrak p^1(x,y),
\end{aligned}
\]
where we have used the notation
\[
\langle [y,\omega_1(x,y)],\partial_y\mathrm{BCH}_\mathfrak p^1(x,y)\rangle=\frac{\mathrm d}{\mathrm d t}\mathrm{BCH}_\mathfrak p^1(x,y+t[y,\omega_1(x,y)])\big\vert_{t=0},
\]
and analogously for other similar expressions in the previous identities.

Finally, we consider the function~\eqref{eq-1-def} on $\mathcal C_{0,1,1}^+$, which is a smooth, compact $1$-dimensional manifold with corners: its two boundary strata of codimension $1$ correspond to the approach of the point either on $i\mathbb R^+$ or on $\mathbb R^+$ to the origin.
If we choose the smooth section of $\mathcal C_{0,1,1}^+$ which corresponds to fixing the point on $i\mathbb R^+$ to $i$, then $\mathcal C_{0,1,1}^+\cong [0,\infty]$: the boundary point $\{0\}$, resp.\ $\{\infty\}$, corresponds to the approach of the point on $\mathbb R^+$, resp.\ on $i\mathbb R^+$, to the origin.
Taking into account the ``normalization'' with respect to $j_\mathcal A(x)$ in the function~\eqref{eq-1-def} and the computations at the beginning of the Subsubsection, the value of the function~\eqref{eq-1-def} at $0$ yields the values of both $\mathrm{BCH}_\mathfrak p^1$ and $D_\mathfrak p^1(x,y)$ at $0$, which are precisely $\mathrm{BCH}_\mathfrak p(x,y)$ and $D_\mathfrak p(x,y)$, whence the function~\eqref{eq-1-def} produces a genuine deformation of the BCH formula and corresponding density for the symmetric pair $(\mathfrak g,\sigma)$.
The value of the ``normalized'' function~\eqref{eq-1-def} at $\infty$ is also readily computed, namely it is simply $e^x\star_R e^y$.

\subsubsection{Second deformation}\label{sss-3-4-2}
We begin by considering, for $(x,y)$ as in the previous Subsubsection, the two expressions $(1\star_R e^x)\star_R e^y$ and $1\star_R(e^x\star_B e^y)$.
We observe that neither the pairing $\star_B$ is associative, nor it is compatible with the pairing $\star_R$ in both expressions.
Since, as observe before, $j_\mathcal B\equiv 1$, the first expression equals $e^x\star_R e^y$, while the second equals $e^x\star_B e^y$.

The fact that $\star_R$ and $\star_B$ are the $(0,1)$- and $2$-Taylor components of the $A_\infty$-structures on $K$ and on $B$, we get
\begin{equation}\label{eq-deform-2}
(1\star_R e^x)\star_R e^y-1\star_R(e^x\star_B e^y)=e^x\star_R e^y-e^x\star_B e^y=\mathrm d_K^{1,1}(1|\mathrm d_B(e^x)|e^y)+\mathrm d_K^{1,1}(1|e^x|\mathrm d_B(e^y)).
\end{equation}

We observe that, due to the choice of the triple $X=U_1=\mathfrak g^*$ and $U_2=\delta-1/4\ \mathrm{tr}_\mathfrak g\circ\mathrm{ad}+\mathfrak k^\perp$, the pairing $\star_B$ depends on the (natural) character $\delta-1/4\ \mathrm{tr}_\mathfrak g\circ\mathrm{ad}$.

We now define, for $(x,y)$ as above, a smooth function on $\mathcal C_{0,0,2}^+$ {\em via}
\[
\sum_{n\geq 0}\frac{1}{n!}\sum_{\Gamma\in\mathcal G_{n,3}}\mu_{n+3}^K\left(\pi_{n,0,2,*}\!\left(\prod_{e\in E(\Gamma)}\omega^K_e\right)(\underset{n}{\underbrace{\pi|\cdots|\pi}}|1|e^x|e^y)\right),
\]
where, for $n\geq 1$, $\pi_{n,0,2}$ is the forgetful projection from $\mathcal C_{n,0,2}^+$ onto $\mathcal C_{0,0,2}^+$, and $\pi_{n,0,2,*}$ denotes the push-forward with respect to $\pi_{n,0,2}$.
Again, because of the fact that $\pi$ is a linear Poisson bivector, the arguments of~\cite[Chapter 2]{BCKT} or~\cite{Kath} yield an explicit expression for the previous function on $\mathcal C_{0,0,2}^+$ in the shape
\begin{equation}\label{eq-2-def}
D_\mathfrak p^2(x,y)e^{\mathrm{BCH}_\mathfrak p^2(x,y)}.
\end{equation}
The $\mathfrak p$-valued function $\mathrm{BCH}_\mathfrak p^2(x,y)$ and the $\mathbb K$-valued function $D_\mathfrak p^2(x,y)$ are as in Formula~\eqref{eq-1-def}, with due modifications in the integral weights.

As in the previous Subsubsection, the exterior derivative of the function~\eqref{eq-2-def} may be computed by means of the generalized Stokes Theorem in a way similar to the sketch of the computation of the exterior derivative of the function~\eqref{eq-1-def}.
The relevant contributions come from the boundary strata of codimension $1$ of the generic fiber of $\pi_{n,0,2}$, for $n\geq 1$: such strata correspond to either $i)$ the collapse of points in $Q^{+,+}$ labeled by a subset $A$ of $[n]$ of cardinality $2\leq |A|\leq n$ to a single point in $Q^{+,+}$, or $ii)$ the approach of points in $Q^{+,+}$ labeled by $A\subseteq [n]$, $0\leq |A|\leq n$, to $i\mathbb R^+$, the origin, or $\mathbb R^+$.
As we are considering a generic fiber, we do not consider strata, where the two points on $\mathbb R^+$ approach to each other or where the first point on $\mathbb R^+$ approaches the origin.
Standard arguments imply that the only non-trivial contributions come from boundary strata of type $ii)$, corresponding to the approach of points in $Q^{+,+}$ either to the first or to the second point on $\mathbb R^+$.

More precisely, the sum over all such contributions yields a smooth $\widehat{\mathrm S}(\mathfrak p)$-valued $1$-form on $\mathcal C_{0,0,2}^+$, whose integral over $\mathcal C_{0,0,2}^+$ is precisely the rightmost expression in the chain of identities~\eqref{eq-deform-2}.
The previous $1$-form yields in turn a $\mathfrak k\times\mathfrak k$-valued $1$-form on $\mathcal C_{0,0,2}^+$ $\omega_2(x,y)=(\omega_2^1(x,y),\omega_2^2(x,y))$, which obeys the identities
\[
\begin{aligned}
\mathrm d\mathrm{BCH}_\mathfrak p^2(x,y)&=\langle \left[x,\omega_2^1(x,y)\right],\partial_x\mathrm{BCH}_\mathfrak p^2(x,y)\rangle+\langle \left[y,\omega_2^2(x,y)\right],\partial_y\mathrm{BCH}_\mathfrak p^2(x,y)\rangle,\\
\mathrm d D_\mathfrak p^1(x,y)&=\langle \left[x,\omega_2^1(x,y)\right],\partial_x D_\mathfrak p^2(x,y)\rangle+\langle \left[y,\omega_2^2(x,y)\right],\partial_y D_\mathfrak p^2(x,y)\rangle+\\
&\phantom{=}+\mathrm{tr}_\mathfrak p\!\left(\mathrm{ad}(x)\partial_x\omega_2^1(x,y)+\mathrm{ad}(y)\partial_y\omega_2^2(x,y)\right)D_\mathfrak p^2(x,y).
\end{aligned}
\]

The previous results imply due modifications of the results of~\cite[Subsubsection 4.4.3]{CT}: we observe that changes are caused by the fact that we have chosen at the beginning the modified triple $X=U_1=\mathfrak g^*$ and $U_2=\delta-1/4\ \mathrm{tr}_\mathfrak g\circ\mathrm{ad}+\mathfrak k^\perp$, thus introducing in many formul\ae\ the (natural) character $\delta-1/4\ \mathrm{tr}_\mathfrak g\circ\mathrm{ad}$ of $\mathfrak k$.

%%%%%%%%%%%%%%%%%%%%%%%%%%%%%%%%%%%%%%%

\section{Harish-Chandra homomorphism in diagrammatical terms re-visited}\label{s-5}
Once again, we borrow notation and conventions from~\cite[Section 5]{CT}, in particular for what concerns the generalized Iwasawa decomposition.
Geometrically, to the decomposition $\mathfrak g=\mathfrak k\oplus\mathfrak p_0\oplus\mathfrak n_+$, and $\mathfrak k=\mathfrak k_0\oplus \mathfrak k/\mathfrak k_0$, we associate $U_1=\mathfrak k^\perp$ and $U_2=(\mathfrak k_0\oplus\mathfrak n_+)^\perp$, whence
\[
A=\mathrm S(\mathfrak p_0)\otimes\mathrm S(\mathfrak n_+)\otimes\wedge(\mathfrak k_0^*)\otimes\wedge(\mathfrak r^*),\ B=\mathrm S(\mathfrak p_0)\otimes\mathrm S(\mathfrak r)\otimes\wedge(\mathfrak k_0^*)\otimes\wedge(\mathfrak n_+^*),\ K=\mathrm S(\mathfrak p_0)\otimes\wedge(\mathfrak k_0^*),
\]
where $\mathrm r=\mathfrak k/\mathfrak k_0$.

We observe that, in this situation, we have all $4$-colored propagators coming into play; thus, in admissible graphs, multiple edges may appear.

\subsection{Harish-Chandra graphs and reduction spaces re-visited}\label{ss-5-1} 
The $A_\infty$-$A_\hbar$-$B_\hbar$-bimodule $K_\hbar$ identifies now, as a vector space, with $\mathrm S(\mathfrak p_0)\otimes\wedge(\mathfrak k_0^*)[\![\hbar]\!]$, and has a differential $\mathrm d_{K_\hbar}^{0,0}$, due to the flatness of both $A_\hbar$, $B_\hbar$.

The next proposition is new, and we need it to identify correctly the $0$-th cohomology of $\mathrm d_{K_\hbar}^{0,0}$: this result is implicitly used in~\cite{CT}, but we have realized that its proof needs a more involved argument which requires Lemma~\ref{l-vanish-4}.
\begin{Prop}\label{p-red-bimod}
The reduction space $\mathrm H^0(K_\hbar)$ identifies with $\mathrm S(\mathfrak p_0)^{\mathfrak k_0}[\![\hbar]\!]$. 
\end{Prop}
\begin{proof}
By its very construction, $\mathrm d_{K_\hbar}^{0,0}$ on $\mathrm S(\mathfrak p_0)$ is determined by admissible graphs in $\mathcal G_{n,1}$; we observe that the only vertex of the second type is the origin.
We consider an admissible graph $\Gamma$ in $\mathcal G_{n,1}$, for $n\geq 2$: there is a vertex of the first type, from which departs one edge to $\infty$, and such an edge is colored by $(-,-)$.
We concentrate on the remaining $n-1$ vertices: every edge hitting the origin is colored by $(+,+)$, therefore, from each one of the $n-1$ vertices of the first type can depart at most one edge to the origin.
We denote by $p$ the number of edges departing from the $n$ vertices of the first type of $\Gamma$ and hitting the origin: by the previous arguments, we know that $p\leq (n-1)+1=n$, because the vertex with an edge to $\infty$ has an additional edge, which may or may not hit the origin.
The polynomial degree of the differential operator associated to $\Gamma$ equals $n-(2n-1-p)=-n+1+p$, which must be greater or equal than $0$: this forces immediately $p\geq n-1$.
This fact, combined with the previous condition on $p$, yields that either $p=n$ or $p=n-1$.

In the case $p=n$, the corresponding differential operator has polynomial degree $1$, and from each vertex of the first type departs exactly one edge to the origin: this implies that the graph $\Gamma$ must have a vertex of the form 
\bigskip
\begin{center}
\resizebox{0.3 \textwidth}{!}{\input{Iwasa_triv.pstex_t}}\\
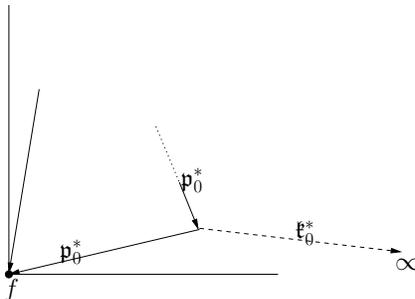
\text{Figure 12 - The only possible configuration at a bivalent vertex} \\
\end{center}
\bigskip
{}Of course, the vertex from which departs the arrow to $\infty$ may be also isolated, {\em i.e.} no edge has it as the endpoint: in this case, as $n\geq 2$, standard dimensional arguments imply that the corresponding weight is trivial.
If it is not isolated, the property $[\mathfrak k_0,\mathfrak p_0]\subset \mathfrak p_0$ implies that the two consecutive edges correspond to propagators of type $(+,+)$, and Lemma~\ref{l-vanish-4} yields triviality of the corresponding weight.

We now consider the case $p=n-1$: the corresponding differential operator has polynomial degree $0$, {\em i.e.} it is a translation-invariant differential operator. 
We first assume that the vertex from which departs an edge to $\infty$ has the other edge hitting the origin: the fact that the differential operator has constant coefficients implies that this vertex cannot be isolated, otherwise Lemma~\ref{l-vanish-4} would imply triviality of its weight.
Therefore, $\Gamma$ must be a disjoint union of wheel-like graphs, one of which must look like as follows:
\bigskip
\begin{center}
\resizebox{0.3 \textwidth}{!}{\input{wheel_infty.pstex_t}}\\
\text{Figure 13 - A wheel with an edge to $\infty$} \\
\end{center}
\bigskip

Lemma~\ref{l-vanish-4} implies that the two edges meeting at the vertex with the edge to $\infty$ must have distinct colors, otherwise the corresponding weight vanishes.
Therefore, the relations $[\mathfrak k_0,\mathfrak k_0]\subseteq \mathfrak k_0$, $[\mathfrak k_0,\mathfrak p_0]\subseteq \mathfrak p_0$, $[\mathfrak k_0,\mathfrak n_+]\subseteq \mathfrak n_+$ imply immediately that the outgoing edge from this special vertex admits only the color $\mathfrak r$.
But then again, the previous relations imply that all edges in the cycle of the wheel-like graph are colored by $\mathfrak r$, hence we may apply once again Lemma~\ref{l-vanish-4}.

We have thus proved that the only non-trivial graph corresponds exactly to the adjoint action of $\mathfrak k_0$ on $\mathrm S(\mathfrak p_0)$, whence the claim follows. 
\end{proof} 
Therefore, the three reduction algebras associated to $A$, $B$ and $K$, are exactly as in~\cite[Subsection 5.1]{CT}.

The discussion in~\cite[Subsection 5.2]{CT} is not modified by the previously discussed changes: in fact, it deals only the reduction algebra placed on the positive imaginary axis, which is left untouched by the changes occurring in biquantization.

\section{Construction of characters re-visited}\label{s-6}
As before, we consider a symmetric pair $(\mathfrak g,\sigma)$, with a Cartan decomposition $\mathfrak g=\mathfrak k\oplus\mathfrak p$.

The techniques of biquantization are applied, in this framework, to the case of two coisotropic subvarieties $f+\mathfrak k^\perp$ and $\mathfrak k^\perp$, where $f$ is an element of $\mathfrak k^\perp$, and $\mathfrak b$ is a polarization for $f$ ({\em i.e.} $f$ determines a skew-symmetric form on $\mathfrak g$, and a polarization $\mathfrak b$ for $f$ is a subalgebra of $\mathfrak g$, which is isotropic for the said skew-symmetric form and of maximal dimension among the isotropic subalgebras).
Thus, we observe that in this situation, too, all $4$-colored propagators are involved.

The aim of~\cite[Subsections 2.5 and 6.1]{CT} is the construction of characters for the reduction algebra $(\mathrm H^0_{\hbar,\mathfrak b}(\mathfrak k^\perp),\star_\mathrm{CF})$, where $\mathrm H^0_{\hbar,\mathfrak b}(\mathfrak k^\perp)$ denotes the $0$-th cohomology of $A_\hbar$, where $A=\mathcal O_{\mathrm N^\vee_{\mathfrak k^\perp,\mathfrak g^*}[-1]}$.
In order to perform explicit computations using biquantization, a decomposition of $\mathfrak g$ into direct summands, two of them being complements of $\mathfrak k$  and $\mathfrak b$, is needed, because in this situation, all $4$-colored propagators are involved: hence, the suffix denotes an explicit dependence (in the explicit computations) of the polarization $\mathfrak b$.
As has been proved in detail in~\cite[Subsection 1.5]{CT}, there is an explicit (non canonical) isomorphism $\mathrm H^0(A_\hbar)\cong \mathrm H^0_{\hbar,\mathfrak b}(\mathfrak k^\perp)$: thus, {\em via} biquantization, it is possible to construct characters for $\mathrm H^0(A_\hbar)$, which still depend from a choice of a polarization $\mathfrak b$ of some element $f$ in $\mathfrak k^\perp$.

Therefore, to really deal with a canonical construction of characters for $H^0(A_\hbar)$ {\em via} biquantization, one needs to prove independence of the choice of polarizations of the characters constructed on $H^0_{\hbar,\mathfrak b}(\mathfrak k^\perp)$.
The main technical tool in the proof of independence of polarizations is a Stokes' argument, reminiscent of the Stokes' argument leading to biquantization, in presence of three coisotropic submanifolds, {\em i.e.} for $f$ in $\mathfrak k^\perp$ as before, $f+\mathfrak b_i^\perp$, $i=1,2$, and $\mathfrak k^\perp$, where now $\mathfrak b_i$ is a polarization for $f$, $i=1,2$.

Of course, in this new situation, we need the $8$-colored propagators $\theta_{j_1,j_2,j_3}$, whose explicit construction and relative discussion of the main properties is the content of~\cite[Subsubsection 6.2.1]{CT}.

Roughly speaking, the $8$-colored propagators interpolate the $4$-colored ones (in a sense that will be made precise later on): therefore, as already pointed out in the Introduction, the appearance of ``regular terms'' in the $4$-colored propagators is likely to cause the appearance of similar ``regular terms'' in the $8$-colored propagators.
This is the main novelty in the present discussion.

\subsection{Construction of the $8$-colored propagators re-visited}\label{ss-6-1}
We will write down in more detail the construction outlined in~\cite[Subsubsection 6.2.1]{CT}, from which we borrow notation and conventions.

The $8$-colored propagators $\theta_{j_1,j_2,j_3}$, $j_k\in\{1,2\}$, $k=1,2,3,$, are $8$ distinct smooth, closed $1$-forms on the compactified configuration space $\mathcal C^+_{2,0,0,0}(\sqsubset)$, where $\sqsubset=\left\{z=x+\mathrm i y\in \mathbb C:\ x\geq 0,\ -\frac{\pi}2\leq y\leq \frac{\pi}2\right\}$.
We observe that the half-strip $\sqsubset$ is a smooth manifold with corners of dimension $2$ with three boundary components $\sqsubseteq_i$, $i=1,2,3$ of codimension $1$, where $\sqsubseteq_1$ is the lower horizontal half-line, $\sqsubseteq_2$ is the vertical segment and $\sqsubseteq_3$ is the upper horizontal half-line.
The space $\mathcal C^+_{2,0,0,0}(\sqsubset)$, which is a smooth manifold with corners of dimension $4$, is the correct generalization in the framework of $3$ branes of Kontsevich's eye $\mathcal C_{2,0}^+$ (the compactified configuration space of $2$ points in the complex upper half-plane), needed to define the propagators in the case of no branes and one brane, see~\cite{K,CF}, and of the I-cube $\mathcal C_{2,1}^+$, needed to define the $4$-colored propagators, see~\cite{CFFR}.
Without going here into further details, the manifold $\mathcal C_{2,0,0,0}^+(\sqsubset)$ is diffeomorphic, as a smooth manifold with corners, to $\mathcal C_{2,2}^+$: the latter space is a quotient with respect to rescalings and translations, which can be used to fix the two ordered vertices on the real axis to $\{0,1\}$, and then we may in a conformal way map the complex upper half-plane to the half-strip $\sqsubset$, the point $1$ to $-\mathrm i\frac{\pi}2$, $0$ to $\frac{\pi}2$, the half-line right to $1$ to $\sqsubseteq_1$, the segment $[0,1]$ to $\sqsubseteq_2$ and the negative real axis to $\sqsubseteq_3$. 
This has been done explicitly in~\cite{F}, where propagators for the Poisson $\sigma$-model in presence of many branes have been discussed, and was sketched in the seminal paper~\cite{CFb}.

The construction of the $8$-colored propagators has been split in two pieces, namely, the case $j_1\neq j_3$ and $j_1=j_3$.

\subsubsection{The case $j_1\neq j_3$}\label{sss-6-1-1}
There are $4$ propagators which fall into this class, namely $\theta_{112}$, $\theta_{122}$, $\theta_{211}$ and $\theta_{221}$, which are constructed starting from the normalized, closed $1$-form on $\mathcal C_2\cong S^1$, the compactified configuration space of $2$ points in $\mathbb C$ by using the reflections with respect to $\sqsubseteq_i$, $i=1,2,3$.

We consider the boundary stratum of codimension $1$ of $\mathcal C_{2,0,0,0}^+(\sqsubset)$ corresponding to the collapse of the two points in the interior of the half-strip $\sqsubset$: more precisely, such a stratum is $\mathcal C_2\times \mathcal C_{1,0,0,0}^+(\sqsubset)\cong \mathcal C_2\times \mathcal C_{1,2}^+$, and local coordinates near such a stratum are given by
\[
\mathcal C_2\times \mathcal C_{1,0,0,0}^+(\sqsubset)\cong S^1\times C_{1,0,0,0}^+(\sqsubset)\ni (\varphi, z)\mapsto (z,z+\varepsilon \mathrm e^{\mathrm i\varphi})\in \mathcal C_{2,0,0,0}^+(\sqsubset),
\]
where the stratum is recovered as $\varepsilon$ tends to $0$.

An easy computation unraveling the formula for $\theta_{j_1,j_2,j_3}$, $j_1\neq j_3$, in~\cite[Subsubsection 6.2.1]{CT} in the same spirit of the proof of~\cite[Lemma 5.4]{CFFR}, yields the following behavior of $\theta_{112}$, $\theta_{122}$, $\theta_{211}$ and $\theta_{221}$, when restricted to the boundary stratum $\mathcal C_2\times \mathcal C_{1,0,0,0}^+(\sqsubset)$ of $\mathcal C^+_{2,0,0,0}(\sqsubset)$:
\[
\begin{aligned}
\theta_{112}\vert_{\mathcal C_2\times \mathcal C_{1,0,0,0}^+(\sqsubset)}&=\pi_1^*(\omega)+\pi_2^*(\widehat\rho), & \theta_{122}\vert_{\mathcal C_2\times \mathcal C_{1,0,0,0}^+(\sqsubset)}&=\pi_1^*(\omega)-\pi_2^*(\widehat\rho),\\
\theta_{221}\vert_{\mathcal C_2\times \mathcal C_{1,0,0,0}^+(\sqsubset)}&=\pi_1^*(\omega)+\pi_2^*(\widehat\rho), & \theta_{211}\vert_{\mathcal C_2\times \mathcal C_{1,0,0,0}^+(\sqsubset)}&=\pi_1^*(\omega)-\pi_2^*(\widehat\rho),\\
\end{aligned}
\]
where $\omega$ is the normalized volume form of $\mathcal C_2\cong S^1$ (the ``singular part'' of the propagator), and $\widehat\rho$ is the smooth, closed $1$-form on $\mathcal C_{1,0,0,0}^+(\sqsubset)$ given by the formula
\begin{equation}
\widehat\rho(z)=\frac{1}{2\pi}\left[\mathrm d\ \mathrm{arg}\!\left(\mathrm i\frac{\pi}2+z\right)-\mathrm d\ \mathrm{arg}\!\left(-\mathrm i\frac{\pi}2+z\right)-\mathrm d\ \mathrm{arg}\!\left(\mathrm i\pi+\mathrm{Re}(z)\right)\right],
\end{equation}
the ``regular part''.
Finally, $\pi_i$, $i=1,2$, is the projection from $\mathcal C_2\times \mathcal C_{1,0,0,0}^+(\sqsubset)$ onto the $i$-th factor.

The configuration space $\mathcal C_{1,0,0,0}^+(\sqsubset)$ is a smooth manifold with corners of dimension $2$: it has six boundary strata of codimension $1$, which correspond to the collapse of the point in the interior of the half-strip $\sqsubset$ to the boundary components $\sqsubseteq_i$ or to the two angle points $\pm\mathrm i\frac{\pi}2$, and the boundary stratum ``at infinity'', where the point in $\sqsubset$ tends to $\infty$ in $\sqsubset$ ({\em i.e.} along a horizontal line in $\sqsubset$).
In the first three cases, the boundary strata are simply $\sqsubseteq_i$, while on the remaining two, they are $\mathcal C_{1,1}^+\times \left\{\pm\mathrm i\frac{\pi}2\right\}$: again, using local coordinates near such boundary strata, we may prove that $\widehat\rho$ vanishes on the first three boundary strata, while on the remaining two we have
\[
\widehat\rho\vert_{\mathcal C_{1,1}^+\times \left\{\mathrm i\frac{\pi}2\right\}}=-\rho,\quad \widehat\rho\vert_{\mathcal C_{1,1}^+\times \left\{-\mathrm i\frac{\pi}2\right\}}=\rho,
\]
where $\rho$ is the (normalized) angle form on $\mathcal C_{1,1}^+$ (or, in previous terminology, the short loop contribution).
Finally, $\widehat\rho$ vanishes on the boundary stratum ``at infinity''.

\subsubsection{The case $j_1=j_3$}\label{sss-6-1-2}
The construction of the propagators $\theta_{111}$, $\theta_{121}$, $\theta_{212}$ and $\theta_{222}$, on the other hand, relies on a trickier argument, which we now review in more details.

We consider the strip $S=\left\{z\in\mathbb C:\ -\frac{\pi}2\leq \mathrm{Im}(z)\leq \frac{\pi}2\right\}$; on its interior, we consider the metric $g=\frac{\mathrm d x^2+\mathrm d y^2}{\cos^2 y}$, where $z=x+\mathrm iy$.
It is not difficult to prove, by a direct computation, that $g$ tends to the standard Poincar\'e hyperbolic metric, when we approach the two boundary lines of $S$: the basic propagator in deformation quantization (Kontsevich's angle form) is constructed {\em via} hyperbolic geometry, and the main idea behind the construction of $\theta_{j_1,j_2,j_1}$ is to use the geometry of $S$ determined by the metric $g$.

A slight variation of the computations leading to the general form of geodesics in the complex upper half-plane $\mathbb H^+$, endowed with the hyperbolic Poincar\'e metric, leads to the following general form of geodesics in the interior of $S$, endowed with the metric $g$:  
\begin{equation}\label{eq-geod}
\sin(y)=A\mathrm e^x+B\mathrm e^{-x},
\end{equation}
where $A$, $B$ are real constants; we additionally have geodesic vertical segments in the interior of $S$, see also~\cite[Figure 6]{CT} for a pictorial description of the geodesics in the interior of $S$ with respect to metric $g$.

It follows immediately from~\eqref{eq-geod} that, for any two points $z_1$, $z_2$ in $S$, there is a unique geodesic passing through them.
%When $\mathrm{Re}(z_1)=\mathrm{Re}(z_2)$, we consider the unique vertical segment joining $z_1$ and $z_2$, while, if $\mathrm{Re}(z_1)\neq \mathrm{Re}(z_2)$, the corresponding linear system for $A$, $B$ admits a unique solution.

Thus, it makes sense to define the geodetic angle function $\widetilde\vartheta$ on the open configuration space $C^+_{2,0,0}(S)$ of two points in the interior of $S$ and no point on the two boundary lines in a way similar to Kontsevich's angle function, {\em i.e.} $\widetilde\vartheta(z_1,z_2)$, for $z_1\neq z_2$ in the interior of $S$, is the angle between the geodesic vertical segment going through $z_1$ and the unique geodesic joining $z_1$ and $z_2$ (obviously, $\widetilde\vartheta$ is well-defined up to the addition of integer multiples of $2\pi$).

More explicitly, the geodetic angle function $\vartheta(z_1,z_2)$ is defined through
\[
\tan\!\left(\widetilde\vartheta(z_1,z_2)+\frac{\pi}2\right)=\frac{\sin(y_1)\cosh(x_1-x_2)-\sin(y_2)}{\cos(y_1)\sinh(x_1-x_2)},\quad z_i=x_i+\mathrm i y_i,\ i=1,2.
\]
It is easily verified that $\widetilde\vartheta(z_1,z_2)$ extends to $S$.

The compactified configuration space $\mathcal C^+_{2,0,0}(S)$ has a boundary stratification: we are particularly interested into the boundary strata of codimension $1$, which correspond to the collapse of exactly one point in the interior of $S$ to one of the two boundary lines of $S$, of both points together in the interior of $S$ and of both points together to a point on one of the two boundary lines of $S$.
The stratum corresponding to the collapse of exactly one point in the interior of $S$ to one of the two boundary axis is represented either by $\mathcal C_{1,0}^+\times\mathcal C_{1,1,0}^+(S)$ or $\mathcal C_{1,0}^+\times\mathcal C_{1,0,1}^+(S)$; the stratum corresponding to the collapse of the two points in the interior of $S$ by $\mathcal C_2\times\mathcal C_{1,0,0}^+(S)$, and the remaining two boundary strata are represented either $\mathcal C_{2,0}^+\times\mathcal C_{0,1,0}^+(S)$ or $\mathcal C_{2,0}^+\times\mathcal C_{0,0,1}^+(S)$.

Using local coordinates near the boundary strata of codimension $1$ previously analyzed, we can prove that the exterior derivative of $\widetilde\vartheta$, which we denote as $\vartheta(z_1,z_2)$ as in~\cite[Subsubsection 6.2.1]{CT}, is a well-defined closed $1$-form on the compactified configuration space $\mathcal C_{2,0,0,0}^+(S)$, which satisfies the following properties:
\begin{enumerate}
\item[$i)$] the restriction of $\vartheta$ to the boundary stratum $\mathcal C_2\times\mathcal C_{1,0,0}^+(S)$ equals $\pi_1^*(\omega)$, $\pi_i$ being the natural projection onto the $i$-th factor, and $\omega$ the normalized volume form on $\mathcal C_2\cong S^1$.
\item[$ii)$] The restriction of $\vartheta$ to either one of the boundary strata $\mathcal C_{1,0}^+\times \mathcal C_{1,1,0}^+(S)$ or $\mathcal C_{1,0}^+\times \mathcal C_{1,0,1}^+(S)$ corresponding to the collapse of the first argument to either one of the two boundary lines of $S$ vanishes.
\item[$iii)$] The restriction of $\vartheta$ to either one of the boundary strata $\mathcal C_{2,0}^+\times\mathcal C_{0,1,0}^+(S)$ or $\mathcal C_{2,0}^+\times\mathcal C_{0,0,1}^+(S)$ equals $\pi_1^*(\omega^\pm)$, where $\pi_i$ is, once again, the projection onto the $i$-th factor, and $\omega^\pm$ denotes the $2$-colored propagators (Kontsevich's angle function and its image with respect to the involution exchanging the two arguments).
\end{enumerate}

We now borrow from~\cite[Subsubsection 6.2.1]{CT} the definition of the propagators $\theta_{1,j_2,1}$ and $\theta_{2,j_2,2}$, $j_2=1,2$, which are smooth, closed $1$-forms on the compactified configuration space $\mathcal C_{2,0,0}^+(\sqsubset)$: they are constructed using the angle form $\vartheta$ on $\mathcal C_{2,0,0}^+(S)$ and the natural involution $\sigma$ of $S$ associated to the reflection with respect to the imaginary axis, see~\cite[Subsubsection 6.2.1]{CT} for the explicit formul\ae.
We observe that the angle form $\vartheta$ is equivariant with respect to the $\mathbb Z_2$-action induced on $\mathcal C_{2,0,0}^+(S)$ by the diagonal action of $\sigma$ and the natural sign action.

Once again, we are interested to the boundary stratum $\mathcal C_2\times \mathcal C_{1,0,0}^+(\sqsubset)$ of $\mathcal C_{2,0,0}^+(\sqsubset)$, where the points in $\sqsubset$ collapse together in $\sqsubset$.
As an example we consider the case $j_1=j_3=1$: then, we have the explicit formul\ae
\[
\theta_{1,1,1}(z_1,z_2)=\frac{1}{2\pi}\left[\vartheta(z_1,z_2)-\vartheta(\sigma(z_1),z_2)\right],\quad \theta_{1,2,1}(z_1,z_2)=\frac{1}{2\pi}\left[\vartheta(z_1,z_2)-\vartheta(z_1,\sigma(z_2))\right],\quad (z_1,z_2)\in\mathcal C_{2,0,0}^+(\sqsubset).
\]

We may use the local coordinates of Subsubsection~\ref{sss-6-1-1} near the said boundary stratum to perform explicit computations: using the properties of the angle form $\vartheta$, we see that the restriction of $\theta_{1,1,1}$ and $\vartheta_{1,2,1}$ to $\mathcal C_2\times \mathcal C_{1,0,0}^+(\sqsubset)$ splits into a singular part and a regular part; more precisely
\[
\theta_{1,1,1}\vert_{\mathcal C_2\times \mathcal C_{1,0,0}^+(\sqsubset)}=\pi_1^*(\omega)+\pi_2^*(\widetilde\rho),\quad \theta_{1,1,1}\vert_{\mathcal C_2\times \mathcal C_{1,0,0}^+(\sqsubset)}=\pi_1^*(\omega)-\pi_2^*(\widetilde\rho),
\]
where $\pi_i$, $i=1,2$, and $\omega$ are as before.
On the other hand, the regular term $\widetilde\rho$ is a smooth, closed $1$-form on $\mathcal C_{1,0,0}^+(\sqsubset)$ defined {\em via}
\[
\widetilde\rho(z)=\frac{1}{2\pi}\vartheta(z,\sigma(z))=\frac{1}{2\pi}\mathrm d\arctan\!\left(\tanh(x)\tan(y)\right),\quad z=x+\mathrm i y\in\mathcal C_{1,0,0}^+(\sqsubset).
\]
Direct computations show that $\widetilde\rho$ has a behavior similar to the regular term of Subsubsection~\ref{sss-6-1-1} on the boundary strata of $\mathcal C_{1,0,0}^+(\sqsubset)$, with the only difference that $\widetilde\rho$ is non-trivial, when restricted to the boundary stratum ``at infinity'' (which is equivalent to a closed interval; in more familiar terms, it is $\mathcal C_{1,1}^+$): in fact, it equals $\rho$, the short loop contribution.

\subsection{Final considerations on polarizations}\label{ss-6-2}
Having constructed the $8$-colored propagators in detail, we now want to discuss, in light of the appearance of ``regular terms'' also in the $8$-colored propagators, how such changes affect the arguments and the computations in the final stages of~\cite{CT}.

The main direct consequence of the construction of the $8$-colored propagators is explained in~\cite[Proposition 21, Subsubsection 6.2.2]{CT}: therein, for the case of a symmetric pair $\mathfrak g=\mathfrak k\oplus \mathfrak p$, an element $f$ of $\mathfrak k^\perp$, and two polarizations $\mathfrak b_i$, $i=1,2$ (we observe that $\mathfrak k$ and $\mathfrak b_i$, $i=1,2$, are assumed to be in a position of normal intersection, {\em i.e.} $\mathfrak k\cap(\mathfrak b_1+\mathfrak b_2)=\mathfrak k\cap\mathfrak b_1+\mathfrak k\cap\mathfrak b_2$), it has been proved that the corresponding characters, depending upon the choice of either $\mathfrak b_1$ or $\mathfrak b_2$, are equal.

The key argument of the proof relies, as the principle of biquantization, on Stokes' Theorem applied to the situation, where we consider sums over admissible graphs in $\mathcal G_{n,3}$ and corresponding integral weights: here we view an admissible graph $\Gamma$ in $\mathcal G_{n,3}$ as an embedded graph in $\sqsubset$, where the first, resp.\ the third, vertex of the second type is $\mathrm i\frac{\pi}2$, resp.\ $-\mathrm i\frac{\pi}2$ (consequently, the second vertex lies on the vertical boundary segment $\sqsubseteq_2$ of $\sqsubset$), graphically
\bigskip
\begin{center}
\resizebox{0.7 \textwidth}{!}{\input{tri_quant_graph.pstex_t}}\\
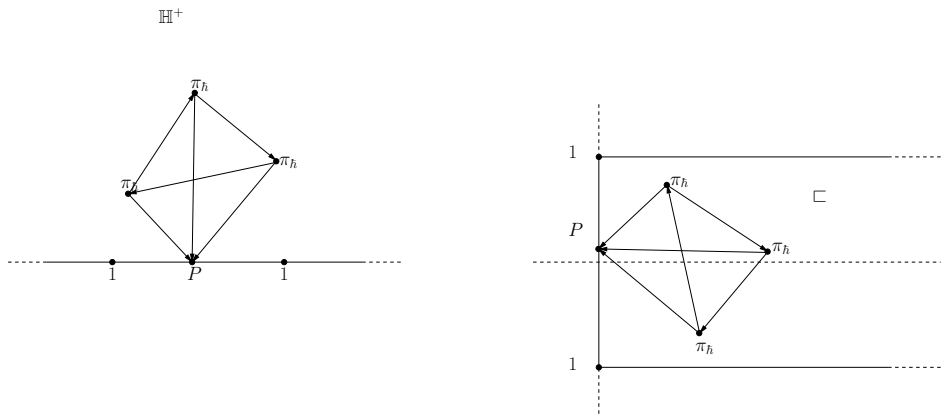
\text{Figure 14 - An admissible graph of type $(3,3)$ in $\mathbb H^+\sqcup \mathbb R$ and in $\sqsubset$} \\
\end{center}
\bigskip
In the previous picture, $P$ is a general element of $\mathrm H^0_{\hbar,\mathfrak b_1,\mathfrak b_2}(\mathfrak k^\perp)$, using notations from~\cite[Subsubsection 6.2.2]{CT}.

As has been already observed 
%in Section~\ref{s-3} 
about biquantization, we need to consider admissible graphs with short loop contributions on vertices of the first type, if we want Stokes' Theorem to do the job, because of the ``regular term'': similarly, the results of Subsubsections~\ref{sss-6-1-1} and~\ref{sss-6-1-2} imply that what may be called ``triquantization'' can be performed using the techniques of Deformation Quantization with some changes, which should keep into account the ``regular terms'' in the $8$-colored propagators (which, by the way, are completely consequent with the ``regular term'' in the $4$-colored propagators).

The first obvious observation about admissible graphs $\Gamma$ in $\mathcal G_{n,3}$ is that between any two vertices there can be {\em at most} eight edges, because to any edge we assign one of the $8$-colored propagators; therefore, we consider admissible graphs with multiple edges.
From the point of view of combinatorics, we will have to take into account in the integral weight the (possible) multiplicity of edges of an admissible graph $\Gamma$ of $\mathcal G_{n,3}$.

To $\Gamma$, we may associate an integral weight by standard prescriptions as in~\cite{K,CF,CFFR}, using the $8$-colored propagators: in this particular situation, as we want to apply Stokes' Theorem, the integral weight $w_\Gamma$ of a general admissible graph $\Gamma$ in $\mathcal G_{n,3}$, for $n\geq 0$, is defined as an integral of the exterior derivative of a product $\omega_\Gamma$ of $8$-colored propagators (specified by the shape of $\Gamma$) over the compactified configuration space $\mathcal C_{n,0,1,0}^+(\sqsubset)$, which is orientable and of dimension $2n+1$.

Such integrals are, on the one hand, obviously trivial (because the $8$-colored propagators are closed); on the other hand, if the degree of the integrand is $2n+1$ ({\em i.e.} if the degree of $\omega_\Gamma$ is $2n$), such a trivial contribution equals the sum over all boundary strata of codimension $1$ of $\mathcal C_{n,0,1,0}^+(\sqsubset)$: of course, this argument is similar to the proof {\em e.g.} of associativity of the $\star$-product in~\cite{K}.

More precisely, we have
\[
0=\int_{\mathcal C_{n,0,1,0}^+(\sqsubset)}\mathrm d\omega_\Gamma=\sum_{i}\pm\int_{\partial_i\mathcal C_{n,0,1,0}^+(\sqsubset)}\omega_\Gamma,
\]
where the sum is over all boundary strata of codimension $1$ of $\mathcal C_{n,0,1,0}^+(\sqsubset)$, and the signs are dictated by their orientations.
\begin{Rem}\label{r-triquant}
We are purposefully sketchy here, but we plan to return to these issues somewhere else, as the $8$-colored propagators are the central tool in the construction of a ``Formality Theorem in presence of $3$ branes'' generalizing the main result of~\cite{CFFR}; triquantizazion should be then related to the evaluation of the corresponding formality quasi-isomorphism at a Poisson structure on some linear space $X$.
\end{Rem}

Let us consider only the boundary strata of codimension $1$ of $\mathcal C_{n,0,1,0}^+(\sqsubset)$ corresponding to the collapse of a subset $A$ of $\{1,\dots,n\}$ of cardinality $2\leq |A|\leq n$ in the interior of $\sqsubset$: such boundary strata are simply $\mathcal C_A\times \mathcal C_{n-|A|+1,0,1,0}^+(\sqsubset)$, where $\mathcal C_A$ is the compactified configuration space of $|A|$ points in $\mathbb C$ (modulo rescalings and complex translations, see~\cite{K}), which is an orientable manifold with corners of dimension $2|A|-3$.
We denote by $\Gamma_A$, resp.\ $\Gamma^A$, the subgraph of $\Gamma$, whose vertices are labeled by $A$ and whose edges are all edges connecting vertices labeled by $A$, resp.\ the graph obtained by contracting $\Gamma_A$ to a single vertex (necessarily of the first type).

The restriction of $\Gamma_A$ to $\mathcal C_A\times \mathcal C_{n-|A|+1,0,1,0}^+(\sqsubset)$ is a product of forms splitting into the sum of a ``singular term'' (living on $\mathcal C_A$) and of ``regular terms'' $\widehat\rho$ or $\widetilde\rho$ (living, on the other hand, on $\mathcal C_{n-|A|+1,0,1,0}^+(\sqsubset)$).
Recalling~\cite[Lemma 6.6]{K} and the fact that only the ``singular part'' of the integrand ({\em i.e.} the product of all ``singular terms'' in $\omega_\Gamma$) is to be integrated over $\mathcal C_A$, we reduce to the case $|A|=2$.
Further, the ``singular term'' and both ``regular terms'' are all $1$-forms: in particular, the shape of $\omega_\Gamma$ forces that possible non-trivial factors in the restriction of $\omega_\Gamma$ on such a boundary stratum are associated to two vertices of the first type, labeled by $A$, which are joined by {\em at most} three edges (either multiple edges or not).
\bigskip
\begin{center}
\resizebox{0.7 \textwidth}{!}{\input{triq_reg-sing.pstex_t}}\\
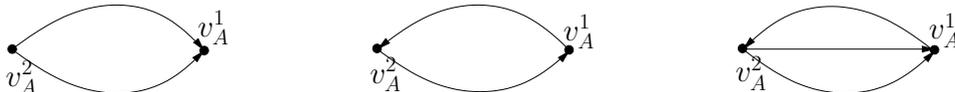
\text{Figure 15 - Possible non-trivial contributions from the stratum $\mathcal C_2\times \mathcal C^+_{n-1,0,1,0}(\sqsubset)$} \\
\end{center}
\bigskip
If there are only two edges between $v_A^1$ and $v_A^2$, then the weight contribution after integrating over $\mathcal C_2\cong S^1$ can be a sum of $\widehat\rho$ and $\widetilde\rho$ of the form $\pm a\widehat\rho\pm b\widetilde\rho$, $a$, $b$ in $\{0,1\}$.
On the other hand, if there are three edges between $v_A^1$ and $v_A^2$, the weight contribution after integration is the product of $\widehat\rho$ and $\widetilde\rho$.

Further, to $\Gamma$ in $\mathcal G_{n,3}$ we associate a polydifferential operator, which acts on the ($\hbar$-shifted) Poisson bivector $\pi_\hbar$ (a copy of which is placed at an vertex of the first type) and to the triple $(1|P|1)$, where $1$ is regarded as a constant function and $P$ is as above.
The rule associates to each oriented edge of $\Gamma$ a derivation operator (in the graded sense, as we deal here with graded vector spaces).
 Finally, we multiply the end result as an element of $\mathfrak g^*$, and take its restriction to $f+(\mathfrak k^\perp+\mathfrak b_1+\mathfrak b_2)$.

We now recall that $X=\mathfrak g^*$ is endowed with a ($\hbar$-shifted) linear Poisson structure, whence no multiple edges are possible between vertices of the first type.
As a consequence, every vertex of the first type admits {\em at most} one incoming edge, {\em i.e.} $\Gamma_A$ can be only of the second type in Figure 15; any of the vertices of $\Gamma_A$ may further have {\em at most} one outgoing edge.

We briefly discuss the coloring of an admissible graph.
We choose a system of coordinates on $\mathfrak g$ which is adapted to the coisotropic submanifolds $\mathfrak k^\perp$, $f+\mathfrak b_1$ and $f+\mathfrak b_2$: in other words, we assume there is a partition of $\left\{1,\dots,d\right\}$, where $d$ is the dimension of $\mathfrak g$ of the form
\[
\begin{aligned}
\{1,\dots,d\}=&(I_1\cap I_2\cap I_3)\sqcup(I_1^c\cap I_2\cap I_3)\sqcup(I_1\cap I_2^c\cap I_3)\sqcup(I_1\cap I_2\cap I_3^c)\sqcup\\
&(I_1^c\cap I_2^c\cap I_3)\sqcup(I_1^c\cap I_2\cap I_3^c)\sqcup(I_1\cap I_2^c\cap I_3^c)\sqcup(I_1^c\cap I_2^c\cap I_3^c),
\end{aligned}
\] 
labeling a basis of $\mathfrak g$.
This means that {\em e.g.} the elements of the basis indexed by $i$ in $I_1\cap I_2\cap I_3$ constitute a coordinate system for the intersection of $\mathfrak k^\perp$, $f+\mathfrak b_1$ and $f+\mathfrak b_2$, the elements indexed by $i$ in $I_1^c\cap I_2\cap I_3$ a coordinate system for the intersection of $\mathfrak k^*$ with $f+\mathfrak b_1$ and $f+\mathfrak b_2$, {\em et similiter}.

The choice of labeling the coordinates on $\mathfrak g^*$ with respect to the above partition yields an obvious coloring of an admissible graph $\Gamma$ in $\mathcal G_{n,3}$: in fact, to any edge we may associate a triple $(j_1,j_2,j_3)$, $j_k$ in $\{1,2\}$, by the rule that the ``parity'' of $I_k$, resp.\ $I_k^c$, is $1$, resp.\ $2$, $k=1,2,3$.
This rule simultaneously determines, for any colored edge of $\Gamma$, the coordinate set with respect to which the edge takes derivation, and the labeling of the edge by one of the $8$-colored propagators.

Therefore, adjusting the arguments of the discussion at the beginning of~\cite[Subsection 7.1]{CFFR}, and taking into account the polydifferential operator associated to the only possibly non-trivial subgraph $\Gamma_A$ as in Figure 15, we see that these problems may be corrected by allowing for the presence of admissible graphs with short loops: the main difference between the results of~\cite{CFFR} and the triquantization discussed here is the fact that triquantization requires the presence of two distinct short loops, namely, one taking care of the ``regular term'' $\widehat\rho$ and one of the ``regular term'' $\widetilde\rho$.

Of course, there is a geometric counterpart to regular terms (which we discuss here very briefly, hoping to return to a more precise statement in the context of a formality theorem in presence of $3$ branes): as in~\cite{CFFR}, the geometric counterpart is played by a partial divergence operator, whose main property is Leibniz' rule with respect to the Schouten--Nijenhuis bracket, which is essential in the computations.

Since both ``regular terms'' are closed $1$-forms on $\mathcal C_{1,0,0,0}^+(\sqsubset)$, an admissible graph $\Gamma$ may admit {\em at most} two short loops of different type at each vertex of the first type: in fact, a vertex of the first type admits in this situation {\em at most} one short loop contribution, because of the linearity of the Poisson structure and because a divergence operator is of order $1$.
In particular, a vertex of the first type with a short loop (either $\widehat\rho$ or $\widetilde\rho$) must have an outgoing edge, otherwise dimensional argument imply triviality of the corresponding integral weight. 

Previous computations imply that both short loop contributions vanish, when restricted to boundary strata $\mathcal C_{A,1}^+\times\mathcal C_{n-|A|,0,1,0}^+(\sqsubset)$, for subset $A$ of $\{1,\dots,n\}$ of cardinality $1\leq |A|\leq n$ of $\mathcal C_{n,0,1,0}^+(\sqsubset)$.
Furthermore, when we consider restrictions to boundary strata of the form $\mathcal C_{A,2}^+\times \mathcal C_{n-|A|,0,0,0}^+(\sqsubset)$, the short loop contributions $\widehat\rho$ and $\widetilde\rho$ restrict to the short loop contribution in biquantization.

Pictorially we now have to consider admissible graphs of the following shape:
\bigskip
\begin{center}
\resizebox{0.5 \textwidth}{!}{\input{triquant_adm.pstex_t}}\\
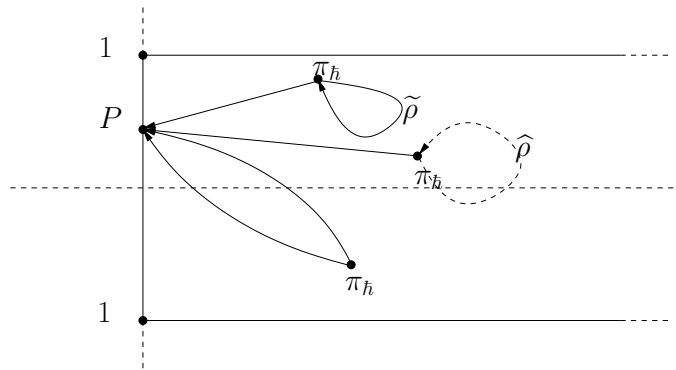
\text{Figure 16 - An admissible graph $\Gamma$ in $\mathcal G_{n,3}$ with multiple edges and short loops in $\subseteq$} \\
\end{center}
\bigskip

Summarizing the discussion so far, Proposition 21 in~\cite[Subsubsection 6.2.2]{CT} remains valid, provided we enlarge the set $\mathcal G_{n,3}$ of admissible graphs, for $n\geq 0$, so as to contain graphs with multiple edges and short loop contributions of two distinct types at vertices of the first type.

Finally, the results of~\cite[Subsubsections 6.3.2 and 6.3.3]{CT} are easily proved to be still valid with the previous modification of admissible graphs.


\begin{bibdiv}
\begin{biblist}

\bib{AM}{article}{
   author={Alekseev, A.},
   author={Meinrenken, E.},
   title={On the Kashiwara-Vergne conjecture},
   journal={Invent. Math.},
   volume={164},
   date={2006},
   number={3},
   pages={615--634},
   issn={0020-9910},
   review={\MR{2221133 (2007g:17017)}},
   doi={10.1007/s00222-005-0486-4},
}

\bib{AS}{article}{
   author={Axelrod, Scott},
   author={Singer, I. M.},
   title={Chern-Simons perturbation theory},
   conference={
      title={ Geometric Methods in Theoretical Physics, Vol.\ 1, 2},
      address={New York},
      date={1991},
   },
   book={
      publisher={World Sci. Publ., River Edge, NJ},
   },
   date={1992},
   pages={3--45},
   review={\MR{1225107 (94g:58244)}},
}

\bib{BCKT}{article}{
   author={Cattaneo, Alberto},
   author={Keller, Bernhard},
   author={Torossian, Charles},
   author={Brugui{\`e}res, Alain},
   title={Introduction},
   language={French},
   note={Dual French-English text},
   conference={
      title={D\'eformation, quantification, th\'eorie de Lie},
   },
   book={
      series={Panor. Synth\`eses},
      volume={20},
      publisher={Soc. Math. France},
      place={Paris},
   },
   date={2005},
   pages={1--9, 11--18},
   review={\MR{2274223}},
}


\bib{CFFR}{article}{
  author={Calaque, Damien},
    author={Felder, Giovanni},
      author={Ferrario, Andrea},
        author={Rossi, Carlo A.},
  title={Bimodules and branes in deformation quantization},
  eprint={arXiv:0908.2299},
  journal={Comp.\ Math.\ (to appear)}
  date={2009}
}

\bib{CFb}{article}{
    AUTHOR = {Cattaneo, Alberto S. and Felder, Giovanni},
     TITLE = {Coisotropic submanifolds in {P}oisson geometry and branes in
              the {P}oisson sigma model},
   JOURNAL = {Lett. Math. Phys.},
  FJOURNAL = {Letters in Mathematical Physics. A Journal for the Rapid
              Dissemination of Short Contributions in the Field of
              Mathematical Physics},
    VOLUME = {69},
      YEAR = {2004},
     PAGES = {157--175},
      ISSN = {0377-9017},
     CODEN = {LMPHDY},
   MRCLASS = {81T45 (22A22 53D17 53D20 53D55)},
  MRNUMBER = {MR2104442 (2005m:81285)},
MRREVIEWER = {Stefan Waldmann},
       DOI = {10.1007/s11005-004-0609-7},
       URL = {http://dx.doi.org/10.1007/s11005-004-0609-7},
}
		
\bib{CF}{article}{
   author={Cattaneo, Alberto S.},
   author={Felder, Giovanni},
   title={Relative formality theorem and quantisation of coisotropic
   submanifolds},
   journal={Adv. Math.},
   volume={208},
   date={2007},
   number={2},
   pages={521--548},
   issn={0001-8708},
   review={\MR{2304327 (2008b:53119)}},
   doi={10.1016/j.aim.2006.03.010},
}

\bib{CT}{article}{
   author={Cattaneo, Alberto S.},
   author={Torossian, Charles},
   title={Quantification pour les paires sym\'etriques et diagrammes de
   Kontsevich},
   language={French, with English and French summaries},
   journal={Ann. Sci. \'Ec. Norm. Sup\'er. (4)},
   volume={41},
   date={2008},
   number={5},
   pages={789--854},
   issn={0012-9593},
   review={\MR{2504434 (2010g:22031)}},
}

\bib{Dix}{book}{
   author={Dixmier, Jacques},
   title={Enveloping algebras},
   series={Graduate Studies in Mathematics},
   volume={11},
   note={Revised reprint of the 1977 translation},
   publisher={American Mathematical Society},
   place={Providence, RI},
   date={1996},
   pages={xx+379},
   isbn={0-8218-0560-6},
   review={\MR{1393197 (97c:17010)}},
}

\bib{F}{article}{
   author={Ferrario, Andrea},
   title={Poisson sigma model with branes and hyperelliptic Riemann
   surfaces},
   journal={J. Math. Phys.},
   volume={49},
   date={2008},
   number={9},
   pages={092301, 23},
   issn={0022-2488},
   review={\MR{2455835 (2009j:81140)}},
   doi={10.1063/1.2982234},
}

\bib{FMcP}{article}{
   author={Fulton, William},
   author={MacPherson, Robert},
   title={A compactification of configuration spaces},
   journal={Ann. of Math. (2)},
   volume={139},
   date={1994},
   number={1},
   pages={183--225},
   issn={0003-486X},
   review={\MR{1259368 (95j:14002)}},
   doi={10.2307/2946631},
}

%\bib{FRW}{article}{
%  author={Ferrario, Andrea},
%    author={Rossi, Carlo A.},
%      author={Willwacher, Thomas},
%  title={A note on the Koszul complex in deformation quantization},
%  eprint={arXiv:1002.2561},
%  date={2010}
%}

\bib{Kath}{article}{
   author={Kathotia, Vinay},
   title={Kontsevich's universal formula for deformation quantization and
   the Campbell-Baker-Hausdorff formula},
   journal={Internat. J. Math.},
   volume={11},
   date={2000},
   number={4},
   pages={523--551},
   issn={0129-167X},
   review={\MR{1768172 (2002h:53154)}},
   doi={10.1142/S0129167X0000026X},
}


\bib{Kel}{article}{
   author={Keller, Bernhard},
   title={$A$-infinity algebras, modules and functor categories},
   conference={
      title={Trends in representation theory of algebras and related topics},
   },
   book={
      series={Contemp. Math.},
      volume={406},
      publisher={Amer. Math. Soc.},
      place={Providence, RI},
   },
   date={2006},
   pages={67--93},
   review={\MR{2258042 (2007g:18002)}},
}


\bib{K}{article}{
   author={Kontsevich, Maxim},
   title={Deformation quantization of Poisson manifolds},
   journal={Lett. Math. Phys.},
   volume={66},
   date={2003},
   number={3},
   pages={157--216},
   issn={0377-9017},
   review={\MR{2062626 (2005i:53122)}},
   doi={10.1023/B:MATH.0000027508.00421.bf},
}

\bib{Lef-Has}{article}{
  author={Lef\`evre-Hasegawa, Kenji},
  title={Sur les $A_\infty$-cat\'egories},
  eprint={http://people.math.jussieu.fr/~keller/lefevre/TheseFinale/tel-00007761.pdf},
  date={2003}
}

\bib{Rouv}{article}{
   author={Rouvi{\`e}re, Fran{\c{c}}ois},
   title={Espaces sym\'etriques et m\'ethode de Kashiwara-Vergne},
   language={French, with English summary},
   journal={Ann. Sci. \'Ecole Norm. Sup. (4)},
   volume={19},
   date={1986},
   number={4},
   pages={553--581},
   issn={0012-9593},
   review={\MR{875088 (88d:43010)}},
}


\bib{Sh}{article}{
   author={Shoikhet, Boris},
   title={Vanishing of the Kontsevich integrals of the wheels},
   note={EuroConf\'erence Mosh\'e Flato 2000, Part II (Dijon)},
   journal={Lett. Math. Phys.},
   volume={56},
   date={2001},
   number={2},
   pages={141--149},
   issn={0377-9017},
   review={\MR{1854132 (2002j:53119)}},
   doi={10.1023/A:1010842705836},
}


\bib{T1}{article}{
   author={Torossian, Charles},
   title={Op\'erateurs diff\'erentiels invariants sur les espaces
   sym\'etriques. I. M\'ethodes des orbites},
   language={French, with English and French summaries},
   journal={J. Funct. Anal.},
   volume={117},
   date={1993},
   number={1},
   pages={118--173},
   issn={0022-1236},
   review={\MR{1240263 (95c:22016a)}},
   doi={10.1006/jfan.1993.1124},
}

\bib{T2}{article}{
   author={Torossian, Charles},
   title={Op\'erateurs diff\'erentiels invariants sur les espaces
   sym\'etriques. II. L'homomorphisme d'Harish-Chandra g\'en\'eralis\'e},
   language={French, with English and French summaries},
   journal={J. Funct. Anal.},
   volume={117},
   date={1993},
   number={1},
   pages={174--214},
   issn={0022-1236},
   review={\MR{1240264 (95c:22016b)}},
   doi={10.1006/jfan.1993.1125},
}

\bib{T3}{article}{
   author={Torossian, Charles},
   title={La conjecture de Kashiwara-Vergne (d'apr\`es Alekseev et
   Meinrenken)},
   language={French, with French summary},
   note={S\'eminaire Bourbaki. Vol. 2006/2007},
   journal={Ast\'erisque},
   number={317},
   date={2008},
   pages={Exp. No. 980, ix, 441--465},
   issn={0303-1179},
   isbn={978-2-85629-253-2},
   review={\MR{2487742 (2010d:22011)}},
}


\bib{VdB}{article}{
   author={Van den Bergh, Michel},
   title={The Kontsevich weight of a wheel with spokes pointing outward},
   journal={Algebr. Represent. Theory},
   volume={12},
   date={2009},
   number={2-5},
   pages={443--479},
   issn={1386-923X},
   review={\MR{2501196 (2010e:53157)}},
   doi={10.1007/s10468-009-9161-6},
}


\end{biblist}
\end{bibdiv}
\end{document}